\documentclass[a4paper,12pt,leqno]{article}

\usepackage[latin1]{inputenc} 
\usepackage[T1]{fontenc}      
\usepackage{geometry}         
\usepackage[francais, english]{babel} 
\usepackage{amsmath}
\usepackage{amsfonts}
\usepackage{amssymb}
\usepackage{fancyhdr}
\usepackage{url}
\usepackage{mathrsfs}
\usepackage{euscript}
\usepackage{enumerate}
\usepackage{cite}
\usepackage{textcomp}
                              
\newtheorem{df}{Definition}[section]

\newtheorem{prop}[df]{Proposition}
\newtheorem{thm}[df]{Theorem}
\newtheorem{dfn}{Definition}
\newtheorem{thrm}[dfn]{Theorem}
\newtheorem{conj}[dfn]{Conjecture}
\newtheorem{lem}[df]{Lemma}

\newcommand{\prf}{\noindent \textit{Proof}}
\newtheorem{rmk}[df]{Remark}

\newcommand{\ddbar}{\partial\overline{\partial}}

\newcommand{\R}{\mathbb{R}}

\newcommand{\C}{\mathbb{C}}

\newcommand{\N}{\mathbb{N}}
\newcommand{\XD}{X\backslash D}
\newcommand{\YE}{Y\backslash E}
\newcommand{\UD}{U\backslash D}

\newcommand{\lich}{\mathbb{L}}

\newcommand{\vol}{\operatorname{vol}}

\newcommand{\vareps}{\varepsilon}

\newcommand{\vl}{\operatorname{Vol}}

\newcommand{\scal}{\mathbf{s}}

\newcommand{\loc}{\operatorname{loc}}
\newcommand{\cqfd}{ \hfill $\square$ }

\newcommand{\omegadelta}{\omega_{\Delta^*}}
\renewcommand{\Re}{\mathfrak{Re}}
\renewcommand{\Im}{\mathfrak{Im}}
\newcommand{\vthet}{ \vartheta }
\newcommand{\chern}{\operatorname{c}} 
\newcommand{\sbarD}{\overline{\scal}_D}
\newcommand{\sbarDj}{\overline{\scal}_{D_j}}
\newcommand{\sbar}{\overline{\scal}}
\newcommand{\spn}{\operatorname{span}}
\newcommand{\kphi}{\mathsf{K}_{\varphi}}

\title{\large ASYMPTOTIC PROPERTIES OF EXTREMAL KÄHLER METRICS OF POINCARÉ TYPE} 

\vspace{5pt}
\author{\normalsize\textsc{Hugues AUVRAY%
          \footnote{This work was started during the author's stay at the MPIM Bonn (EPDI post-doc, 2013), 
                    and completed at his arrival at ENS Cachan.}}}
\date{}

\begin{document}

\makeatletter
\renewcommand%
   {\section}%
   {%
   \@startsection{section}%
      {1}
      {0mm}%
      {\baselineskip}
      {0.5\baselineskip}%
      {\sc\large\centering}
   }%
\makeatother

\makeatletter
\renewcommand%
   {\subsubsection}%
   {%
   \@startsection{subsubsection}%
      {1}
      {0mm}%
      {1.25\baselineskip}
      {0.25\baselineskip}%
      {\sf\normalsize}
   }%
\makeatother

\makeatletter
\renewcommand%
   {\paragraph}%
   {%
   \@startsection{paragraph}%
      {4}%
      {0mm}%
      {0mm}%
      {-5pt}%
      {\it\normalsize}%
   }%
\makeatother

\maketitle

\renewcommand{\abstractname}{Abstract}
 \begin{abstract}
  Consider a compact Kähler manifold $X$ with a  simple normal crossing divisor $D$, 
  and define Poincaré type metrics on $\XD$ as Kähler metrics on $\XD$ with cusp singularities along $D$. 
  We prove that the existence of a constant scalar curvature (resp. an extremal)
  Poincaré type Kähler metric on $\XD$
  implies the existence of a constant scalar curvature (resp. an extremal) Kähler metric, 
  possibly of Poincaré type, on every component of $D$.  
  We also show that when the divisor is smooth, 
  the constant scalar curvature/extremal metric on $X\backslash D$ is asymptotically a product near the divisor. 
 \end{abstract}

 \begin{center}
   \rule{3cm}{0.25pt}
 \end{center}

~

\section*{Introduction}
 In his search for canonical representants of Kähler classes on compact Kähler manifolds,  
 generalising the Kähler-Einstein problem, 
 E. Calabi introduced \textit{extremal Kähler metrics}, 
 defined as the minimisers of the $L^2$-norm of the Ricci tensor among a fixed class \cite{cal}. 
 
 Extremal metrics turn out to satisfy rich geometric properties, e.g. 
 maximality of the group of isometric automorphisms among connected compact Lie groups of automorphisms \cite{cal2}. 
 Conversely though, these properties may be viewed as \textit{obstructions} to the existence of extremal metrics;  
 see for instance the example produced by M. Levine \cite{lev} of a complex Kähler surface admitting no extremal metric.  
 The subsequent (counter)examples produced by D. Burns and P. de Bartolomeis \cite{bdb} revealed moreover deeper links 
 between the (non-)existence of extremal metrics, and algebro-geometric conditions on the manifold.  
 
 In this direction, the so-called Yau-Tian-Donaldson conjecture  
 predicts that in the algebraic case, 
 the existence of extremal Kähler metrics is equivalent to a stability condition, 
 close to the Geometric Invariant Theory, on the polarised manifold: 
  \begin{conj}    \label{conj_YTD}
   Let $(X,L)$ be a compact polarised manifold. 
   Then there exists an extremal Kähler metric in $\chern_1(L)$ if and only if $(X,L)$ 
   is K-stable relatively to a maximal torus of ${\rm Aut}^0(X,L)$. 
  \end{conj} 
 This conjecture, first designed for Kähler-Einstein metrics on Fano manifolds \cite{yau,tia}, 
 was reformulated \cite{don} for constant scalar curvature Kähler metrics (an important special case of extremal metrics), 
 and finally adapted to extremal Kähler metrics \cite{mab,sze2}. 
 This problem is still widely open in the ``if'' direction, 
 except for the notable case of its specialisation to Kähler-Einstein metrics on Fano manifolds,  
 see \cite{cds1,cds2,cds3} and \cite{tia2}. 
 
 Within the scope of finding necessary conditions for the existence of extremal metrics, 
 this article provides constraints to the existence of extremal Kähler metrics with \textit{cusp singularities} 
 along a divisor in a compact Kähler manifold.  
 Cusp singularities are compatible with the extremal condition, 
 in the sense that such singular canonical metrics have already been produced \cite{sze1};   
 they appear more specifically along a continuity path between \textit{stable} and \textit{unstable} polarisations,  
 when following \textit{smooth} extremal metrics.   
 We believe in this respect that extremal Kähler metrics with cusp singularities might be of crucial interest in the study of Conjecture \ref{conj_YTD}, 
 as particular degenerations of smooth extremal metrics. 
 
 Following \cite{auv1,auv2} 
 for the definition of the class of metrics we are investigating, 
 fix a simple normal crossing divisor $D$ in a 
 compact Kähler manifold $(X,J,\omega_X)$, $\dim_{\C}X=m$, of $X$ 
 of polydiscs $U$ of holomorphic coordinates $(z^1,\dots,z^m)$ of radius $\frac{1}{2}$, 
 such that $U\cap D=\{z^1\cdots z^k=0\}$ for some $k=k(U)\in\{0,\dots, m\}$. 
  \begin{dfn}   \label{df_PKmtrcs}
   Let $\omega$ be a smooth $(1,1)$-form on $\XD$. 
   We say that $\omega$ is a \emph{Poincaré type Kähler metric} if for all $U$ and $k$ as above, 
   $\omega$ is quasi-isometric to 
   the product $\sum_{j=1}^k\frac{idz^j\wedge d\overline{z^j}}{|z^j|^2\log^2(|z^j|)}+\sum_{j=k+1}^m idz^j\wedge d\overline{z^j}$, 
   and has bounded derivatives at any order with respect to this model on $\UD$. 
   
   We say moreover that \emph{$\omega$ has class $[\omega_X]$} if $\omega=\omega_X+dd^c\varphi$, where $\varphi$ is smooth on $\XD$,  
   $\varphi=\mathcal{O}\big(1+\sum_{j=1}^k\log[-\log(|z^j|)]\big)$, 
   and $\varphi$ has bounded derivatives at any positive order for the model metric, in the above charts. 
  \end{dfn}

 Notice that this definition allows a rather loose behaviour near the divisor, in the sense that one can easily produce Poincaré type metrics 
 such that their restrictions to directions parallel to the divisor \textit{does not} converge near the divisor. 
 Our first main result states nonetheless that such a convergence does occur for extremal Poincaré type Kähler metrics, when $D\subset X$ is smooth: 
  \begin{thrm}   \label{thm_mainthm1}
   Assume that $\omega$ is an extremal Kähler Poincaré type metric of class $[\omega_X]$ on the complement of a smooth divisor 
   $D=\sum_{j=1}^N D_j$ in a compact Kähler manifold $(X,\omega_X)$. 
   Then for all $j$ there exist $a_j>0$, $\delta>0$, and a metric $\omega_j\in [\omega_X|_{D_j}]$ such that 
   on any open subset $U$ of coordinates $(z^1,z^2\dots,z^m)$ such that $U\cap D_j=\{z^1=0\}$, 
   then $\omega=\frac{a_j idz^1\wedge d\overline{z^1}}{|z^1|^2\log^2(|z^1|)}+p^*\omega_j+\mathcal{O}\big(\big|\log(|z^1|)\big|^{-\delta}\big)$ 
   as $z^1\to 0$. 
  \end{thrm}
 
 Here $p(z^1,z^2,\dots,z^m)=(z^2,\dots,z^m)$ in $U$, 
 and the $\mathcal{O}$ is understood at any order with respect to $\frac{idz^1\wedge d\overline{z^1}}{|z^1|^2\log^2(|z^1|)}+p^*\omega_j$. 
 
 One easily sees from this that the induced metrics $\omega_j$ are extremal, 
 and even have constant scalar curvature if $\omega$ does; 
 in particular, the existence of a \textit{canonical} Poincaré type metric on $\XD$ implies the existence of a \textit{canonical} metric on the components of $D$, 
 canonical meaning either \textit{extremal} or \textit{with constant scalar curvature}. 
 This implication actually holds when $D$ is no longer assumed smooth: 
  \begin{thrm}    \label{thm_mainthm2}
   Assume that there exists an extremal (resp. a constant scalar curvature) Poincaré type Kähler metric $\omega$ of class $[\omega_X]$ on the complement of a simple normal crossing divisor 
   $D=\sum_{j=1}^N D_j$ in a compact Kähler manifold $(X,\omega_X)$. 
   Then for all $j$, there exists an extremal (resp. a constant scalar curvature) Kähler metric on $D_j\backslash\sum_{\ell\neq j}D_{\ell}$ of class $[\omega_X|_{D_j}]$, 
   of Poincaré type if $D_j\cap\sum_{\ell\neq j}D_{\ell}\neq \varnothing$. 
  \end{thrm}
  
 ~
 
  Theorem \ref{thm_mainthm1} states that extremal Kähler metrics of Poincaré type are \textit{asymptotically products} 
  near the divisor. 
  Similar results for Kähler-Einstein metrics were already known \cite{sch,wu}; 
  these previous approaches differ fundamentally to ours though. 
  Indeed, in that case, 
  the Kähler-Einstein analogue of Theorem \ref{thm_mainthm2} follows from 
  topological reasons and Tian-Yau's extension \cite{ty} of Aubin-Yau theorem. 
  Hence, starting with a Poincaré type metric with asymptotically product behaviour, 
  \textit{inducing on the divisor the Kähler-Einstein metrics},
  and running Tian-Yau's continuity method towards the Kähler-Einstein metric on $\XD$, 
  G. Schumacher and D. Wu prove, roughly speaking, that the asymptotics of the metrics are preserved under the continuity path. 
  In the wider extremal case, 
  the schematic implication 
  ``existence of a canonical metric on $\XD$ $\Rightarrow$ existence of a canonical metric on $D$'' must be proven by different means, 
  as there is no such construction as Tian-Yau's for extremal metrics.  
  This illustrates the interest of Theorem \ref{thm_mainthm2}; 
  this also suggests why our proof of Theorem \ref{thm_mainthm1}, based on a good understanding of a model (divisor)$\times$(punctured unit disc in $\C$), 
  and a weighted analysis of a Lichnerowicz fourth-order operator near the divisor, is essentially different from Schumacher and Wu's proofs. 
  Let us specify also here that Theorem \ref{thm_mainthm1} is limited to the smooth divisor case so far,  
  due to the weighted analysis not transposing clearly to the normal crossing case.  

  One can interpret Theorems \ref{thm_mainthm1} and \ref{thm_mainthm2} as giving constraints on extremal Kähler metrics of Poincaré type; 
  in this way, a conjecture analogous to Conjecture \ref{conj_YTD} on Poincaré type metrics should keep track of this heredity property in the stability conditions; 
  see nonetheless the conjecture in \cite[\S3.2]{sze1} in the constant scalar curvature case. 
  Still in this particular case, in view of Theorem \ref{thm_mainthm2}, the topological constraint obtained in \cite{auv2} propagates to higher codimensional crossings, 
  giving further obstructions to the existence of constant scalar curvature Kähler metrics of Poincaré type. 
  Finally, Theorem \ref{thm_mainthm1} provides sharp asymptotic analytic properties of extremal Poincaré type metrics; 
  besides indicating what is the ``right'' class of ``metrics with cusp singularities'' in the extremal case, 
  this analytic prerequisite leads one to try and transpose analytic constructions of extremal metrics such as that of \cite{aps}, 
  crucial in the treatment of the ``only if'' direction of Conjecture \ref{conj_YTD}, to the Poincaré framework; 
  this will be addressed in a future paper. 
 
 ~
 
  \noindent
  \textit{Organisation of the article. ---} 
   This paper is composed of four parts. 
   In the first three parts, we focus on the constant scalar curvature case, 
   which already requires most of the techniques used in proving Theorems \ref{thm_mainthm1} and \ref{thm_mainthm2}. 
   More specifically, we analyse in Part \ref{sctn_mdl} the model for Poincaré type Kähler metrics, 
   i.e. $\mathbb{S}^1$-invariant Kähler metrics on products 
   (punctured unit disc)$\times$(complement of a divisor), 
   and prove for such metrics, with constant scalar curvature, a splitting theorem (Theorem \ref{thm_splttng}). 
   
   In Part \ref{sctn_almcscKfamily}, 
   we introduce the notion of a \textit{family of Kähler metrics of almost constant scalar curvature} 
   on a compact manifold, and construct a parametrisation in terms of automorphisms of the manifolds for such families 
   (Proposition \ref{prop_asympparam}). 
   
   Coming back to the complement of a simple normal crossing divisor in Part \ref{sctn_asympttcs}, 
   we use the results of Parts \ref{sctn_mdl} and \ref{sctn_almcscKfamily} 
   to prove the constant scalar curvature cases of Theorems \ref{thm_mainthm1} and \ref{thm_mainthm2} 
   (Theorems \ref{thm_asmpttcsPKcsc} and \ref{thm_cscKdivisor}). 
   For this we recall in Section \ref{subsec_fib} fibrations used in \cite{auv2}; 
   the link with the model of Part \ref{sctn_mdl} and the families of almost constant scalar curvature is made in Section \ref{subsec_link}, 
   where is proved Theorem \ref{thm_cscKdivisor}, 
   and the last three sections of Part \ref{sctn_asympttcs} are devoted to the weighted analysis needed for Theorem \ref{thm_asmpttcsPKcsc}. 
   
   In Part \ref{sctn_ext} we generalise what precedes to extremal Kähler metrics, 
   first on the product model in Section \ref{subsec_extmodel} where is proven the splitting theorem \ref{thm_splttngext}, 
   then on the complement of a simple normal crossing divisor in Section \ref{subsec_PKext}.

 
\section{Constant scalar curvature Kähler metrics of Poincaré type: the model case}   \label{sctn_mdl}

\noindent
\textit{Set-up and splitting theorem. ---} 
As a model of Poincaré type metrics near a divisor, 
we consider a compact Kähler manifold $(Y,J_Y,\omega_Y)$ together with a (possibly empty) simple normal crossing divisor $E=\sum_{j=1}^N E_j\subset Y$, 
and take its product with the punctured unit disc $\Delta^*\subset \C$ endowed with the standard complex structure $J_{\C}$.  
The product $\Delta^*\times (\YE)$ inherits the natural holomorphic $\mathbb{S}^1$-action on $\Delta^*$, 
if one merely declares that $\mathbb{S}^1$ acts trivially on $\YE$. 

We see $\Delta^*$ with its hyperbolic geometry:  
the \textit{reference metric} 
(or, more exactly, \textit{Kähler form} -- we shall exchange them often without more specification when there is no risk of confusion) 
is the Poincaré metric 
 \begin{equation*} \label{eqn_PcrMet}
  \omega_{\Delta^*}:= -dd^c\log\big(-\log(|z|^2)\big) = \frac{2idz\wedge d\bar{z}}{|z|^2\log^2(|z|^2)}
 \end{equation*}
-- notice that this equation makes it explicit that $\omega_{\Delta^*}$ is Einstein with negative scalar curvature -2. 
It is convenient to describe $\omega_{\Delta^*}$ with help of ``logarihtmic polar coordinates'' $(t,\vthet)\in \R\times \mathbb{S}^1$ defined 
via the writing $z=\exp\big(-\tfrac{1}{2}e^{t}-i\vthet) \in \Delta^*$, that is: 
$\vthet$ is the opposite of the standard angular coordinate $\theta$ on $\mathbb{S}^1$, and  
 \begin{equation*}  \label{eqn_defn_t}
  t := \log\big(-\log(|z|^2)\big). 
 \end{equation*}
This way $J_{\C}dt=2e^{-t}d\vthet$, and thus 
  $\omega_{\Delta^*} = -dd^c t = -d\big( 2e^{-t}d\vthet\big) = dt\wedge 2e^{-t}d\vthet$.

On the factor $\YE$, we fix a Poincaré type Kähler metric $\omega_{\YE}$ of class $[\omega_Y]$ according to Definition \ref{df_PKmtrcs} -- 
such an $\omega_{\YE}$ always exists, take for instance $\omega_{\YE}=\omega_Y-dd^c\mathfrak{u}_Y$, 
where $\mathfrak{u}_Y=\sum_{j=1}^N \log\big(-\log(|\sigma_j|_j^2)\big)$, 
with $\sigma_j$ a section of $\mathscr{O}\big([E_j]\big)$ canonically associated to $E_j$,  and 
$|\cdot|_j$ a well-chosen smooth hermitian metric on $[E_j]$ such that $|\sigma_j|_j^2\leq e^{-1}$ on $Y$; 
see \cite[\S1.1]{auv1} for precisions. 

We now endow $\Delta^*\times(\YE)$ with $\omega_0:=\omegadelta + \omega_{\YE}$, 
and consider the set of $\mathbb{S}^1$-invariant potentials of Kähler metrics on $\Delta^*\times(\YE)$ 
quasi-isometric to $\omega_0$, 
and whose derivatives at any order with respect to this model metric are bounded; 
we restrict more specifically to those potentials uniformly dominated by $1+\mathfrak{u}_Y$  
($\mathfrak{u}_Y$ extended constantly along $\Delta^*$), with bounded derivatives of positive order for $\omega_0$. 
In a nutshell, we look at the space:
 \begin{align*}
  \mathscr{K}(\omega_0) := \big\{ \varphi \in \mathscr{E}_0\big(\Delta^*\times&(\YE)\big)\big| \\
                                   C^{-1}&\omega_0\leq \omega_{\varphi}:=\omega_0+dd^c\varphi\leq C\omega_0 
                                    \text{ for some constant } C >0\big\}, 
 \end{align*}
with $\mathscr{E}_0\big(\Delta^*\times(\YE)\big)$ the set of $\mathbb{S}^1$-invariant 
-- emphasised through the 0 index --
smooth functions $v$ on $\Delta^*\times(\YE)$ such that 
  $|v|\leq C(1+\mathfrak{u}_Y)$ for some constant $C$,  
and for all $k, \ell\geq 0$ such that $k+\ell\geq 1$, $|\nabla^k\partial_t^{\ell} v|_{\omega_0} \leq C_{k,\ell}$, 
with $\nabla$ the Levi-Civita connexion of $\omega_{\YE}$, 
for some constant $C_{k,\ell}$.

For $\varphi \in \mathscr{K}$ -- from now on, the reference metric, fixed, is omitted --, 
we use as above and along all this part the notation $\omega_{\varphi}=\omega_0+dd^c\varphi$; 
we refer to the resulting metrics as 
\textit{Poincaré type Kähler metrics on} 
$\Delta^*\times(Y\backslash E)$, by analogy with Poincaré type Kähler metrics on 
complements of divisors in compact Kähler manifolds. 
The main result of this part deals with those $\omega_{\varphi}$ with \textit{constant scalar curvature}:  
 \begin{thm}   \label{thm_splttng}
  Assume that there exists $\varphi\in \mathscr{K}$ such that $\omega_{\varphi}$ has constant scalar curvature. 
  Then $\varphi$ does not depend on $t$, 
  and is a Poincaré type potential $\psi$ for $\omega_{\YE}$;   
  therefore, $\omega_{\varphi}$ splits as a product $\omega_{\Delta^*}+\omega_{Y}^{\psi}$, 
  with $\omega_{Y}^{\psi}:=\omega_{\YE}+dd^c_Y\psi$ a constant scalar curvature metric, 
  of Poincaré type if $E\neq\varnothing$, and of class $[\omega_Y]$.  
 \end{thm}
 
 Recall that given any Kähler metric $\omega$ on an $m$-complex dimensional manifold, 
 its scalar curvature $\scal(\omega)$ is given by the formula
  \begin{equation*} \label{eqn_CScurv}
   2m\varrho(\omega)\wedge\omega^{m-1} = \scal(\omega)\omega^{m},
  \end{equation*} 
 In the situation of Theorem \ref{thm_splttng}, if $m=\dim_{\C}(Y)+1$ and $\varrho_{\varphi}$ is the Ricci form of $\omega_{\varphi}$, 
 one thus has $ \varrho_{\varphi}\wedge\omega_{\varphi}^{m-1} = \frac{1}{2m} \sbar\omega_{\varphi}^{m}$ for some constant $\sbar$. 

Theorem \ref{thm_splttng} states a \textit{splitting principle} for constant scalar curvature metrics on products $\Delta^*\times(Y\backslash E)$, 
with Poincaré behaviour in the $\Delta^*$ direction, 
as well as in the $(Y\backslash E)$ direction when $E$ is non-trivial,  
and can thus be viewed in the same scope as the main results of \cite{ah, hua}.  
Notice that no existence of constant scalar curvature Kähler metrics (of Poincaré type) 
of class $[\omega_Y]$ on $\YE$ is a priori assumed in the statement; 
notice also that we make an implicit use of the general equivalence ``a product metric has constant scalar curvature if and only if its components do'', 
automatic in Riemannian geometry. 

Moreover, 
if one thinks to $Y$ as some component, $D_1$ say, of a simple normal crossing divisor $D=\sum_{j=1}^N D_j$ in a compact Kähler manifold $X$, 
and $E$ as the induced divisor $E_1:=\sum_{j=2}^N (D_j\cap D_1)$, then Poincaré type Kähler metrics on $\Delta^*\times(Y\backslash E)$ 
are roughly speaking asymptotic models for Poincaré type Kähler metrics on $X\backslash D$ near $D_1$. 
Heuristically, constant scalar curvature Poincaré type metrics on $\XD$ are thus modelled on products near the $D_j$, 
which thus admit constant scalar curvature metrics; 
as is seen in Part \ref{sctn_asympttcs}, the first property indeed holds if $D$ is smooth, 
and the second one holds in general (Theorems \ref{thm_asmpttcsPKcsc} and \ref{thm_cscKdivisor}). 

Our last comment concerns the class of potentials $\mathscr{K}$; 
we could have chosen, in order to respect more closely the analogy with Definition \ref{df_PKmtrcs}, 
a similar definition but with a $C^0$-bound of type $|\varphi|\leq C(|t|+\mathfrak{u}_Y)$. 
However, starting with an $\omega_{\varphi}$ with constant scalar curvature $\sbar$ and using the same integral techniques as in \cite{auv2}, 
we would have ended up with $|\varphi-at|\leq C(|t|+\mathfrak{u}_Y)$ for some $a<1$, 
\textit{completely determined by the data}: $\frac{2}{1-a}=\sbar_{\YE}-\sbar$, 
with $\sbar_{\YE}$ the mean scalar curvature attached to Poincaré type Kähler metrics of class $[\omega_Y]$ on $\YE$. 
Up to a replacing $\omega_0$ by $\frac{1}{1-a}\omega_0$, there is thus no loss of generality with our choice for $\mathscr{K}$. 
In the extremal case, one has to establish such a priori asymptotics for the potential, 
which is thus taken in a larger space as sketched above, see Section \ref{subsec_extmodel}. 

The rest of this part is devoted to the proof of Theorem \ref{thm_splttng}.

 ~

 \noindent
 \textit{A fourth order equation on $\partial_t\varphi$. ---} 
  The first step towards Theorem \ref{thm_splttng} is:
   \begin{lem}         \label{lem_lich_v}
    For $\varphi$ as in Theorem \ref{thm_splttng}, 
    denote by $\lich_{\varphi}$ the Lichnerowicz operator of order 4 associated to $\omega_{\varphi}$, 
    and set $v_{\varphi}=\dot{\varphi}-1$. 
    Then: 
   \begin{equation}  \label{eqn_lich_v}
    \lich_{\varphi}(v_{\varphi}) = 0. 
   \end{equation}
 \end{lem}

We use here the notation $\dot{\,}$ for the $t$ derivative; 
we use it again frequently, as well as its twice iterated version $\ddot{\,}$, 
in what follows. 
Besides, \eqref{eqn_lich_v} is of course equivalent to $\lich_{\varphi}(\dot{\varphi}) = 0$, but  
the ``$v_{\varphi}$-shape'' is more convenient, as we shall see below. 

~

\noindent
\prf \textit{of Lemma \ref{lem_lich_v}. ---}
For any $\phi\in \mathscr{E}_0$, 
since $\omega_{\varphi}$ has constant scalar curvature, 
the Lichnerowicz operator describes the variation of scalar curvature along a deformation of the metric in the $dd^c\phi$ direction: 
for $\vareps$ small, 
 $\scal\big(\omega_{\varphi}+dd^c(\vareps\phi)\big)=
  \scal(\omega_{\varphi})+\vareps\lich_{\varphi}(\phi)+\mathcal{O}(\vareps^2)$,  
where $\scal(\omega_{\varphi})=\sbar$ is constant; 
more generally, 
if $(\varphi_{\vareps})$ is a path in $\mathscr{K}$ with $\varphi_0=\varphi$ and 
$\phi=\frac{d\varphi_{\vareps}}{d\vareps}\big|_{\vareps=0}$, 
then $\frac{d\scal(\omega_{\varphi_{\vareps}})}{d\vareps}\big|_{\vareps=0}=\lich_{\varphi}(\phi)$. 
Notice that this holds locally if $\phi$ is only locally defined. 

Recall the complex coordinate $z=\exp\big(-\tfrac{1}{2}e^{t}-i\vthet\big)$ on $\Delta^*$, 
and consider the (locally defined) real holomorphic vector field 
 \begin{equation*}   \label{eqn_defZ}
  Z:=\Re\Big[z(\log z)\frac{\partial }{\partial z}\Big]; 
 \end{equation*}
then $Z=\Re(\log z)\Re\big(z\frac{\partial }{\partial z}\big)-\Im(\log z)\Im\big(z\frac{\partial }{\partial z}\big)$, 
$\Re(\log z)= \log|z|=-\frac{1}{2}e^{t}$, $\Im(\log z)=-\vthet$ up to $2\pi$, 
and $z\frac{\partial }{\partial z}=-e^{-t}\frac{\partial}{\partial t}+\frac{i}{2}\frac{\partial}{\partial \vthet}$, 
thus: 
 \begin{equation*}   \label{eqn_Z}
  Z=  \frac{1}{2}\frac{\partial}{\partial t} 
       + \frac{1}{2}\vthet\frac{\partial}{\partial \vthet}
       \qquad\text{up to}\quad \pi\frac{\partial}{\partial \vthet}. 
 \end{equation*}
In particular, $Z\cdot f$ 
makes global sense as $\frac{1}{2}\dot{f}$ for any $\mathbb{S}^1$-invariant $f$.  
Now as $\omega_{\varphi}=\omega_{\YE}+dd^c(\varphi-t)$, 
$\mathcal{L}_Z\omega_{\varphi}=\mathcal{L}_Z\omega_{\YE}+dd^c\big(Z\cdot(\varphi-t)\big)=dd^c\big(Z\cdot(\varphi-t)\big)$, 
as $Z$ is normal to $Y$. 
Moreover, $(\varphi-t)$ is $\mathbb{S}^1$-invariant, 
and $Z\cdot(\varphi-t)$, globally defined, equals $\frac{1}{2}(\dot{\varphi}-1)=\frac{1}{2}v_{\varphi}$. 
Therefore by the preliminary remark, denoting by $\Phi_{\vareps}^Z$ the flow of $Z$, 
one has, for $\vareps$ small: 
 \begin{align*}
  \bar{\scal}           = (\Phi_{\vareps}^Z)^*\scal(\omega_\varphi)
                        = \scal\big((\Phi_{\vareps}^Z)^*\omega_{\varphi}\big) 
                        = \scal\big(\omega_{\varphi}+\frac{\vareps}{2}dd^c v_{\varphi}+\mathcal{O}(\vareps^2)\big) 
                        = \bar{\scal} + \frac{\vareps}{2}\lich_{\varphi}(v_{\varphi})+\mathcal{O}(\vareps^2), 
 \end{align*}
and thus 
$\lich_{\varphi}(v_{\varphi})$ 
vanishes identically. 
\cqfd

 ~

\noindent
\textit{A useful holomorphic gradient. ---}  
Recall that the Lichnerowicz operator is self-adjoint by construction, 
as it can be defined 
-- independently of $\omega_{\varphi}$ having constant scalar curvature -- 
as $\mathcal{D}_{\varphi}^*\mathcal{D}_{\varphi}$, 
with $\mathcal{D}_{\varphi}=(\nabla^{\varphi})^{-}d$, 
where $(\nabla^{\varphi})^{-}$ is the $J$-anti-invariant part of the Levi-Civita connection of $\omega_{\varphi}$, 
and $\mathcal{D}_{\varphi}^*$ the formal adjoint of $\mathcal{D}_{\varphi}$ for $\omega_{\varphi}$. 
On compact manifolds, $\mathbb{L}$ and $\mathcal{D}$ thus have the same kernel;   
this comes at once from an integration by parts, 
and cannot therefore be applied directly on $\Delta^*\times(Y\backslash E)$ in general. 
Our aim is to prove, however, after Lemma \ref{lem_lich_v}, 
that indeed, $\mathcal{D}_{\varphi}(v_{\varphi})=0$. 
An important intermediate step for this is: 
 \begin{lem} \label{lem_D_u}
  For $\varphi\in \mathscr{K}$, 
  set $u_{\varphi}=e^{-t}(\dot{\varphi}-1)=e^{-t}v_{\varphi}$. 
  Then: 
   \begin{equation}  \label{eqn_D_u}
    \mathcal{D}_{\varphi}(u_{\varphi}) = 0;  
   \end{equation}
  in particular, $\lich_{\varphi}(u_{\varphi})=0$, 
  and $\mathcal{D}_{\varphi}(v_{\varphi})=\frac{1}{2}\big(du_{\varphi}\cdot de^t- d^cu_{\varphi}\cdot d^ce^t\big)+u_{\varphi}\mathcal{D}_{\varphi}(e^t)$.   
 \end{lem}
  In this statement and from now on, 
  we adopt the convention that $\alpha\cdot\beta=\alpha\otimes\beta+\beta\otimes\alpha$ for any 1-forms $\alpha$ and $\beta$; 
  $\alpha^2$ always means $\alpha\otimes\alpha$ though. 

~

 \prf \textit{of Lemma \ref{lem_D_u}. ---}  
  When $\scal(\omega_{\varphi})$ is constant, the equation $\lich_{\varphi}(u_{\varphi})=0$ can be obtained 
  in the same way as the equation $\lich_{\varphi}(v_{\varphi})=0$ of Lemma \ref{lem_lich_v}, 
  using the real holomorphic vector field $\Re\big(z\frac{\partial}{\partial z}\big)$ 
  instead of 
  $\Re\big(z\log(z)\frac{\partial}{\partial z}\big)$. 
  Now in our context $\lich_{\varphi}(u_{\varphi})=0$ is not enough to deduce $\mathcal{D}_{\varphi}(u_{\varphi})=0$, 
  as for instance $u_{\varphi}$ a priori has size $e^{-t}$ for $t$ going to $-\infty$, 
  which brings up problematic boundary terms if one tries and performs the usual integrations by parts. 
  Equation \eqref{eqn_D_u} actually comes from a more direct computation, and holds in general, 
  i.e. \textit{independently of $\omega_{\varphi}$ having constant scalar curvature}.  
  
  For any (twice differentiable, say) function $f$, the equation $\mathcal{D}_{\varphi}(f)=0$ is indeed 
  equivalent to $\nabla^{\varphi}f$ being a real holomorphic vector field, 
  where $\nabla^{\varphi}f$ denotes the gradient of $f$ computed with respect to $g_{\varphi}=\omega_{\varphi}(\cdot,J\cdot)$ 
  -- there should be no confusion between our two different uses of $\nabla^{\varphi}$, 
  as it refers to a gradient only when used with functions, 
  and as we always denote differentials by $d$. 
  So if we check the structural equation  
   \begin{equation}   \label{eqn_nablau=holo}
    \nabla^{\varphi}u_{\varphi}=e^{-t}\frac{\partial}{\partial t}, 
   \end{equation}
  then we are done, since $e^{-t}\frac{\partial}{\partial t}=-\Re\big(z\frac{\partial}{\partial z}\big)$,  
  as seen in the previous proof. 
  
  According to the splitting $J=J_{\C}\oplus J_Y$ and the rule 
  $J_{\C}dt=2e^{-t}d\vthet$, and since $\varphi$ is $\mathbb{S}^1$-invariant, 
  we have: 
   \begin{equation} \label{eqn_omegavarphi1}
      \omega_{\varphi} = (1+\ddot{\varphi}-\dot{\varphi})dt\wedge2e^{-t}d\vthet 
                          + dt\wedge d^c_Y\dot{\varphi}                            
                          + d_Y\dot{\varphi}\wedge 2e^{-t}d\vthet
                          + \big(\omega_{\YE}+dd^c_Y\varphi\big),
   \end{equation}
  where $d_Y$, $d^c_Y$ and $dd^c_Y$ are respectively the operators $d$, $d^c$ and $dd^c$ acting on functions on $Y$ 
  -- or, for instance: $(d_Yf)(t,\vthet,\cdot)=d\big(f(t,\vthet,\cdot)\big)$, 
  and so on.
  
  Given any 1-form $\alpha$ 
  and any function $f$, 
  one has: 
   \begin{equation*}
    \alpha(\nabla^{\varphi}f)\frac{\omega_{\varphi}^m}{m!} 
                   = \langle\alpha, df\rangle_{g_{\varphi}}\frac{\omega_{\varphi}^m}{m!}
                   = \alpha\wedge d^cf \wedge \frac{\omega_{\varphi}^{m-1}}{(m-1)!}. 
   \end{equation*}
 Now observe that if one takes $f=e^{-t}(\dot{\varphi}-1)=:u_{\varphi}$, 
 giving thus 
  \begin{equation*}
   d^cu_{\varphi}=\frac{\partial u_{\varphi}}{\partial t}2e^{-t}d\vthet+d^c_Yu_{\varphi}
                =e^{-t}\big((1+\ddot{\varphi}-\dot{\varphi})2e^{-t}d\vthet + d^c_Y\dot{\varphi}\big), 
  \end{equation*}
 one has, in view of \eqref{eqn_omegavarphi1} and using $\mathbb{S}^1$-invariance to set $\omega_{t}^{\varphi}=(\omega_{\YE}+dd^c_Y\varphi)|_{t,\vthet}$:
  \begin{equation*}  \label{eqn_gradu} 
   \begin{aligned}
    d^cu_{\varphi}\wedge \frac{\omega_{\varphi}^{m-1}}{(m-1)!}
       = e^{-t}\Big[
        (1+\ddot{\varphi}-\dot{\varphi}&)2e^{-t}d\vthet\wedge dt\wedge d^c_Y\dot{\varphi}\wedge\frac{(\omega_{t}^{\varphi})^{m-2}}{(m-2)!} \\
        &+(1+\ddot{\varphi}-\dot{\varphi})2e^{-t}d\vthet\wedge\frac{(\omega_{t}^{\varphi})^{m-1}}{(m-1)!}                                   \\
        &+d^c_Y\dot{\varphi}\wedge(1+\ddot{\varphi}-\dot{\varphi})dt\wedge2e^{-t}d\vthet\wedge \frac{(\omega_{t}^{\varphi})^{m-2}}{(m-2)!}  \\
        &+d^c_Y\dot{\varphi}\wedge d_Y\dot{\varphi}\wedge 2e^{-t}d\vthet\wedge\frac{(\omega_{t}^{\varphi})^{m-2}}{(m-2)!}    
                                 \Big]                                                                                                       \\
       =e^{-t}\Big[(1+\ddot{\varphi}-\dot{\varphi}&)2e^{-t}d\vthet\wedge\frac{(\omega_{t}^{\varphi})^{m-1}}{(m-1)!}                         \\
                  &-2e^{-t}d\vthet\wedge d_Y\dot{\varphi}\wedge d^c_Y\dot{\varphi}\wedge \frac{(\omega_{t}^{\varphi})^{m-2}}{(m-2)!}
                   \Big],                         
   \end{aligned}
  \end{equation*}
 since the first and third lines of the right-hand side of the first equality cancel each other. 
 Therefore, for any 1-form $\alpha$ on $\Delta^*\times(Y\backslash E)$ written as $\alpha_tdt+\alpha_{\vthet}d\vthet+\alpha_Y$, 
  \begin{equation*}
   \begin{aligned}
    \alpha(\nabla^{\varphi}u_{\varphi})\frac{\omega_{\varphi}^m}{m!} 
      = e^{-t}\alpha_tdt\wedge\Big[(1+\ddot{\varphi}-\dot{\varphi})&2e^{-t}d\vthet
                                   \wedge\frac{(\omega_{t}^{\varphi})^{m-1}}{(m-1)!}               \\                      
                                    -&2e^{-t}d\vthet\wedge d_Y\dot{\varphi}\wedge d^c_Y\dot{\varphi}
                                   \wedge \frac{(\omega_{t}^{\varphi})^{m-2}}{(m-2)!}
                              \Big]. 
   \end{aligned}
  \end{equation*}
 On the other hand, a direct computation yields  
  \begin{align*}
   \frac{\omega_{\varphi}^m}{m!}=  (1+\ddot{\varphi}-\dot{\varphi})dt\wedge 2e^{-t}d\vthet
                                   \wedge\frac{(\omega_{t}^{\varphi})^{m-1}}{(m-1)!}                        
                                   -dt\wedge 2e^{-t}d\vthet\wedge d_Y\dot{\varphi}\wedge d^c_Y\dot{\varphi}
                                   \wedge \frac{(\omega_{t}^{\varphi})^{m-2}}{(m-2)!}, 
  \end{align*}
 hence $\alpha(\nabla^{\varphi}u_{\varphi})
          =e^{-t}\alpha_t=\alpha\big(e^{-t}\tfrac{\partial}{\partial t}\big)$ for 
 any 1-form $\alpha$: 
 equation is \eqref{eqn_nablau=holo} verified.  

 Knowing that $\mathcal{D}_{\varphi}(u_{\varphi})=0$, 
 the assertion on $\mathcal{D}_{\varphi}(v_{\varphi})$ now directly comes from the definitions 
 of $\mathcal{D}_{\varphi}=(\nabla^{\varphi})^-d$ and $v_{\varphi}=e^tu_{\varphi}$, and Leibniz rule. 
 \cqfd
  
  ~

\noindent
\textit{Finiteness of a weighted $L^2$ norm of $\mathcal{D}_{\varphi}(v_{\varphi})$. ---}  
 For general $\phi\in \mathscr{K}$, 
 we only know that $\mathcal{D}_{\phi}(\dot{\phi})$ is bounded on $\Delta^*\times(Y\backslash E)$;  
 since $e^t\frac{\omega_{\phi}^m}{m!}$ is mutually bounded with the \textit{cylindrical} volume form 
 $dt\wedge d\vthet\wedge \tfrac{\omega_{\YE}^{m-1}}{(m-1)!}$,  
 $\int_{\Delta^*\times(Y\backslash E)}e^t\big|\mathcal{D}_{\phi}(\dot{\phi})\big|^2_{\phi}\frac{\omega_{\phi}^m}{m!}$ 
 has no reason to be finite -- here $|\cdot|_{\phi}$ denotes the norm computed with $g_{\phi}$. 
 Combining equalities \eqref{eqn_lich_v} and \eqref{eqn_D_u}, 
 we claim that this is indeed the case for $\varphi$ such that $\scal(\omega_{\varphi})$ 
 is constant:
  \begin{lem}  \label{lem_Dv_isL2} 
   For $\varphi$ as in Theorem \ref{thm_splttng}, 
    \begin{equation*}  \label{eqn_Dv_isL2}
     \int_{\Delta^*\times(Y\backslash E)}
             e^t\big|\mathcal{D}_{\varphi}(v_{\varphi})\big|^2_{\varphi}\vol^{\varphi} <\infty,
    \end{equation*}
   where $\vol^{\varphi}=\frac{\omega_{\varphi}^m}{m!}$, and where we recall the notation: $v_{\varphi}=\dot{\varphi}-1$.
  \end{lem}
 \prf\textit{ of Lemma \ref{lem_Dv_isL2}. ---} 
  For $s\geq 0$, set $\Delta_s:= \{|t|\leq s\}\subset \Delta^*$. 
  First we relate 
   $\int_{\Delta_s\times(Y\backslash E)}
             e^t\big|\mathcal{D}_{\varphi}(v_{\varphi})\big|^2_{\varphi}\vol^{\varphi}$ 
  to $\int_{\Delta_s\times(Y\backslash E)}e^t v_{\varphi}\lich_{\varphi}(v_{\varphi})\vol^{\varphi}$, 
  where $\vol^{\varphi}=\frac{\omega_{\varphi}^m}{m!}$; 
  notice though that the latter integral is always 0, by Lemma \ref{lem_lich_v}. 
  Assume momentarily that $E$ is empty, so that we work on $\Delta_s\times Y$. 
  From \eqref{eqn_lich_v} and as $\lich_{\varphi}=\delta^{\varphi}\delta^{\varphi}\mathcal{D}_{\varphi}$, 
  one has: 
   \begin{align*}
    0=\int_{\Delta_s\times Y }e^tv_{\varphi}\lich_{\varphi}(v_{\varphi})\vol^{\varphi}  
        =&  \int_{\Delta_s\times Y }
               e^tv_{\varphi}\delta^{\varphi}\delta^{\varphi}\mathcal{D}_{\varphi}(v_{\varphi})\vol^{\varphi}                                   \\
        =& \int_{\Delta_s\times Y }
            \big\langle d(e^tv_{\varphi}),\delta^{\varphi}\mathcal{D}_{\varphi}(v_{\varphi})\big\rangle_{\varphi}\vol^{\varphi}     \\                  
         &+\int_{\{t=s\}}e^{s}
                 \mathfrak{D}_{s}^{\varphi}\wedge 2e^{-s}d\vthet
            -\int_{\{t=-s\}}e^{-s}
                 \mathfrak{D}_{-s}^{\varphi} \wedge 2e^{s}d\vthet,      
   \end{align*}
  where 
   \begin{equation*}   \label{eqn_gothD}
    \mathfrak{D}_{t}^{\varphi}:=\delta^{\varphi}\mathcal{D}_{\varphi}(v_{\varphi})(\partial_t)\frac{(\omega_{t}^{\varphi})^{m-1}}{(m-1)!}
                 -\delta^{\varphi}\mathcal{D}_{\varphi}(v_{\varphi})|_Y
                 \wedge d^c_Yv_{\varphi}\wedge\frac{(\omega_{t}^{\varphi})^{m-2}}{(m-2)!},
   \end{equation*}
    with $\omega_{s}^{\varphi}=\omega_{\YE}+dd^c_Y\big(\varphi(s,\cdot)\big)$.  
    Here we use  Stokes' theorem, and the Kähler identities 
    $\langle\alpha,\beta\rangle_{\omega}\frac{\omega^m}{m!}=\alpha\wedge J\beta\wedge\frac{\omega^{m-1}}{(m-1)!}$ 
    and $(\delta\alpha) \frac{\omega^m}{m!}= - d^c\alpha\wedge\frac{\omega^{m-1}}{(m-1)!}$ for 1-forms, 
    as well as $\omega_{\varphi}=\omega_{\pm s}^{\varphi}+d_Y\dot{\varphi}\wedge 2e^{\mp s}d\vthet$ on slices $\{t=\pm s\}$, 
    to compute the boundary integrals. 
    Observe that the $e^s$ and $e^{-s}$ cancel each other in these terms, 
    and therefore the integrands are bounded -- for the metric $g_{\YE}+d\vthet^2$, say -- 
    independently of $s$.  
    Consequently the boundary integrals are $\mathcal{O}(1)$, that is:
     \begin{equation}   \label{eqn_ibp1}
      0=\int_{\Delta_s\times Y }e^t v_{\varphi}\lich_{\varphi}(v_{\varphi})\vol^{\varphi}
       = \int_{\Delta_s\times Y }
            \big\langle d(e^tv_{\varphi}),\delta^{\varphi}\mathcal{D}_{\varphi}(v_{\varphi})\big\rangle_{\varphi}\vol^{\varphi}
         +\mathcal{O}(1),
     \end{equation}
    the $\mathcal{O}(1)$ being understood with respect to the variable $s$. 
    
    We proceed to a further integration by parts, 
    using that by definition the $\delta^{\varphi}$ in the second integral of the right-hand side of \eqref{eqn_ibp1} 
    is the adjoint of the projection of the Levi-Civita connexion $\nabla^{\varphi}$ from 1-forms and to symmetric 2-forms: 
     \begin{equation}  \label{eqn_ibp2}
      0=\int_{\Delta_s\times Y }e^t v_{\varphi}\lich_{\varphi}(v_{\varphi})\vol^{\varphi}
       = \int_{\Delta_s\times Y }
            \big\langle \nabla^{\varphi}d(e^tv_{\varphi}),
                 \mathcal{D}_{\varphi}(v_{\varphi})\big\rangle_{\varphi}\vol^{\varphi}
         +\mathcal{O}(1); 
     \end{equation}
    here we have included the boundary in the $\mathcal{O}(1)$, 
    since they are bounded independently of $s$ for the same reasons as for the first integration by parts above.

    As $\mathcal{D}_{\varphi}(v_{\varphi})$ is $J$-anti-invariant by construction, 
    $\nabla^{\varphi}d(e^tv_{\varphi})$ can be replaced by is $J$-anti-invariant part $\mathcal{D}_{\varphi}(e^tv_{\varphi})$
    in the inner product in the right-hand side of \eqref{eqn_ibp2}: 
     \begin{equation}  \label{eqn_ibp3}
      0=\int_{\Delta_s\times Y }e^t v_{\varphi}\lich_{\varphi}(v_{\varphi})\vol^{\varphi}
       = \int_{\Delta_s\times Y }
            \big\langle\mathcal{D}_{\varphi}(e^tv_{\varphi}),
                 \mathcal{D}_{\varphi}(v_{\varphi})\big\rangle_{\varphi}\vol^{\varphi}
         +\mathcal{O}(1); 
     \end{equation}
    Expand now $\mathcal{D}_{\varphi}(e^tv_{\varphi})$ with Leibniz rule, using that $v_{\varphi}=e^tu_{\varphi}$  
    and that $\mathcal{D}_{\varphi}(u_{\varphi})=0$: 
     \begin{equation*}
      \mathcal{D}_{\varphi}(e^tv_{\varphi}) = \mathcal{D}_{\varphi}(e^{2t}u_{\varphi})
           = 2e^{2t}(dt\cdot du_{\varphi})^- + u_{\varphi}\mathcal{D}_{\varphi}(e^{2t}); 
     \end{equation*}
    here and further on, $h^-$ denotes the $J$-anti-invariant part of any symmetric 2-form $h$. 
    Moreover, $\mathcal{D}_{\varphi}(e^{2t})=2e^{2t}(dt^2)^-+2e^t\mathcal{D}_{\varphi}(e^{t})$, 
    and 
    $\mathcal{D}_{\varphi}(v_{\varphi})= \mathcal{D}_{\varphi}(e^{t}u_{\varphi})
                                       = e^{t}(dt\cdot du_{\varphi})^- + u_{\varphi}\mathcal{D}_{\varphi}(e^{t})$
    by Lemma \ref{lem_D_u}, 
    hence: 
     \begin{equation*}
      \mathcal{D}_{\varphi}(e^tv_{\varphi})= 2e^t\mathcal{D}_{\varphi}(v_{\varphi}) + 2v_{\varphi}(dt^2)^-. 
     \end{equation*}

    From this and \eqref{eqn_ibp3} we thus infer: 
     \begin{equation}   \label{eqn_ibpeq2}
      \begin{aligned}
      \int_{\Delta_s\times Y }
             e^t|\mathcal{D}_{\varphi}(v_{\varphi})|_{\varphi}^2&\vol^{\varphi} 
              = -\frac{1}{2}\int_{\Delta_s\times Y }
                  v_{\varphi}e^t\big\langle(dt^2)^-,\mathcal{D}_{\varphi}(v_{\varphi})\big\rangle_{\varphi}\vol^{\varphi} +\mathcal{O}(1) \\
      \end{aligned}
     \end{equation}
    We shall make explicit the computations involved in the right-hand side of this estimate, 
    replacing the $\mathcal{D}_{\varphi}(v_{\varphi})$ there by its expansion given in Lemma \ref{lem_D_u}, 
    that is: $(de^t\cdot du_{\varphi})^-+u_{\varphi}\mathcal{D}_{\varphi}(e^t)$.  
    We also replace $(dt^2)^-$ by $dt^2$ in the inner product. 

    First, 
    $\big\langle dt^2,(de^t\cdot du_{\varphi})^-\big\rangle_{\varphi}$ 
    is merely equal to 
    $\frac{e^t}{2} dt(\nabla^{\varphi}t)du_{\varphi}(\nabla^{\varphi}t)\cdot 2
        =|dt|^2_{\varphi}$, 
    as $du_{\varphi}(\nabla^{\varphi}t)=dt(\nabla^{\varphi}u_{\varphi})=e^{-t}$, 
    and $d^ct(\nabla^{\varphi}t)=d^cu_{\varphi}(\nabla^{\varphi}t)=0$;   
    to see these vanishings, 
    write for example $d^ct(\nabla^{\varphi}t)\frac{\omega_{\varphi}^m}{m!}
           =\langle dt,d^ct\rangle_{\varphi}\frac{\omega_{\varphi}^m}{m!}=dt\wedge (-dt)\wedge\frac{\omega_{\varphi}^{m-1}}{(m-1)!}$ for the first one, 
    and $d^cu_{\varphi}(\nabla^{\varphi}t)=d^ct(\nabla^{\varphi}u_{\varphi})
                                             =2e^{-t}d\vthet\big(e^{-t}\frac{\partial}{\partial t}\big)=0 $ for the second one. 
                                             
    Now $\mathcal{D}_{\varphi}(e^t)=e^t(dt^2)^-+e^t\mathcal{D}_{\varphi}(t)$, 
    thus 
     \begin{align*}
      \big\langle dt^2,\mathcal{D}_{\varphi}(e^t)\big\rangle_{\varphi}
          =&\frac{e^t}{2}\big(dt(\nabla^{\varphi}t)^2+d^ct(\nabla^{\varphi}t)^2\big)
            +e^t\mathcal{D}_{\varphi}(t)(\nabla^{\varphi}t,\nabla^{\varphi}t)                               \\
          =&\frac{e^t}{2}|dt|^4_{\varphi}+ e^t\mathcal{D}_{\varphi}(t)(\nabla^{\varphi}t,\nabla^{\varphi}t); 
     \end{align*}
    we are left with the computation of $\mathcal{D}_{\varphi}(t)(\nabla^{\varphi}t,\nabla^{\varphi}t)$, 
    hence those of
    $\big(\nabla^{\varphi}_{\nabla^{\varphi}t}dt\big)(\nabla^{\varphi}t)$ and 
    $\big(\nabla^{\varphi}_{J\nabla^{\varphi}t}dt\big)(J\nabla^{\varphi}t)$.  

   \begin{lem}  \label{lem_Dt}
    One has: $\big(\nabla^{\varphi}_{\nabla^{\varphi}t}dt\big)(\nabla^{\varphi}t)
               =\frac{1}{2}(\nabla^{\varphi}t)\cdot|dt|_{\varphi}^2$, 
    and: $\big(\nabla^{\varphi}_{J\nabla^{\varphi}t}dt\big)(J\nabla^{\varphi}t)= 
            -|dt|^4_{\varphi}-\frac{1}{2}(\nabla^{\varphi}t)\cdot|dt|^2_{\varphi}$. 
   
   \end{lem}
   \prf \textit{of Lemma \ref{lem_Dt}}. --- 
    Rewrite the first quantity to compute as 
    $\big\langle\nabla^{\varphi}_{\nabla^{\varphi}t}dt,dt\big\rangle_{\varphi}$, to see that it is indeed nothing but 
    $\frac{1}{2}(\nabla^{\varphi}t)\cdot|dt|_{\varphi}^2$. 
    For the second quantity we proceed as follows: 
    $\big(\nabla^{\varphi}_{J\nabla^{\varphi}t}dt\big)(J\nabla^{\varphi}t)
     = (J\nabla^{\varphi}t)\cdot\big(dt(J\nabla^{\varphi}t)\big)
       -dt\big(\nabla^{\varphi}_{J\nabla^{\varphi}t}(J\nabla^{\varphi}t)\big)
     = - \big\langle\nabla^{\varphi}t,\nabla^{\varphi}_{J\nabla^{\varphi}t}(J\nabla^{\varphi}t) \big\rangle_{\varphi}$, 
    since $dt(J\nabla^{\varphi}t)\equiv 0$. 
    Now by Koszul formula, 
     \begin{align*}
      2\big\langle\nabla^{\varphi}t,\nabla^{\varphi}_{J\nabla^{\varphi}t}(J\nabla^{\varphi}t) \big\rangle_{\varphi}
          = &2(J\nabla^{\varphi}t)\cdot\langle \nabla^{\varphi}t,J\nabla^{\varphi}t\rangle_{\varphi}
             -(\nabla^{\varphi}t)\cdot|J\nabla^{\varphi}t|^2_{\varphi}                                  \\
            &+\big\langle[J\nabla^{\varphi}t,J\nabla^{\varphi}t],\nabla^{\varphi}t\big\rangle_{\varphi}
             +2\big\langle[\nabla^{\varphi}t,J\nabla^{\varphi}t],J\nabla^{\varphi}t\big\rangle_{\varphi}            \\
          = &-(\nabla^{\varphi}t)\cdot|dt|^2_{\varphi}+2\big\langle[\nabla^{\varphi}t,J\nabla^{\varphi}t],J\nabla^{\varphi}t\big\rangle_{\varphi}, 
     \end{align*}
    since $\langle\nabla^{\varphi}t,J\nabla^{\varphi}t\rangle_{\varphi}=0$, 
    $|J\nabla^{\varphi}t|^2_{\varphi}=|\nabla^{\varphi}t|^2_{\varphi}=|dt|^2_{\varphi}$ 
    (and of course $[J\nabla^{\varphi}t,J\nabla^{\varphi}t]=0$). 
    Now 
     \begin{align*}
      \big\langle[\nabla^{\varphi}t,J\nabla^{\varphi}t],J\nabla^{\varphi}t\big\rangle_{\varphi}
          =& d^ct\big([\nabla^{\varphi}t,J\nabla^{\varphi}t]\big) =2e^{-t}d\vthet\big([\nabla^{\varphi}t,J\nabla^{\varphi}t]\big)                                        \\
          =& 2e^{-t}\big[(\nabla^{\varphi}t)\cdot\big((J\nabla^{\varphi}t)\cdot\vthet\big)-(J\nabla^{\varphi}t)\cdot\big((\nabla^{\varphi}t)\cdot\vthet\big)\big] \\
          =& 2e^{-t}\big[(\nabla^{\varphi}t)\cdot\big(\frac{1}{2}e^t|dt|^2_{\varphi}\big)- 0 \big], 
     \end{align*}
    since $(J\nabla^{\varphi}t)\cdot\vthet=d\vthet(J\nabla^{\varphi}t)=\frac{1}{2}e^tdt(\nabla^{\varphi}t)=\frac{1}{2}e^t|dt|^2_{\varphi}$, 
    and $(\nabla^{\varphi}t)\cdot\vthet=d\vthet(\nabla^{\varphi}t)=0$ as already seen. 
    So finally $\big\langle[\nabla^{\varphi}t,J\nabla^{\varphi}t],J\nabla^{\varphi}t\big\rangle_{\varphi}=|dt|^4_{\varphi}+(\nabla^{\varphi}t)\cdot|dt|^2_{\varphi}$, 
    and thus 
    $\big(\nabla^{\varphi}_{J\nabla^{\varphi}t}dt\big)(J\nabla^{\varphi}t)
       =-\big\langle\nabla^{\varphi}t,\nabla^{\varphi}_{J\nabla^{\varphi}t}(J\nabla^{\varphi}t) \big\rangle_{\varphi}
       =-|dt|^4_{\varphi}-\frac{1}{2}(\nabla^{\varphi}t)\cdot|dt|^2_{\varphi} $.  
   \hfill $\blacksquare$

    ~

    We deduce from Lemma \ref{lem_Dt} that: 
     \begin{equation*}  \label{eqn_Dt}
      \mathcal{D}(t)(\nabla^{\varphi}t,\nabla^{\varphi}t)
       =\frac{1}{2}\big[\big(\nabla^{\varphi}_{\nabla^{\varphi}t}dt\big)(\nabla^{\varphi}t)
                        -\big(\nabla^{\varphi}_{J\nabla^{\varphi}t}dt\big)(J\nabla^{\varphi}t)\big]
       =\frac{1}{2}(\nabla^{\varphi}t)\cdot|dt|^2_{\varphi}+\frac{1}{2}|dt|^4_{\varphi}.   
     \end{equation*}
    This yields  
    $\big\langle dt^2,\mathcal{D}_{\varphi}(e^t)\big\rangle_{\varphi}=\frac{e^t}{2}(\nabla^{\varphi}t)\cdot|dt|^2_{\varphi}+e^t|dt|^4_{\varphi}$, 
    and therefore 
     \begin{equation}   \label{eqn_vDvdt2}
      v_{\varphi}\big\langle(dt^2)^-,\mathcal{D}_{\varphi}(v_{\varphi})\big\rangle_{\varphi} 
         = v_{\varphi}|dt|^2_{\varphi}+v_{\varphi}^2|dt|^4_{\varphi}+\frac{1}{2}v_{\varphi}^2(\nabla^{\varphi}t)\cdot|dt|^2_{\varphi}. 
     \end{equation}
    Do not use \eqref{eqn_vDvdt2} with \eqref{eqn_ibpeq2} yet; instead, focus on its last summand, 
    and notice that: 
     \begin{align*}
      \int_{\Delta_s\times Y }e^tv_{\varphi}^2(\nabla^{\varphi}t)\cdot|dt|^2_{\varphi}\vol^{\varphi} 
         &= \int_{\Delta_s\times Y } v_{\varphi}^2 d(|dt|^2_{\varphi})\wedge d^c(e^t)\wedge \frac{\omega_{\varphi}^{m-1}}{(m-1)!}      \\
         &=-\int_{\Delta_s\times Y }  |dt|^2_{\varphi}d(v_{\varphi}^2)\wedge d^c(e^t)\wedge \frac{\omega_{\varphi}^{m-1}}{(m-1)!}+\mathcal{O}(1)
     \end{align*}
    by Stokes, since $dd^c(e^t)=0$ -- here again, the boundary terms are bounded independently of $s$. 
    This we rewrite as 
     \begin{equation*}
      \int_{\Delta_s\times Y }e^tv_{\varphi}^2(\nabla^{\varphi}t)\cdot|dt|^2_{\varphi}\vol^{\varphi}
        = -2\int_{\Delta_s\times Y }  e^tv_{\varphi}|dt|^2_{\varphi}dt(\nabla^{\varphi}v_{\varphi})\vol^{\varphi}+\mathcal{O}(1),
     \end{equation*}
    and as $\nabla^{\varphi}v_{\varphi}=\nabla^{\varphi}(e^tu_{\varphi})
              =e^t\nabla^{\varphi}u_{\varphi}+e^tu_{\varphi}\nabla^{\varphi}t=\frac{\partial}{\partial t}+v_{\varphi}\nabla^{\varphi}t$, 
    we get 
     \begin{equation*}
      \int_{\Delta_s\times Y }e^tv_{\varphi}^2(\nabla^{\varphi}t)\cdot|dt|^2_{\varphi}\vol^{\varphi}
        = -2\int_{\Delta_s\times Y }e^tv_{\varphi}|dt|^2_{\varphi}\big(1+v_{\varphi}|dt|^2_{\varphi}\big)\vol^{\varphi}+\mathcal{O}(1). 
     \end{equation*}
    From this latter equality and \eqref{eqn_vDvdt2}, we thus exactly end up with 
     \begin{align*}
      \int_{\Delta_s\times Y }
                  &v_{\varphi}e^t\big\langle(dt^2)^-,\mathcal{D}_{\varphi}(v_{\varphi})\big\rangle_{\varphi}\vol^{\varphi}                            \\
        &= \int_{\Delta_s\times Y }
                  e^tv_{\varphi}|dt|^2_{\varphi}\big(1+v_{\varphi}|dt|^2_{\varphi}\big)\vol^{\varphi}
           +\frac{1}{2}\int_{\Delta_s\times Y }e^tv_{\varphi}^2(\nabla^{\varphi}t)\cdot|dt|^2_{\varphi}\vol^{\varphi}                    
         = \mathcal{O}(1),
     \end{align*}
    that is, coming back now to \eqref{eqn_ibpeq2}, 
     \begin{equation*}
      \int_{\Delta_s\times Y }e^t|\mathcal{D}_{\varphi}(v_{\varphi})|_{\varphi}^2 \vol^{\varphi}= \mathcal{O}(1).
     \end{equation*}
    As the integrand is nonnegative, 
    this implies that the integral converges as $s$ goes to $\infty$; 
    in other words, $\int_{\Delta\times Y }e^t|\mathcal{D}_{\varphi}(v_{\varphi})|_{\varphi}^2 \vol^{\varphi}$ is finite.

    When the divisor $E$ in $Y$ is not a priori empty, the same arguments apply, 
    by replacing $\Delta_s\times Y$ by $\Delta_s\times (Y\backslash E)$ in the integrals above. 
    One just has to check that the integrations by parts still provide bounded boundary terms, 
    which is indeed the case thanks to the Poincaré assumption on $\omega_{\varphi}$. 
    \cqfd

 ~
 
\noindent 
\textit{Vanishing of the weighted $L^2$ norm of $\mathcal{D}_{\varphi}(v_{\varphi})$. ---} 
  We strengthen Lemma \ref{lem_Dv_isL2} as: 
    \begin{lem}   \label{lem_Dv_is0}
     For $\varphi$ as in Theorem \ref{thm_splttng}, 
     $\int_{\Delta^*\times (Y\backslash E)}|\mathcal{D}_{\varphi}(v_{\varphi})|_{\varphi}^2\,e^t\vol^{\varphi}=0$, and thus $\mathcal{D}_{\varphi}(v_{\varphi})\equiv0$. 
    \end{lem}

   \prf \textit{of Lemma \ref{lem_Dv_is0}}. --- 
    Denote by $\mathcal{F}_{\varphi}$ the function 
     \begin{align*}
      s\longmapsto& \int_{\{t=s\}}
                v_{\varphi}\,\mathfrak{D}_s^{\varphi}\wedge 2d\vthet 
                 =  4\pi \int_{\{s\}\times(\YE)} v_{\varphi}\,\mathfrak{D}_s^{\varphi}
     \end{align*}
    (this holds by $\mathbb{S}^1$-invariance under the first integral);  
    one moment's thought -- use the flow along $Z=\partial_t+\vthet\partial_{\vthet}$ -- gives: 
    $\mathcal{F}_{\varphi(\cdot+a,\cdot)}=\mathcal{F}_{\varphi}(\cdot+a)$ for all $a\in \R$. 
    Now formula \eqref{eqn_ibp1} -- or rather its analogue on $\Delta_s\times (Y\backslash E)$ -- can then be rewritten as: 
     \begin{align*}
      0 =\int_{\Delta_s\times (Y\backslash E)} e^tv_{\varphi}\lich_{\varphi}&(v_{\varphi})\vol^{\varphi}                                              \\
       &=\int_{\Delta_s\times (Y\backslash E)} \big\langle d(e^tv_{\varphi}),\delta^{\varphi}\mathcal{D}_{\varphi}(v_{\varphi})\big\rangle_{\varphi}\vol^{\varphi}
         +\mathcal{F}_{\varphi}(s)-\mathcal{F}_{\varphi}(-s). 
     \end{align*}
    Similarly, 
     \begin{align*}
      \int_{\Delta_s\times (Y\backslash E)} \big\langle d(e^tv_{\varphi}),&\delta^{\varphi}\mathcal{D}_{\varphi}(v_{\varphi})\big\rangle_{\varphi}\vol^{\varphi}            \\
        &= \int_{\Delta_s\times (Y\backslash E)} \big\langle\mathcal{D}_{\varphi}(e^tv_{\varphi}),\mathcal{D}_{\varphi}(v_{\varphi})\big\rangle_{\varphi}\vol^{\varphi}
          +\mathcal{G}_{\varphi}(s)-\mathcal{G}_{\varphi}(-s),
     \end{align*}
    if $\mathcal{G}_{\varphi}$ denotes on $\R$ the function
      \begin{align*}
       s\longmapsto \int_{\{t=s\}}\mathfrak{d}_s^{\varphi} \wedge 2d\vthet  
                   = 4\pi \int_{\{s\}\times(\YE)}\mathfrak{d}_s^{\varphi},                                                                                                                                       \\ 
      \end{align*} 
    where 
      \begin{align*}   \label{eqn_gothd}
       \mathfrak{d}_t^{\varphi}
                :=& e^{-t}\Big[\mathcal{D}_{\varphi}(v_{\varphi})\big(\nabla^{\varphi}(e^tv_{\varphi}),\partial_t\big)
                              \frac{(\omega_{s}^{\varphi})^{m-1}}{(m-1)!}                                                                    \\
                  &\qquad \quad    -\mathcal{D}_{\varphi}(v_{\varphi})\big(\nabla^{\varphi}(e^tv_{\varphi}),\cdot\big)|_Y\wedge d^c_Y v_{\varphi}
                              \wedge \frac{(\omega_{s}^{\varphi})^{m-2}}{(m-2)!}\Big]                                                        \\
                =& \mathcal{D}_{\varphi}(v_{\varphi})\big((v_{\varphi}\nabla^{\varphi}t+\nabla^{\varphi}v_{\varphi}),\partial_t\big)
                              \frac{(\omega_{s}^{\varphi})^{m-1}}{(m-1)!} \\
                 &\qquad \quad -\mathcal{D}_{\varphi}(v_{\varphi})\big(v_{\varphi}\nabla^{\varphi}t+\nabla^{\varphi}v_{\varphi},\cdot\big)|_Y
                              \wedge d^c_Y v_{\varphi}\wedge \frac{(\omega_{s}^{\varphi})^{m-2}}{(m-2)!}
                  \qquad \text{ on }\YE.
      \end{align*}
    This way $\mathcal{G}$ satisfies the translation property $\mathcal{G}_{\varphi(\cdot+a, \cdot)}=\mathcal{G}_{\varphi}(\cdot+a)$, 
    and the exact formulation of \eqref{eqn_ibpeq2} is: 
     \begin{equation*}   
      \begin{aligned}
       \int_{\Delta_s\times (Y\backslash E) } e^t |\mathcal{D}_{\varphi}(v_{\varphi})|_{\varphi}^2\vol^{\varphi}                               
                  = -\frac{1}{2}\bigg(&\int_{\Delta_s\times (Y\backslash E) }
                   v_{\varphi}e^t\big\langle(dt^2)^-,\mathcal{D}_{\varphi}(v_{\varphi})\big\rangle_{\varphi}\vol^{\varphi}                          \\ 
                   &\qquad +\mathcal{F}_{\varphi}(s)+\mathcal{G}_{\varphi}(s)-\mathcal{F}_{\varphi}(-s)-\mathcal{G}_{\varphi}(-s) \bigg).
      \end{aligned}
     \end{equation*}
    Next, in our explicit computation of 
      $\int_{\Delta_s\times (Y\backslash E) } v_{\varphi}e^t\big\langle(dt^2)^-,\mathcal{D}_{\varphi}(v_{\varphi})\big\rangle_{\varphi}\vol^{\varphi}$, 
    the only integration by parts occurs when transforming 
      $\frac{1}{2}\int_{\Delta_s\times Y }e^tv_{\varphi}^2(\nabla^{\varphi}t)\cdot|dt|^2_{\varphi}\vol^{\varphi}$ 
    into $-\int_{\Delta_s\times Y }  e^tv_{\varphi}|dt|^2_{\varphi}dt(\nabla^{\varphi}v_{\varphi})\vol^{\varphi}$; 
    the resulting boundary term is $\mathcal{H}_{\varphi}(s)- \mathcal{H}_{\varphi}(-s)$, 
    where $\mathcal{H}_{\varphi}(s):= -4\pi\int_{\{s\}\times(Y\backslash E)}v_{\varphi}^2|dt|^2_{\varphi}\frac{\omega_{\varphi,s}^{m-1}}{(m-1)!}$ 
    -- again, $\mathcal{H}_{\varphi(\cdot+a, \cdot)}=\mathcal{H}_{\varphi}(\cdot+a)$  --, 
    and thus 
     \begin{equation*}
       \int_{\Delta_s\times (Y\backslash E) }
             e^t |\mathcal{D}_{\varphi}(v_{\varphi})|_{\varphi}^2\vol^{\varphi}
        = \mathcal{K}_{\varphi}(s)- \mathcal{K}_{\varphi}(-s),
      \end{equation*}
    with $\mathcal{K}_{\varphi}=-\frac{1}{2}(\mathcal{F}_{\varphi}+\mathcal{G}_{\varphi}+\mathcal{H}_{\varphi})$, 
    and more generally 
      $\int_{\Delta_{a,b}\times (Y\backslash E) }
             e^t |\mathcal{D}_{\varphi}(v_{\varphi})|_{\varphi}^2\vol^{\varphi}
        = \mathcal{K}_{\varphi}(a)- \mathcal{K}_{\varphi}(b)$
    for all $a\geq b$, if $\Delta_{a,b}:=\{b\leq t\leq a\}$; 
    $\mathcal{K}_{\varphi}$ is thus non-decreasing, and has limits at $\pm\infty$. 
    We will thus be done if we prove that these limits are identical.  
    
    Let us consider any increasing sequence $(t_j)_{j\geq0}$ going to $+\infty$;   
    we also assume 
    that $(t_{j+1}-t_j)_{j\geq0}$ increases to $+\infty$. Set $\alpha_j=\frac{t_{j+1}-t_j}{2}$,  
    and denote by $\varphi_{j}$ the function $\varphi\big(\cdot+t_{j+1}-\alpha_j\big)$. 
    Then $\scal(\omega_{\varphi_j})$ is constant, equal to $\scal(\omega_{\varphi})$, and thus
     \begin{equation*}
       \int_{\Delta_{\alpha_j}\times (Y\backslash E) }
             e^t |\mathcal{D}_{\varphi_j}(v_{\varphi_j})|_{\varphi_j}^2\vol^{\varphi_j}
        = \mathcal{K}_{\varphi_j}(\alpha_j)- \mathcal{K}_{\varphi_j}(-\alpha_j)=\mathcal{K}_{\varphi}(t_{j+1})-\mathcal{K}_{\varphi}(t_j).
     \end{equation*}
    Up to considering a subsequence, 
    one can assume that $(\varphi_j)$, which is uniformly dominated by $\mathfrak{u}_Y$, 
    and has uniformly bounded derivatives at any positive order on $\Delta^*\times(Y\backslash E)$ as $\varphi$ does, 
    converges to some $\varphi_{\infty}\in \mathscr{K}(\omega_0)$ in $C^{\infty}$ on every compact subset of $\Delta^*\times(Y\backslash E)$. 
    And as for $\varphi$, as $\scal(\omega_{\varphi_{\infty}})\equiv\sbar$, 
    $\int_{\Delta^*\times(Y\backslash E)}e^t|\mathcal{D}_{\varphi_{\infty}}(v_{\varphi_{\infty}})|_{\varphi_{\infty}}^2\vol^{\varphi_{\infty}}$ 
    is then finite, 
    and given some compact exhaustive sequence $(K_j)_{j\geq0}$ of $Y\backslash E$, 
    is thus equal to 
    $\lim_{j\to\infty}\int_{\Delta_{\alpha_j}\times K_j}e^t|\mathcal{D}_{\varphi_{\infty}}(v_{\varphi_{\infty}})|_{\varphi_{\infty}}^2\vol^{\varphi_{\infty}}$. 

    Now, 
    $\int_{\Delta_{\alpha_j}\times K_j}e^t|\mathcal{D}_{\varphi_{\infty}}(v_{\varphi_{\infty}})|^2_{\varphi_{\infty}}\vol^{\varphi_{\infty}}
     = \lim_{k\to\infty}\int_{\Delta_{\alpha_j}\times K_j}e^t|\mathcal{D}_{\varphi_{k}}(v_{\varphi_{k}})|^2_{\varphi_{k}}\vol^{\varphi_{k}}$ 
    for any fixed $j$. 
    But 
     $\int_{\Delta_{\alpha_j}\times K_j}e^t|\mathcal{D}_{\varphi_{k}}(v_{\varphi_{k}})|_{\varphi_{k}}^2\vol^{\varphi_{k}}
      \leq \int_{\Delta_{\alpha_k}\times (Y\backslash E)}e^t|\mathcal{D}_{\varphi_{k}}(v_{\varphi_{k}})|_{\varphi_{k}}^2\vol^{\varphi_{k}}
      =\mathcal{K}_{\varphi}(t_{k+1})-\mathcal{K}_{\varphi}(t_k)$ for $k\geq j$, 
    and therefore $\int_{\Delta_{\alpha_j}\times K_j}e^t|\mathcal{D}_{\varphi_{\infty}}(v_{\varphi_{\infty}})|^2\vol^{\varphi_{\infty}}=0$, 
    as $\mathcal{K}_{\varphi}(t_k)$ tends to $\lim_{+\infty}\mathcal{K}_{\varphi}$ as $k$ goes to $\infty$. 

    We can thus conclude that $\int_{\Delta^*\times(Y\backslash E)}e^t|\mathcal{D}_{\varphi_{\infty}}(v_{\varphi_{\infty}})|^2\vol^{\varphi_{\infty}}=0$, 
    which is equivalent to: $\mathcal{D}_{\varphi_{\infty}}(v_{\varphi_{\infty}})$ vanishes identically, 
    or: $\nabla^{\varphi_{\infty}}v_{\varphi_{\infty}}$ is real holomorphic. 
    Since $\nabla^{\varphi_{\infty}}v_{\varphi_{\infty}}=\tfrac{\partial}{\partial t}+\nabla^{\varphi_{\infty}}t$, 
    it has shape $\gamma\tfrac{\partial}{\partial t}+Z$ with $Z$ tangent to $Y\backslash E$, 
    as $d\vartheta(\nabla^{\varphi_{\infty}}t)=\langle dt,d\vartheta\rangle_{\varphi_{\infty}}=0$. 
    Moreover 
    $\gamma=dt(\nabla^{\varphi_{\infty}}v_{\varphi_{\infty}})
           =1+v_{\varphi_{\infty}}|dt|_{\varphi_{\infty}}^2$;   
    furthermore, one has: 
     \begin{lem}  \label{lem_holovf}
      Let $Z$ be a real holomorphic vector field on $\Delta^*\times(Y\backslash E)$, 
      bounded up to order 1 for $\omega_0$. Then $Z$ is tangent to $Y\backslash E$, and constant along $\Delta^*$. 
     \end{lem}
    The proof of this lemma is postponed after the current proof. 
    For now we get, as $\nabla^{\varphi_{\infty}}v_{\varphi_{\infty}}$ is bounded at any order with respect to $\omega_0$, 
    that $\gamma\equiv0$, i.e.
    $\dot{\varphi_{\infty}}-1=v_{\varphi_{\infty}}=-|dt|_{\varphi_{\infty}}^{-2}
                     =-\big(1+\ddot{\varphi_{\infty}}-\dot{\varphi_{\infty}}-|d_Y\dot{\varphi_{\infty}}|_{\infty,t}^2\big)$, 
   or: $\ddot{\varphi_{\infty}}=|d_Y\dot{\varphi_{\infty}}|_{\infty,t}^2\geq0$. 
   Since $\varphi_{\infty}(\cdot,y)$ is bounded for all $y\in Y\backslash E$, this implies that $\varphi_{\infty}$ is constant in the $\Delta^*$-direction, 
   thus $\dot{\varphi_{\infty}}=0$, and in particular $\ddot{\varphi_{\infty}}=|d_Y\dot{\varphi_{\infty}}|_{\infty,t}^2=0$. 
   In other words, $\varphi_{\infty}$ is a function on $(Y\backslash E)$, $\psi_{\infty}$ say, 
   \textit{independent of $t$}.  

   We interpret this by saying that $v_{\varphi_j}$ converges to $-1$, 
   and that $\omega_{\varphi_j}$ converges to 
   $\omega_{\varphi_{\infty}}=dt\wedge 2e^{-t}d\vthet+\omega_Y^{\psi_{\infty}}$, $\omega_Y^{\psi_{\infty}}=\omega_{\YE}+dd^c_Y\psi_{\infty}$, 
   in $C^{\infty}_{\loc}$ topology. 
   Hence by dominated convergence, $\mathcal{K}_{\varphi}\big(\frac{t_{j+1}+t_j}{2}\big)=\mathcal{K}_{\varphi_j}(0)$ 
   tends to $\mathcal{K}_{\varphi_{\infty}}(0)$, and, as $v_{\varphi_{\infty}}=-1$,
    \begin{align*}
     \mathcal{K}_{\varphi_{\infty}}(0)
       &=  -2\pi \int_{\{0\}\times (Y\backslash E)}v_{\varphi_{\infty}}\mathfrak{D}_0^{\varphi_{\infty}}
              +\mathfrak{d}_0^{\varphi_{\infty}}
              -v_{\varphi_{\infty}}^2|dt|^2_{\varphi_{\infty}}\frac{(\omega_{Y}^{\psi_{\infty}})^{m-1}}{(m-1)!} \\
       &=-2\pi \int_{\{0\}\times(Y\backslash E)}0+0-(-1)^2\cdot1 \, \frac{(\omega_{Y}^{\psi_{\infty}})^{m-1}}{(m-1)!}                               
        =2\pi\vl(Y\backslash E),
    \end{align*}
   as $|dt|^2_{\varphi_{\infty}}=1$, since $\omega_{\varphi_{\infty}}=dt\wedge 2e^{-t}d\vthet+\omega_{Y}^{\psi_{\infty}}$.  
   We recall that the volume of $Y\backslash E$, 
   even if computed with respect to $\omega_{Y}^{\psi_{\infty}}$, 
   depends only on $[\omega_Y]$. 
   On the other hand, $\mathcal{K}_{\varphi}\big(\frac{t_{j+1}+t_j}{2}\big)$ converges to $\lim_{+\infty}\mathcal{K}_{\varphi}$; 
   we hence get: $\lim_{+\infty}\mathcal{K}_{\varphi}=2\pi\vl(Y\backslash E)$. 

   These arguments apply symmetrically, and thus $\lim_{-\infty}\mathcal{K}_{\varphi}=2\pi\vl(Y\backslash E)$. 
   Therefore $\lim_{-\infty}\mathcal{K}_{\varphi}=\lim_{+\infty}\mathcal{K}_{\varphi}$, 
   and finally $\int_{\Delta^*\times(Y\backslash E)}e^t|\mathcal{D}_{\varphi}(v_{\varphi})|_{\varphi}^2\vol^{\varphi}=0$. 
  \cqfd
 
~
  
 \noindent
 \textit{End of proof of Theorem \ref{thm_splttng}. --- } 
  We have: $\mathcal{D}_{\varphi}(v_{\varphi})\equiv 0$; 
 as seen in the above proof for $\varphi_{\infty}$, this implies that $\varphi$ \textit{does not depend on }$t$: 
 Theorem \ref{thm_splttng} is proved. 
 \cqfd
      
 ~
 
 \noindent
 \prf \textit{of Lemma \ref{lem_holovf}. ---}   
  Any holomorphic function $f$ on $\Delta^*$ which is $\mathcal{O}(1-|z|)$ near $\partial\Delta$, 
  and $\mathcal{O}\big(|z|\big|\log(|z|)\big|\big)$ near $0$, vanishes identically on $\Delta^*$.  
  Indeed, extend $f$ through 0, and for $r\in(0,1)$, $\vareps>0$, pick $\delta\in(0,1-r)$ so that $|f|\leq\vareps$ 
  on $\partial\Delta_{1-\delta}$. 
  By the maximum principle, $|f|\leq\vareps$ on $\Delta_r$; since $r$ and $\vareps$ are arbitrary, $f\equiv0$. 
    
  Now given $Z$ as in the statement, take a open subset $U$ of holomorphic coordinates $(z^2,\dots,z^m)$ on $Y\backslash E$, 
  and write $Z^{1,0}=f\frac{\partial}{\partial z}+\sum_{j=2}^mf_j\frac{\partial}{\partial z^m}$. Fix $x\in U$. 
  As $Z$ is bounded, $f(\cdot,x):\Delta^*\to\C$, which is holomorphic, 
  is $\mathcal{O}\big(|z|\big|\log(|z|)\big|\big)$ near $0$, and $\mathcal{O}(1-|z|)$ near $\partial\Delta$, and thus $f(\cdot,x)\equiv0$. 
  The same holds for the $\frac{\partial f_j(\cdot,x)}{\partial z}$, 
  and as $\frac{\partial f_j(\cdot,x)}{\partial \overline{z}}\equiv0$, 
  the $f_j$ are constant along $\Delta^*$. 
 \cqfd

\section{Parametrisation of Kähler metrics of almost constant scalar curvature}   \label{sctn_almcscKfamily}

 \subsection{Family of Kähler metrics of almost constant scalar curvature}
 
\textit{Definition. ---} 
We consider a compact Kähler manifold $(Y,\omega_Y,J_Y)$ of dimension $n\geq 1$, 
and define: 
 \begin{df}
  Let $(\omega_t)_{t\geq0}$ be a smooth family of Kähler metrics in $[\omega_Y]$ such that:
  \begin{enumerate}
   \item $\scal(\omega_t)$ converges at any order to 
        $\overline{\scal}_Y:=-4\pi n\frac{\chern_1(K_Y)\cdot [\omega_Y]^{n-1}}{[\omega_Y]^{n}}$, i.e.  
        for any $\kappa\geq0$, $\scal(\omega_t)\to \overline{\scal}_Y$
        and for any positive $\ell$, $\partial^{\ell}_t\scal(\omega_t)\to 0$, in $C^{\kappa}(Y)$, as $t$ goes to $\infty$; 
   \item $(\omega_t)_{t\geq0}$ is bounded in $C^{\kappa}$ for any $\kappa$, 
         and there is some positive constant $c$ such that for all $t\geq 0$, $\omega_t\geq c\omega_Y$;  
  \end{enumerate} 
  we then say that $(\omega_t)_{t\geq0}$ is a \emph{family of Kähler metrics of almost constant scalar curvature}. 
  
  We say moreover that such a family has \emph{extinguishing variation} if for all positive $\ell$, 
  $\partial_t^{\ell}(\omega_t)$ tends to 0 in all $C^{\kappa}(Y)$ as $t$ goes to $\infty$. 
 \end{df}
 
In this definition, we assume of course that all the metrics are Kähler with respect to the \textit{fixed} complex structure $J_Y$.  
Notice moreover the existence of a family $(\omega_t)$ of almost constant scalar curvature in $[\omega_Y]$, 
implies that of a constant scalar curvature metric in this class:  
take any $C^{\infty}$-limit point of $(\omega_t)$.

 ~

\noindent
\textit{Basic example. --- } 
Assuming that $\omega_Y$ has constant scalar curvature, and is the unique such metric in its Kähler class 
(as is the case when for instance $Y$ has no non-trivial holomorphic vector fields \cite{ChTi}), 
it is rather straightforward to see that a family of almost constant scalar curvature tends to $\omega_Y$ in $C^{\infty}$-topology. 
When there exist several constant scalar curvature metrics in a same Kähler class,  
the situation might be more delicate, even with extinguishing variation, as the following example illustrates: 
consider a smooth family of holomorphic vector fields $(Z_t)$ such that for all 
$\ell\geq0$, $\partial_t^{\ell} Z_t$ tends to 0 in $C^{\infty}$-topology. 
Assume moreover that the family of automorphisms $\sigma_t$ of $Y$ such that $\partial_t\sigma_t=Z_t$ for all $t$ 
is bounded, in the sense that $\sigma_t^*\omega_Y$ stays bounded below, and bounded in $C^{\kappa}$ for all $\kappa$. 
Then $(\omega_t):=(\sigma_t^*\omega_Y)$ is clearly a family of metrics of (almost) constant scalar curvature, 
since for all $t$, $\scal(\omega_t)=\sbar_Y$.  
On the other hand, we arrange easily 
so that $(\omega_t)$ does not converge:  
take $Z$ a vector field such that $\mathcal{L}_Z\omega_Y\neq 0$ 
-- this exists as soon as the uniqueness for constant scalar curvature metrics fails --, 
and set $Z_t=f(t)Z$ with $f$ a smooth function of $t$ 
tending to 0 at any order, but with $\int_{s=0}^t f$ keeping oscillating between two fixed values, 
e.g. $f(t)=\frac{\cos[\log(1+t)]}{1+t}$.

 \subsection{Parametrisation}

The following proposition, which is the technical core of this part, 
tells us that the previous example is asymptotically the only possible type of situation for 
a family of almost constant scalar curvature with extinguishing variation: 
 \begin{prop}   \label{prop_asympparam}
  Let $(\omega_t)$ be a family of almost constant scalar curvature in $[\omega_Y]$, 
  with extinguishing variation. 
  Then there exists 
  a smooth family of constant scalar curvature Kähler metrics $(\varpi_t)$ 
  in $[\omega_Y]$ such that:
   \begin{itemize}
    \item $\omega_t-\varpi_t=o(1)$ in $C^{\kappa}(Y)$ as $t$ goes to $\infty$, for all $\kappa\geq0$; 
          in particular, $(\varpi_t)_{t\geq 0}$ is uniformly bounded below, 
          i.e. $\varpi_t\geq c_1\omega_Y$ for some positive $c_1$; 
    \item there exists a smooth family of holomorphic vector fields $(\mathsf{Z}_t)$ such that if $(\sigma_t)_{t}$ 
          is the associated flow, then $\varpi_t=\sigma_t^*\varpi_0$ for all $t$, 
          and for all $\ell\geq0$, $\partial_t^{\ell}\mathsf{Z}_t=o(1)$ in $C^{\kappa}(Y)$ as $t$ goes to $\infty$, for all $\kappa\geq0$; 
    \item $\varpi_0$ can be taken as any limit point of $(\omega_t)$ in $C^{\infty}(Y)$. 
   \end{itemize}

 \end{prop}
 
~
 
In other terms, 
a family of almost constant scalar curvature with extinguishing variation 
can be parametrised, up to a small error in $C^{\infty}$-topology, as the pull-back of a \textit{fixed} constant scalar curvature metric 
by some automorphism flow with asymptotically vanishing time derivatives of positive order. 

~

\prf \textit{of Proposition \ref{prop_asympparam}. }--- 
 We consider a family $(\omega_t)$ 
 as in the statement of the proposition. 
 The following strategy will guide us: 
  \begin{enumerate}
   \item we fix $\kappa\geq 2$ and $\alpha\in(0,1)$, call $\mathscr{L}$ the set of limit points of $(\omega_t)$ in $C^{\kappa,\alpha}$-topology, 
         show that $\mathscr{L}$ is a nonempty set of smooth metrics with constant scalar curvature, 
         and is actually the set of $C^{\infty}$-limit points of $(\omega_t)$; 
         in particular, $\mathscr{L}$ does not depend on $\kappa$ nor on $\alpha$, hence is bounded at any order;
   \item similarly,  we see  that $\mathscr{L}$ is connected in $C^{\kappa,\alpha}$, hence in $C^{\infty}$, topologies; 
   \item we use a result due to Calabi to parametrise a well-chosen sequence of metrics in $\mathscr{L}$ as inductive pull-backs of some arbitrary fixed metric 
         by ``small'' reduced automorphisms of $Y$, 
         interpolate this sequence into a family of constant scalar curvature metrics, and compare it to $(\omega_t)$. 
  \end{enumerate}
  
 Points 1 and 2 require very usual arguments, which we include for the sake of completeness;  
 Point 3, less standard, needs a little more care. 

 ~
 
 \noindent
 \textit{Point 1: Regularity of elements of $\mathscr{L}$, and independence from $(\kappa,\alpha)$.} 
 The family $(\omega_t)$ is bounded in, say, $C^{\kappa,\beta}$ with $\beta\in (\alpha,1)$, 
 and therefore admits a subsequence converging in $C^{\kappa,\alpha}$: $\mathscr{L}\neq \varnothing$. 
 
 Let us consider some $\varpi_{\infty}\in \mathscr{L}$, 
 which is the $C^{\kappa,\alpha}$-limit of a subsequence $(\omega_{t_j})$ where we can assume that $(t_j)$ increases to $\infty$. 
 As $\omega_{t_j}\geq c\omega_Y$ for all $j$, $\varpi_{\infty}$ is positively bounded below, and is therefore a metric;  
 it moreover has constant scalar curvature, since along our subsequence, 
 the scalar curvatures converge to $\scal(\varpi_{\infty})$ in $C^{\kappa-2,\beta}$, 
 while $\scal(\omega_t)$ converges to $\sbar_Y$ as $t$ goes to $\infty$: $\scal(\varpi_{\infty})=\sbar_Y$. 
 
 Given moreover any $\lambda\geq \kappa+1$, $(\omega_{t_j})$ is bounded in $C^{\lambda,\alpha}$, 
 and thus admits a subsequence $(\omega_{t_{j_{\ell}}})$ converging in $C^{\lambda}$, 
 necessarily to $\varpi_{\infty}$, as $C^{\lambda}$-convergence implies $C^{\kappa,\alpha}$-convergence. 
 Hence $\varpi_{\infty}$, which is thus indeed $C^{\lambda}$, is a $C^{\lambda}$-limit point of $(\omega_t)$, 
 and this holds for all $\lambda\geq\kappa+1$: 
 $\varpi_{\infty}$ is smooth (this is also deducible from $\scal(\varpi_{\infty})=\sbar_Y$), 
 and is a $C^{\infty}$-limit point of $(\omega_t)$. 
This settles Point 1. 
 
 ~
 
 \noindent
 \textit{Point 2: Connectedness of $\mathscr{L}$.}
 The connectedness assertion on $\mathscr{L}$ endowed with the induced $C^{\kappa,\alpha}$-topology 
 can be viewed as the analogue of the similar statement on the set of limit points of 
 a sequence in a compact metric space such that the distance between two consecutive terms goes to 0. 
 Now the bound on $(\omega_t)$ in $C^{\kappa, \alpha}$  
 does not provide compactness in general; 
 we nonetheless bypass this lack thanks to higher order bounds, and to the (almost) constant scalar curvature property.  
 
 Let us thus assume that $\mathscr{L}=\mathscr{L}_0\sqcup\mathscr{L}_1$, 
 with the $\mathscr{L}_i$ non-empty closed subsets of $\mathscr{L}$, for the $C^{\kappa,\alpha}$-distance 
 ${\rm d}_{C^{\kappa,\alpha}}$; 
 notice that as a set of limit points, $\mathscr{L}$ is closed for this distance in the set of $C^{\kappa,\alpha}$ metrics, 
 hence so are the $\mathscr{L}_i$. 
 And as $\mathscr{L}$ is bounded at any higher order, so are the $\mathscr{L}_i$. 
 We claim that there exist $\varpi_{\infty}^i\in\mathscr{L}_i$, $i=0,1$, 
 such that 
 $\vareps:={\rm d}_{C^{\kappa,\alpha}}(\mathscr{L}_0,\mathscr{L}_1)={\rm d}_{C^{\kappa,\alpha}}(\varpi_{\infty}^0,\varpi_{\infty}^1) $, 
 distance which is thus $>0$. 
 Consider indeed two sequences $(\varpi_{j}^i)_{j\geq 0}$ of elements of $\mathscr{L}_i$, $i=0,1$, 
 such that ${\rm d}_{C^{\kappa,\alpha}}(\varpi_{j}^0,\varpi_{j}^1)\to \vareps$ as $j\to \infty$. 
 Then as the $(\varpi_{j}^i)_{j\geq 0}$ are bounded in $C^{\kappa,\beta}$ for any $\beta\in (\alpha,1)$, 
 they admit subsequences $(\varpi_{j_k}^i)_{k\geq 0}$, $(j_k)$ independent of $i\in\{0,1\}$,  
 converging in $C^{\kappa,\alpha}$ to respective limits $\varpi_{\infty}^i$; 
 moreover, $\varpi_{\infty}^i\in\mathscr{L}_i$ closed in $C^{\kappa,\alpha}$ topology, 
 and $\vareps= \lim_k {\rm d}_{C^{\kappa,\alpha}}(\omega_{j_k}^0,\omega_{j_k}^1)= {\rm d}_{C^{\kappa,\alpha}}(\varpi_{\infty}^0,\varpi_{\infty}^1)$. 
 
 We consider now an increasing sequence of ``times'' $(t_j)$ admitting two disjoint subsequences 
 $(t_{j_{k}^0})$ and $(t_{j_{k}^1})$ such that 
 $\omega_{t_{j_{k}^i}}\to \varpi_{\infty}^i$ as $k\to\infty$ in $C^{\kappa,\alpha}$. 
 Up to adding intermediate times, we can furthermore assume that $\|\omega_{t_{j+1}}-\omega_{t_j}\|_{C^{\kappa,\alpha}}\to0$ as $j\to\infty$, 
 since $\|\partial_t\omega_t\|_{C^{\kappa,\alpha}}$ is bounded.  
 Let us show that an infinite number of $\omega_{t_j}$ are a distance at least $\vareps/3$ from both $\mathscr{L}_0$ and $\mathscr{L}_1$. 
 
 Take $J$ large enough so that $\|\omega_{t_{j+1}}-\omega_{t_j}\|_{C^{\kappa,\alpha}} < \frac{\vareps}{3}$; 
 take also $k,\ell$ large enough so that, if $j_0:=j_{k}^0<j_{\ell}^1=:j_1$, 
 then $\|\omega_{t_{j_0}}-\varpi_{0,\infty}\|_{C^{\kappa,\alpha}}$, $\|\omega_{t_{j_1}}-\varpi_{1,\infty}\|_{C^{\kappa,\alpha}}\leq\frac{\vareps}{3}$. 
 Define $j_2$ to be the smallest integer $\geq j_0$ so that $\omega_{t_{j_2}}$ is at distance $\geq\frac{\vareps}{3}$ from $\mathscr{L}_0$; 
 $j_2$ is well-defined, and $j_2\in\{j_0+1,\dots,j_1\}$, 
 since $\omega_{t_{j_1}}$ is at least at distance $\frac{2\vareps}{3}$ from $\mathscr{L}_0$. 
 One also has ${\rm d}_{C^{\kappa,\alpha}}(\omega_{t_{j_2}},\mathscr{L}_1)
                \geq  {\rm d}_{C^{\kappa,\alpha}}(\mathscr{L}_0,\mathscr{L}_1) 
                     - {\rm d}_{C^{\kappa,\alpha}}(\omega_{t_{j_2}},\omega_{t_{j_2-1}})
                     - {\rm d}_{C^{\kappa,\alpha}}(\omega_{t_{j_2-1}},\mathscr{L}_1)
                \geq \frac{\vareps}{3}$. 
                
 We can repeat this argument with $j_0$ and $j_1$ as large as wanted, 
 and thus end up with a sequence $(t_{j_{k}})_{k\geq0}$ going to $\infty$ such that for any $k$,  
 $\omega_{t_{j_{k}}}$ is at distance at least $\frac{\vareps}{3}$ from both $\mathscr{L}_i$. 
 Now $(\omega_{t_{j_{k}}})$ is bounded in $C^{\kappa,\beta}$, 
 hence admits a subsequence converging in $C^{\kappa,\alpha}$ to some $\varpi_{\infty}$, 
 necessarily at distance at least $\frac{\vareps}{3}$ from the $\mathscr{L}_i$. 
 But by definition, $\varpi_{\infty}\in \mathscr{L}$, hence a contradiction with the assumption $\mathscr{L}=\mathscr{L}_0\cup\mathscr{L}_1$. 
 
 We are left with the connectedness assertion in $C^{\infty}$-topology. 
 We actually settle this by the more general statement that the induced $C^{\kappa,\alpha}$ and $C^{\infty}$-topologies coincide on $\mathscr{L}$. 
 As the elements of $\mathscr{L}$ are smooth constant scalar curvature metrics 
 with uniform lower bound and bounds at any order (coming from such bounds on $(\omega_t)$), 
 this statement merely comes from the iterated observation that given smooth $\varpi$ and 
 $\tilde{\varpi}$,  
 then $\|\varpi-\widetilde{\varpi}\|_{C^{\lambda,\alpha}}$ 
 is bounded by $C_{\lambda,\alpha}\big(\|\varpi-\widetilde{\varpi}\|_{C^{\kappa,\alpha}}
                                           +\|\scal(\varpi)-\scal(\widetilde{\varpi})\|_{C^{\lambda-2,\alpha}}\big)$, 
 where $C_{\lambda,\alpha}$ depends only on lower bounds and $C^{\lambda,\alpha}$-bounds on $\varpi$ and $\widetilde{\varpi}$. 
 Let us detail how this goes for $\lambda=\kappa+1$. 
 Let $\psi$ so that $\widetilde{\varpi}=\varpi+dd^c\psi$, normalized by $\int_Y\psi\,\varpi^n=0$; 
 this way $\|\psi\|_{C^{\kappa+2,\alpha}}\leq C\|\widetilde{\varpi}-\varpi\|_{C^{\kappa,\alpha}}$ with $C$ as announced. 
 Moreover, in local coordinates, if $g$ is the metric $\varpi(\cdot,J_Y\cdot)$ and $\tilde{g}$ is $\widetilde{\varpi}(\cdot,J_Y\cdot)$,
  \begin{align*}
   \scal(\widetilde{\varpi})-\scal(\varpi)
                          =& -g^{p \bar{q}}\big(g^{r \bar{s}}\partial_p\partial_{\bar{q}}(\tilde{g}_{r\bar{s}}+\partial_r\partial_{\bar{s}}\psi)
                                        +(\partial_p g^{r \bar{s}})\partial_{\bar{q}}(\tilde{g}_{r\bar{s}}+\partial_r\partial_{\bar{s}}\psi)\big)                   \\
                           & +\tilde{g}^{p \bar{q}}\big(\tilde{g}^{r \bar{s}}\partial_p\partial_{\bar{q}}(\tilde{g}_{r\bar{s}})
                                        +(\partial_p \tilde{g}^{r \bar{s}})\partial_{\bar{q}}(\tilde{g}_{r\bar{s}}\big)                                                  \\    
                          =& -g^{p \bar{q}}g^{r \bar{s}}\partial_{p\bar{q}r\bar{s}}\psi 
                             -g^{p \bar{q}}(\partial_p g^{r \bar{s}})\partial_{\bar{q}r\bar{s}}\psi                                                            \\
                           & -(\partial_{p\bar{q}} \tilde{g}_{r\bar{s}})(g^{p \bar{q}}g^{r \bar{s}}-\tilde{g}^{p \bar{q}}\tilde{g}^{r \bar{s}})                
                             -(\partial_{\bar{q}} \tilde{g}_{r\bar{s}})(g^{p \bar{q}}\partial_pg^{r \bar{s}}-\tilde{g}^{p \bar{q}}\partial_p\tilde{g}^{r \bar{s}})
  \end{align*}
 We rewrite the latter equation as: 
  \begin{equation}   \label{eqn_eta}
   g^{p \bar{q}}g^{r \bar{s}}\partial_{p\bar{q}r\bar{s}}\psi 
    = \eta   -\big(\scal(\varpi)-\scal(\widetilde{\varpi})\big)
  \end{equation}
 with 
  $\eta = -g^{p \bar{q}}(\partial_p g^{r \bar{s}})\partial_{\bar{q}r\bar{s}}\psi
       -(\partial_{p\bar{q}} \tilde{g}_{r\bar{s}})(g^{p \bar{q}}g^{r \bar{s}}-\tilde{g}^{p \bar{q}}\tilde{g}^{r \bar{s}})                
       -(\partial_{\bar{q}} \tilde{g}_{r\bar{s}})(g^{p \bar{q}}\partial_pg^{r \bar{s}}-\tilde{g}^{p \bar{q}}\partial_p\tilde{g}^{r \bar{s}})$.

 A local $C^{\kappa-1,\alpha}$-bound on $\eta$ now easily follows from a $C^{\kappa+2,\alpha}$-bound on $\psi$, 
 lower bounds plus $C^{\kappa+1,\alpha}$-bounds on $\varpi$ and $\widetilde{\varpi}$, 
 and a $C^{\kappa,\alpha}$-bound on $\varpi-\widetilde{\varpi}$; 
 more precisely, 
 $\|\eta\|_{C^{\kappa-1,\alpha}}\leq C_1\|\psi\|_{C^{\kappa+2,\alpha}}+C_2\|\varpi-\widetilde{\varpi}\|_{C^{\kappa,\alpha}}$ with $C_1,C_2$ 
 as announced, hence $\|\eta\|_{C^{\kappa-1,\alpha}}\leq C\|\varpi-\widetilde{\varpi}\|_{C^{\kappa,\alpha}}$ with $C$ as announced, 
 by the previous control on $\|\psi\|_{C^{\kappa+2,\alpha}}$. 
 The conclusion follows from Schauder elliptic estimates applied to \eqref{eqn_eta}, 
 together with $\|\psi\|_{C^{0}}\leq C\|\varpi-\widetilde{\varpi}\|_{C^{\kappa,\alpha}}$: 
 as $g^{p \bar{q}}g^{r \bar{s}}\partial_{p\bar{q}r\bar{s}}$ is elliptic in the open subset $\Omega$ of work, 
 with lower bounds and $C^{\kappa-1,\alpha}$ bounds on $(g^{p \bar{q}}g^{r \bar{s}})_{pqrs}$ coming from such bounds on $\varpi$, 
 one has $\|\psi\|_{C^{\kappa+3,\alpha}(\Omega')}  \leq C\big(\|\varpi-\widetilde{\varpi}\|_{C^{\kappa,\alpha}(\Omega)}
                                                     +\|\scal(\varpi)-\scal(\widetilde{\varpi})\|_{C^{\kappa-1,\alpha}(\Omega)}\big)$, 
 thus $\|\varpi-\widetilde{\varpi}\|_{C^{\kappa+1,\alpha}(\Omega')}     
                                                   \leq C\big(\|\varpi-\widetilde{\varpi}\|_{C^{\kappa,\alpha}(\Omega)}
                                                     +\|\scal(\varpi)-\scal(\widetilde{\varpi})\|_{C^{\kappa-1,\alpha}(\Omega)}\big)$, 
 with $C$ as announced, on a slightly smaller open subset $\Omega'$, 
 hence the desired estimate as $Y$ is compact. 
 
 ~

 \noindent
 \textit{Point 3: Interpolation of a sample of limit points.} 
  We proceed to the sampling mentioned above according to the following principle: 
  \emph{given $\vareps>0$, there exists $J=J(\vareps)\geq0$ such that 
  for any $j\geq J$, there exists $\varpi_{j}\in \mathscr{L}$ such that $\|\omega_{t=j}-\varpi_{j}\|_{C^{\kappa,\alpha}}<\vareps$}. 
  Indeed, assume this does not hold, and pick $\vareps>0$  
  and a sequence $(j_k)$ of nonnegative integers going to $\infty$, such that for all $k$ and all $\varpi\in\mathscr{L}$, 
  $\|\omega_{t=j_k}-\varpi\|_{C^{\kappa,\alpha}}\geq \vareps$. 
  As $(\omega_{j_k})_k$ is bounded in $C^{\kappa,\beta}$ ($\beta\in(\alpha,1)$), 
  it admits a subsequence converging in $C^{\kappa,\alpha}$ to some $\omega_{\infty}$. 
  Hence for all $\varpi\in \mathscr{L}$, $\|\omega_{\infty}-\varpi\|_{C^{\kappa,\alpha}}\geq \vareps$;  
  this contradicts the assertion $\omega_{\infty}\in\mathscr{L}$ given by the definition of $\mathscr{L}$. 
  
  We fix now a positive sequence $(\vareps_k)_{k\geq0}$ going to 0,  
  and set $j_k=J(\vareps_k)$ for $k\geq 0$, according to our principle; we can assume that for all $k$, 
  $J(\vareps_{k+1})>J(\vareps_k)$. 
  This enables us, for all $k$ and all $j\in\{j_k,\dots,j_{k+1}-1\}$, 
  to pick some $\varpi_j\in\mathscr{L}$ such that $\|\omega_{t=j}-\varpi_j\|\leq \vareps_k$. 
  We thus constitute a sequence $(\varpi_j)$ in $\mathscr{L}$ which is asymptotic to $(\omega_{t=j})$ 
  (in $C^{\kappa,\alpha}$, thus in all $C^{\lambda}$, by the estimate of the previous point). 
  
  In order to interpolate between the $\varpi_j$, 
  we use the following result due to Calabi \cite{cal2}, see also \cite[Prop. 3.3.4]{gau}: 
  the space of extremal Kähler metrics among $[\omega_Y]$ is a submanifold of the space of Kähler metrics in $[\omega_Y]$, 
  and each connected component of this submanifold is an orbit of the reduced automorphism group $H_{\rm red}(Y,J_Y)$. 
  More precisely, the tangent space of this submanifold at each point coincides with the tangent space of the orbit. 
  Here extremal metrics are constant scalar curvature metrics: 
  indeed 
  as seen above, for any $\varpi\in\mathscr{L}\subset[\omega_Y]$, $\scal(\varpi)=\sbar_Y$, 
  which makes the Futaki character of $[\omega_Y]$ vanish, and forces any extremal metric in $[\omega_Y]$ to have constant scalar curvature. 
  Moreover $\mathscr{L}$ is connected, and is thus contained in \textit{one} connected component of 
  the space $\mathscr{S}_{[\omega_Y]}$ of constant scalar curvature metrics of $[\omega_Y]$; 
  its elements can thus all be written as $\sigma^*\omega_Y$, up to choosing $\omega_Y$ in $\mathscr{L}$, say. 
  Here we can be more precise: 
  the above statement tells us that given any $\omega$ of constant scalar curvature, 
  all the constant scalar curvature metrics of its neighbourhood 
  (for any $C^{\kappa,\alpha}$ topology, or even for $C^{\infty}$ topology, 
  as the submanifold has finite dimension) 
  can be written as $\sigma^*\omega$ with $\sigma$ a reduced automorphism close to identity; 
  up to reducing the neighbourhood, these metrics can thus all be written $(\Phi_1^{Z})^*\omega$ with $Z$ a small (real) holomorphic vector field 
  in $J\mathfrak{k}_0(\omega)$ where $\mathfrak{k}_0(\omega)$ is the set of Hamiltonian Killing fields for $\omega$, 
  as in the constant scalar curvature case, $\mathfrak{h}_0:={\rm Lie}\big(H_{\rm red}(Y,J_Y)\big)$ 
  splits as $\mathfrak{k}_0(\omega)\oplus J\mathfrak{k}_0(\omega)$.  
  
  Now $\mathscr{L}$ is closed and bounded in $\mathscr{S}_{[\omega_Y]}$, hence compact,  
  and therefore there exists a function $\eta=\eta(\vareps)$ which is $o(1)$ 
  such that any ball of radius $\vareps\in(0,\bar{\vareps}_0)$ and centre $\varpi\in \mathscr{L}$ 
  is covered by the $(\sigma_Z)^*\varpi$ with $Z\in J\mathfrak{k}_{0}(\varpi)$, $|Z|<\eta(\vareps)$, 
  where $\sigma_Z$ denotes the flow at time 1 of $Z$. 
  
  Without loss of generality, assume $(\vareps_k)$ is decreasing and takes value in $(0,\bar{\vareps}_0/3)$;  
  Assume likewise that $\gamma(t):=\|\partial_t\omega_t\|_{C^{\kappa,\alpha}}<\frac{\bar{\vareps}_0}{3} $ for all $t$. 
  Then for all $k$, and all $j\geq j_k$,
   \begin{align*}
    \|\varpi_{j+1}-\varpi_{j}\|_{C^{\kappa,\alpha}}
                               &\leq \|\varpi_{j+1}-\omega_{j+1}\|_{C^{\kappa,\alpha}}
                                     +\|\varpi_{j}-\omega_{j}\|_{C^{\kappa,\alpha}}+\Big\|\int_{t=j}^{j+1}\partial_t\omega_t\,dt\Big\|_{C^{\kappa,\alpha}} \\
                               &\leq 2\vareps_k+\delta_j <\bar{\vareps}_0, 
   \end{align*}
  where $\delta_j=\int_{t=j}^{j+1}\gamma(t)\,dt$. 
  Hence for all $k$, and all $j\in\{j_k,\dots,j_{k+1}-1\}$, there exists a holomorphic vector field $Z_j$, $|Z_j|\leq \eta(2\vareps_k+\delta_j)$, 
  such that $\varpi_{j+1}=\sigma_{Z_j}^*\varpi_{j}$. 
  
  We conclude as follows: 
  let $\chi$ be a nonnegative smooth function with compact support in $(0,1)$ and with integral $1$. 
  We define the time-dependent holomorphic vector field $\mathsf{Z}_t:=\sum_{j=j_0}^{\infty}\chi(t-j)Z_j$, 
  and the associated flow $\sigma$. 
  We set $\varpi_t:= \sigma_{t-j_0}^*\varpi_{j_0}$ for all $t$; 
  these are metrics of constant scalar curvature, in some fixed neighbourhood of $\mathscr{L}$. 
  One easily checks that indeed $\varpi_{t=j}=\varpi_j$ for all $j\geq j_0$, 
  as well as the claimed asymptotic properties of $(\varpi_t)_t$. 
  For instance, for $t\geq j_0$,
  $\partial_t\varpi_t= \mathcal{L}_{\mathsf{Z}_t}\varpi_t=\chi(t-[t])\mathcal{L}_{Z_{[t]}}\varpi_{t}
                     =\chi(t-[t])d\big(\varpi_{t}(Z_{[t]},\cdot)\big)$; 
  since $\varpi_t$ is bounded at any order, $\|\partial_t\varpi_t\|_{C^{\kappa,\alpha}}$ is controlled by $|Z_{[t]}|$, 
  which is controlled by $\eta(2\vareps_{k([t])}+\delta_{[t]})$ (where $k(\cdot)$ is defined by $\ell\in\{j_{k(\ell)},\dots,j_{k(\ell)-1}\}$), 
  hence goes to $0$ as $t$ goes to $\infty$. 
  Finally, as $\mathsf{Z}_t=0$ for $t\in[0,j_0]$, $\varpi_0=\varpi_{j_0}\in \mathscr{L}$; 
  now $\varpi_{j_0}$ was arbitrarily chosen in $\mathscr{L}$ so that $\|\omega_{t=j_0}-\varpi_{j_0}\|_{C^{\kappa,\alpha}}\leq \vareps_0$.  
  But given $\varpi\in\mathscr{L}$, there exists $t_0$ so that $\|\omega_{t_0}-\varpi_{j_0}\|_{C^{\kappa,\alpha}}\leq \vareps_0$, 
  and for all $t\geq t_0$, ${\rm d}_{C^{\kappa,\alpha}}(\omega_t,\mathscr{L})<\vareps_0$;  
  the assertion on the geniricity of $\varpi_0$ follows by applying the previous construction to $(\omega_{t+t_0-j_0})$ 
  for which one can keep the same $j_0$. 
  \cqfd

\section{Constant scalar curvature Poincaré type Kähler metrics on the complement of a divisor}  \label{sctn_asympttcs}
 We fix in this part a compact Kähler manifold $(X,\omega_X)$, 
 and a simple normal crossing divisor $D=\sum_{j=1}^N D_j$ (the $D_j$ are the smooth irreducible components). 
 
 \subsection{Basic tools and statements of the results}    \label{subsec_fib}
 \textit{Reminder: fibration near the divisor, and previous results. --- } 
 The details concerning following material, necessary for what follows, can be found in \cite{auv1,auv2}. 
 To fix ideas, assume first that $D$ is smooth, and even reduced to one component. 
 One can endow a tubular neighbourhood $\mathcal{N}_A$ of $D$ with an $\mathbb{S}^1$-action and an $\mathbb{S}^1$-invariant projection $p:\mathcal{N}_A\to D$, 
 and construct an $\mathbb{S}^1$-invariant function $t$ such that $t=\log\big(-\log(|\sigma|^2)\big)+\mathcal{O}(e^{-t})$ at any order for Poincaré type metrics on $\XD$, 
 where $\sigma\in\mathscr{O}([D])$ is such that $D=\{\sigma=0\}$, and $|\cdot|$ any smooth hermitian metric on $[D]$. 
 Up to adjusting $\mathcal{N}_A$ and $A\in \R$, we get this way a circle fibration 
 $q=(t,p):\mathcal{N}_A\backslash D\to [A,+\infty)\times D$, 
 for which we construct a connexion 1-form $\eta$, such that $d^ct=2e^{-t}\eta+\mathcal{O}(e^{-t})$ at any order 
 -- notice the analogy with the formula $J_{\C}dt=2e^{-t}d\vthet$ on the model $\Delta^*$. 
 Setting moreover $\omega:=\omega_X-dd^c\log\big(-\log(|\sigma|^2)\big)$ with a good choice of $|\cdot|$ -- which does not affect what precedes --
 we have that $\omega$ is of Poincaré type, and near $D$, 
  \begin{equation}   \label{eqn_aspttcsomega}
   \omega=dt\wedge 2e^{-t}\eta + p^*\omega_D +\mathcal{O}(e^{-t}) 
  \end{equation}
 at any order, where $\omega_D=\omega_X|_D$. 
 This can be generalised when the divisor has several components, and the fibrations respect the crossings when they exist; 
 we simply add $j$ indices to specify the component $D_j$ near which we work. 
 
 In \cite{auv2}, it is proven, using these fibrations, 
 that if a metric of Poincaré type $\omega_{\varphi}=\omega+dd^c\varphi$ has constant scalar curvature, 
 then $\varphi-\sum_{j=1}^Na_jt_j$ is in $C^{\infty}(\XD)$, i.e. is bounded at any order with respect to 
 the model Poincaré type metric $\omega$, where the $a_j$ are $<1$, 
 and given by $\frac{2}{1-a_j}=\sbar_{D_j}-\sbar$;  
 here $\sbar$ (resp. $\sbar_{D_j}$) denotes the mean scalar curvature of Poincaré type metrics of class $\omega_X$ (resp. $\omega_{D_j}$) 
 on $\XD$ (resp. on $D_j\backslash\sum_{\ell\neq j}D_{\ell}$) -- 
 recall the formulae 
 $\sbar=-4\pi m\frac{\chern_1(K_X[D])\cdot[\omega_X]^{m-1}}{[\omega_X]^{m}}$ 
 and $\sbarDj=-4\pi m\frac{\chern_1(D_j)\cdot \chern_1(K_X[D])\cdot[\omega_X]^{m-2}}{\chern_1(D_j)\cdot[\omega_X]^{m-1}}$, 
 typical of the Poincaré setting. 
 
 Considering again a fixed $D_j$ and given $(\kappa,\alpha)\in \N\times[0,1)$, one can use the circle action on its neighbourhood to decompose 
 functions $f\in C^{\kappa,\alpha}(\XD)$ 
 (``bounded functions at order $(\kappa,\alpha)$ for $\omega$ on $\XD$'' -- see \cite[\S 1.2]{auv1} for the exact definition) 
 as an $\mathbb{S}^1$-invariant part $\Pi_{0,j}f$ and a part $\Pi_{\perp,j}f$ with null mean against $\eta_j$. 
 Of course 
 $\Pi_{0,j}f$ and $\Pi_{\perp,j}f$ remain in $C^{\infty}(\XD)$ -- 
 in particular, seen as a function on $[A,\infty)\times\big(D_j\backslash \sum_{\ell\neq j}D_{\ell}\big)$, $\Pi_{0,j}f$ is  bounded up to order $(\kappa,\alpha)$ 
 for $dt_j^2+\omega|_{D_j}(\cdot,J_{D_j}\cdot)$; 
 moreover, as the fibres have length in $e^{-t_j}$ near $D_j$, if $T_j$ denotes the infinitesimal circle-action vector field, 
 $(e^{-t_j}T_j)^{k}\Pi_{\perp,j}f\in C^{\kappa-k,\alpha}$ near $D_j$ for all $k\leq \kappa$. 
 Thus for instance, if $f\in C^{\infty}(\XD)$, 
 then 
  \begin{equation}   \label{eqn_decomp}
   df=\partial_{t_j}\Pi_0f +p_j^*d(\Pi_{0,j}f_{t_j})+\mathcal{O}(e^{-t_j})
  \end{equation}
 and $\Delta_{\omega}f= (\partial_{t_j}-\partial_{t_j}^2)\Pi_{0,j}f+ p_j^*\big(\Delta_{\omega|_{D_j}}(\Pi_{0,j}f_{t_j})\big)+\mathcal{O}(e^{-t_j})$ 
 near $D_j$ at any order with respect to $\omega$, and where $\Pi_{0,j}f_{t_j}=(\Pi_{0,j}f)(t_j,\cdot)$. 
 
  ~
 
 \noindent
 \textit{Two theorems on Poincaré type Kähler metrics with constant scalar curvature. ---} 
  We can now state the main results of this part: 
  \begin{thm} \label{thm_asmpttcsPKcsc} 
    Assume that $D$ is smooth. 
    Let $\omega_{\varphi}$ be a constant scalar curvature metric of Poincaré type on $\XD$, of class $[\omega_X]$. 
    Then at the level of Riemannian metrics, one has near all component $D_j$ of $D$ the asymptotics: 
     \begin{equation*}
      g_{\varphi} = a_{j}\big(dt_j^2 + (2e^{-t_j}\eta_j)^2\big) + p_j^*h_{j} +\mathcal{O}(e^{-\delta t_j})
     \end{equation*}
    at any order, with $\delta>0$, $h_{j}$ a Kähler metric on $(D_j,J_{D_j})$ such that $[h_{j}(J_{D_j}\cdot,\cdot)]=[\omega_X|_{D_j}]$ 
    and with constant scalar curvature, 
    and where $a_{j}=2(\sbar_{D_j}-\sbar)^{-1}>0$. 
   \end{thm}
  and: 
   \begin{thm}   \label{thm_cscKdivisor}
    Assume that there exists a Poincaré type Kähler metric of class $[\omega_X]$ on $\XD$. 
    Then for all components $D_j$ of $D$ there exists a (Poincaré type) Kähler metric on $D_j\backslash\sum_{\ell\neq j}D_{\ell}$. 
   \end{thm}
  
  The rest of this part is mainly devoted to the proof of Theorem \ref{thm_asmpttcsPKcsc}.  
  The first step is the construction of a \textit{family of almost constant scalar curvature metrics} on a fixed $D_j$, 
  along which a proof of Theorem \ref{thm_cscKdivisor} with $D$ smooth is provided; 
  this is done in next section, 
  as well as the generalisation of Theorem \ref{thm_cscKdivisor} to the general case. 
  In the analytical subsequent sections, we use the produced family and Proposition \ref{prop_asympparam} to end the proof of Theorem \ref{thm_asmpttcsPKcsc}.

  \subsection{A family of almost constant scalar curvature on $D$}   \label{subsec_link}
  Assuming $D$ reduced to one component -- hence smooth -- and starting with a constant scalar curvature $\omega_{\varphi}$ on $\XD$, 
  we use the results of Part \ref{sctn_mdl} us to identify a specific family of almost constant scalar curvature 
  as defined in Part \ref{sctn_almcscKfamily}:
   \begin{prop} \label{prop_almcnstt}
    Assume that $D$ is reduced to one component and that $\scal(\omega_{\varphi})=\sbar$ on $\XD$, with $\omega_{\varphi}=\omega+dd^c\varphi$ 
    of Poincaré type of class $[\omega_X]$. 
    Then for $T$ large enough, 
    $(\omega^{\varphi}_t)_{t\geq T}:=\big(\omega_D+dd^c_D\Pi_0\varphi(t,\cdot)\big)_{t\geq T}$, with $\omega_D=\omega_X|_D$, is a family of Kähler metrics on $D$ 
    of almost constant curvature, with extinguishing variation. 
    Moreover, for any $(\kappa,\alpha)$, and any positive $\ell$, $(\partial_t^{\ell}\Pi_0\varphi)_t\to 0$ in $C^{\kappa,\alpha}(D)$ 
    as $t$ goes to $\infty$. 
   \end{prop}
   
  \prf. --- Observe first that we can assume $\varphi\in C^{\infty}(\XD)$; 
   for this, just replace $\omega_X$ by $\frac{1}{a}\omega_X$, with $a=(\sbar_D-\sbar)^{-1}$, which imposes $\sbarD=\sbar+2$ (see previous section). 
   Now the principle of the proof goes as follows: 
   to check convergences, we proceed by contradiction: assuming the desired convergences do not hold, 
   we use the boundedness of $\varphi$ at any order as well as its almost $\mathbb{S}^1$-invariance to identify some subsequence 
   of $\Pi_0\varphi$ on compact subsets of $[A,\infty)\times D$
   converging to a limit fitting in the framework of Section \ref{sctn_mdl}, 
   and use Theorem \ref{thm_splttng} to contradict the assumption. 
   
   We thus set on $D$, for $t\geq A$, $\omega_t^{\varphi}=\omega_D+dd^c_D(\Pi_0\varphi)_t$; 
   up to increasing $A$, these are indeed metrics, uniformly bounded below, 
   and uniformly bounded in $C^{\kappa,\alpha}(D)$ for all $(\kappa,\alpha)$; 
   more precisely, $t\mapsto \omega_t^{\varphi}$ is bounded in $C^{\kappa,\alpha}\big([A,\infty)\times D\big)$ for all $(\kappa,\alpha)$.  
  
   Let us assume that there exist $\vareps>0$, a sequence $(t_j)$, $\lim_{j\to\infty} t_j=\infty$, and $(z_j)$ in $D$ such that 
   $|\scal(\omega_{t_j}^{\varphi})(z_j)-\sbarD|=\sup_{w\in D}|\scal(\omega_{t_j}^{\varphi})(w)-\sbarD|\geq \vareps$. 
   We consider a subsequence of $(z_j)$, still denoted $(z_j)$, converging to some $z\in D$. 
   Similarly, we pick any $(\kappa,\alpha)\in\N\times(0,1)$; 
   as $\Pi_0\varphi$ is bounded in $C^{\kappa,\alpha}\big([A,\infty)\times D\big)$, 
   we have for any $N$ a sequence $(\Pi_0\varphi)(\cdot+t_j,\cdot)_{j\geq j_0(N)}$ uniformly bounded in $C^{\kappa,\alpha}\big([-N,N]\times D\big)$. 
   A diagonal argument thus gives us $\varphi_{\infty}\in C^{\kappa,\alpha}(\R\times D)$ and a fixed subsequence $(t_j)$ such that $(\Pi_0\varphi)(\cdot+t_j,\cdot)$ converges 
   in $C^{\kappa,\alpha/2}$ on every compact subset of $\R\times D$ to $\varphi_{\infty}$; 
   an extra diagonal argument gives the $C^{\infty}_{\rm loc}$-convergence (of a subsequence), hence $\varphi_{\infty}\in C^{\infty}(\R\times D)$, i.e. 
   is bounded at any order for $dt^2+h$ for any $h$ on $D$. 
   Now we see $\R\times D$ as a factor of $\Delta^*\times D$, endowed with the complex structure $J_{\C}\oplus J_D$; 
   we claim: 
   
    ~

   \noindent 
   \textit{The $(1,1)$-form }$\omega_{\varphi_{\infty}}:=\omega_{D}+dd^c(\varphi_{\infty}-t)$ 
   \textit{is a Poincaré type metric on $\Delta^*\times D$ (in the sense of Part \ref{sctn_mdl}), 
   of constant scalar curvature, equal to} $\sbar$. 
   
    ~
   
   This claim decomposes into several assertions: 
   one has to check that $\omega_{\varphi_{\infty}}\geq c\omega_0$, with $\omega_{0}:=\omega_D-dd^ct$, 
   that $\omega_{\varphi_{\infty}}$ is bounded with respect to this model at any order 
   and that $\scal(\omega_{\varphi_{\infty}})=\sbar$ 
   (by construction, $\varphi_{\infty}$ is $\mathbb{S}^1$-invariant). 
   We start by the positivity assertion; 
   it is actually an easy exercise to prove that it is enough to check it on vectors $\xi$ of type $\xi_D+\partial_t$, 
   with $\xi_D\in TD$ independent of $t$ and $\vartheta$ -- 
   use the $J$-invariance, 
   and the writing 
   $\omega_{\varphi_{\infty}}=(\partial_t^2-\partial_t)(\varphi_{\infty}-t)dt\wedge 2e^{-t}d\vartheta
                              +dt\wedge d_D^c\partial_t\varphi_{\infty}+d_D\partial_t\varphi_{\infty}\wedge 2e^{-t}d\vartheta
                              +(\omega_D+dd^c_D\varphi_{\infty})$.   
                              
   For $\xi$ as above, consider $\zeta:=\big(q^*(\xi^{\sharp_{\omega_0}})\big)^{\flat_{\omega}}$, so that $|(T q)_x\zeta-\xi_{q(x)}|_{\omega_0}=O(e^{-t(x)})$, 
   and in particular $|(T_x p)\zeta_x-(\xi_D)_{q(x)}|_{\omega_0}=O(e^{-t(x)})$. 
   Then for $x$ close to $D$ in $\XD$, 
   as $\omega_{\varphi}= (\partial_t^2-\partial_t)(\Pi_0\varphi-t)dt\wedge 2e^{-t}\eta
                          +dt\wedge p^*d^c(\partial_t\Pi_0\varphi_{t})+p^*d\partial_t\varphi_{t}\wedge 2e^{-t}d\vartheta
                          +p^*\omega_t^{\varphi}
                          +\mathcal{O}(e^{-t})$ (at any order),
    \begin{equation*}
     |\zeta|^2_{\omega_{\varphi},x}=(\partial_t^2-\partial_t)(\Pi_0\varphi-t)_{x}
                                    +2p^{*}(d\partial_t\Pi_0\varphi_t)_x(\zeta_x)
                                    +p^*g_{t(x)}(\zeta_x,\zeta_x)+\mathcal{O}\big(e^{-t(x)}\big)
    \end{equation*}
   with $g_t^{\varphi}=\omega_t^{\varphi}(\cdot,J_D\cdot)$ on $TD$, thus: 
    \begin{align*}
     |\zeta|^2_{\omega_{\varphi},x}=(\partial_t^2-\partial_t)(\Pi_0\varphi-t)_{x}
                                    +2(\xi_D\cdot\partial_t\Pi_0\varphi)_{q(x)}
                                    +g_{t(x)}^{\varphi}(\xi_D,\xi_D)_{q(x)}+\mathcal{O}\big(e^{-t(x)}\big),  
    \end{align*}
  whereas for any $(t_0,\vartheta_0,w_0)\in \Delta^*\times D$, 
    \begin{align*}
     |\xi|^2_{\omega_{\varphi_{\infty}},(t_0,\vartheta_0,w_0)}=&(\partial_t^2-\partial_t)(\varphi_{\infty}-t)_{(t_0,w_0)}   
                                                               +2(\xi_D\cdot\partial_t\varphi_{\infty})_{(t_0,w_0)}
                                                               +g_{t_0}^{\varphi_{\infty}}(\xi_D,\xi_D)_{(t_0,w_0)},
    \end{align*}
  with $g_{t_0}^{\varphi_{\infty}}=\big(\omega_D+dd^c_D[\varphi_{\infty}(t_0,\cdot)]\big)(\cdot,J_D\cdot)$ on $TD$. 
  Choosing now $x_j$ in $q^{-1}(t_0+t_j,w_0)$, 
  as $\varphi_{\infty}$ is the $C^{2}$-limit of $(\Pi_0\varphi)(\cdot+t_j,\cdot)$ on every compact subset of $\R\times D$, 
  we have by the latter two formulae that $|\zeta|^2_{\omega_{\varphi},x_j}$ tends to $|\xi|^2_{\omega_{\varphi_{\infty}},(t_0,w_0)}$. 
  On the other hand since $\omega_{\varphi}\geq c\omega$ on $\XD$ for some $c>0$, 
  for all $x$ close to $D$ one has $|\zeta|^2_{\omega_{\varphi},x}\geq c |\zeta|^2_{\omega,x}$; 
  reasoning as above, one sees moreover that $|\zeta|^2_{\omega,x_j}$ tends to $|\xi|^2_{\omega_0,(t_0,\vartheta_0,w_0)}$, 
  hence the positivity assertion. 
  
  Bounds on $\omega_{\varphi_{\infty}}$ at any order 
  being an immediate consequence of $\varphi_{\infty}\in C^{\infty}(\R\times D)$,  
  we are therefore left with the constant scalar curvature assertion. 
  For this we use the asymptotic decomposition 
  $\omega_{\varphi}$ as above 
  together with that of $\omega$ \eqref{eqn_aspttcsomega}, 
  and the formula 
  $\scal(\omega_{\varphi})
    =2\Lambda_{\omega_{\varphi}}\big[\varrho(\omega)-\frac{1}{2}dd^c\log\big(\frac{\omega_{\varphi}^m}{\omega^m}\big)\big]$ 
  of Kähler geometry. 
  All computations done and dropping the $p^*$ for simplicity, one has:
   \begin{align*}
    \scal(\omega_{\varphi}&)
               =e^{-f^{\varphi}}\Big[2(m-1)\frac{(\omega^{\varphi}_t)^{m-2}\wedge \varrho(\omega_D)}{(\omega_D)^{m-1}}                                                                \\
                &-(m-1)(m-2)\frac{d_D\partial_t(\Pi_0\varphi)\wedge d_D^c\partial_t(\Pi_0\varphi)\wedge(\omega^{\varphi}_t)^{m-3}\wedge\varrho(\omega_D)}{(\omega_D)^{m-1}}            \\
                &-\big(1+(\partial_t^2-\partial_t)f^{\varphi}\big)\frac{(\omega^{\varphi}_t)^{m-1}}{(\omega_D)^{m-1}}     
                 +2(m-1)\frac{d_D\partial_t(\Pi_0\varphi)\wedge d_D^c\partial_t(\Pi_0\varphi)\wedge(\omega^{\varphi}_t)^{m-2}}{(\omega_D)^{m-1}}                                       \\
                &-(m-1)\frac{dd^c_Df^{\varphi}_t\wedge(\omega^{\varphi}_t)^{m-2}}{(\omega_D)^{m-1}}                                                                                    \\
                &+(m-1)(m-2)\frac{dd^c_Df^{\varphi}_t\wedge d_D\partial_t(\Pi_0\varphi)\wedge d_D^c\partial_t(\Pi_0\varphi)\wedge(\omega^{\varphi}_t)^{m-3}}{(p^*\omega_D)^{m-1}}\Big]
                 +\mathcal{O}(e^{-t}).
   \end{align*}
  near $D$ on $\XD$, with $f^{\varphi}=\Pi_0\log\big(\frac{\omega_{\varphi}^m}{\omega^m}\big)
                                      =\log\big(\frac{\omega_{\varphi}^m}{\omega^m}\big)+\mathcal{O}(e^{-t})$, 
  as $\frac{\omega_{\varphi}^m}{\omega^m}\in C^{\infty}(\XD)$; 
  more explicitly, 
  $f^{\varphi}=p^*\log\big[\big(1+(\partial_t^2-\partial_t)\Pi_0\varphi-|d_D\Pi_0\varphi|^2_{\omega_{t}^{\varphi}}\big)\frac{(\omega_{t}^{\varphi})^{m-1}}{\omega_D^{m-1}}\big]
   +\mathcal{O}(e^{-t})$. 
  Similarly, setting $f^{\infty}=\log\big(\frac{\omega_{\varphi_{\infty}}^m}{(dt\wedge2e^{-t}d\vthet+\omega_D)^m}\big)$, 
   \begin{align*}
    \scal(\omega_{\varphi_{\infty}}&)
               =e^{-f^{\infty}}\Big[2(m-1)\frac{(\omega^{\varphi_{\infty}}_t)^{m-2}\wedge\varrho(\omega_D)}{(\omega_D)^{m-1}}                                                                \\
                &-(m-1)(m-2)\frac{d_D\partial_t\varphi_{\infty}\wedge d_D^c\partial_t\varphi_{\infty}\wedge(\omega^{\varphi_{\infty}}_t)^{m-3}\wedge\varrho(\omega_D)}{(\omega_D)^{m-1}}      \\
                &-\big(1+(\partial_t^2-\partial_t)f^{\infty}\big)\frac{(\omega^{\varphi_{\infty}}_t)^{m-1}}{(\omega_D)^{m-1}}     
                 +2(m-1)\frac{d_D\partial_t\varphi_{\infty}\wedge d_D^c\partial_t\varphi_{\infty}\wedge(\omega^{\varphi_{\infty}}_t)^{m-2}}{(\omega_D)^{m-1}}                                 \\
                &-(m-1)\frac{dd^c_Df^{\infty}_t\wedge(\omega^{\varphi_{\infty}}_t)^{m-2}}{(\omega_D)^{m-1}}                                                                                   \\
                &+(m-1)(m-2)\frac{dd^c_Df^{\infty}_t\wedge d_D\partial_t\varphi_{\infty}\wedge d_D^c\partial_t\varphi_{\infty}\wedge(\omega^{\varphi_{\infty}}_t)^{m-3}}{(\omega_D)^{m-1}}\Big]
   \end{align*}
  on $\Delta^*\times D$. 
  Hence for $(t_0,\vartheta_0,w_0)$ and $(x_j)$ as above, 
  $\scal(\omega_{\varphi})_{x_j}$ tends to $\scal(\omega_{\varphi_{\infty}})_{(t_0,\vartheta_0,w_0)}$; 
  since $\scal(\omega_{\varphi})_{x_j}=\sbar$ for all $j$, $\scal(\omega_{\varphi_{\infty}})_{(t_0,\vartheta_0,w_0)}=\sbar$, 
  and this holds for any $(t_0,\vartheta_0,w_0)\in\Delta^*\times D$: 
  $\omega_{\varphi_{\infty}}$ has constant scalar curvature, equal to $\sbar$. 
  
  Now Theorem \ref{thm_splttng} says that $\varphi_{\infty}$ does not depend on $t$, and $\omega_{\varphi_{\infty}}$ 
  is thus a product $dt\wedge 2e^{-t}d\vartheta+\omega^{\psi}_D$, with $\omega^{\psi}_D=\omega_D+dd^c_D\psi$, $\psi\in C^{\infty}(D)$, 
  and with $\scal(\omega^{\psi}_D)$ constant, equal to $\sbar-2=\sbarD$. 
  As $\omega^{\psi}_D$ is the $C^{2}$-limit of $(\omega_{t_j}^{\varphi})$, 
  we have a contradiction with the preliminary assumption $\sup_{D}|\scal(\omega_{t_j}^{\varphi})-\sbarD|\geq\vareps$ for all $j$ large enough. 
  
  Similarly, starting with an assumption such as $\sup_D\big|(\nabla^D)^\kappa(\partial_t^{\ell}\Pi_0\varphi)_{t_j}\big|\geq \vareps>0$ for $j\geq j_0$, 
  with $(t_j)$ tending to $\infty$ and $\kappa\geq0$, $\ell\geq1$, 
  we similarly pass to a diagonal subsequence of the $(\Pi_0\varphi)(\cdot+t_j,\cdot)$ converging in $C^{\infty}$ 
  on every compact subset of $\R\times D$ to a function $\varphi_{\infty}\in C^{\infty}(\R\times D)$. 
  As above, $\varphi_{\infty}$ is associated to a Poincaré type metric with constant scalar curvature on $\Delta^*\times D$, 
  and by Theorem \ref{thm_splttng}, $\varphi_{\infty}$ is independent of $t$, with contradicts the assumption. 
 \cqfd
  
   ~
   
  \noindent
  \prf \textit{of Theorem \ref{thm_cscKdivisor}}. ---
   Notice that the construction of the function(s) $\varphi_{\infty}$ in the above proof does not require any of the contradictory assumptions, 
   hence the result when the divisor is reduced to one component, which readily generalises to the smooth divisor case. 
   In the simple normal case, one still works near one component, 
   but far from the other ones in the sense that the functions $\varphi_j$ are considered on sets of type 
   $[-N,N]\times K_p$, with $(K_p)$ an exhausting sequence of compact subsets of $D_j\backslash \sum_{\ell\neq j}D_{\ell}$,  
   before the use of the diagonal arguments. 
   Notice that in this case, the uniform $C^0$ bound is of type $|\varphi_j|\leq C(1+\sum_{\ell\neq j} |t_{\ell}|)$. 
  \cqfd

 ~
 
 For simplicity, \textbf{\textit{we assume from now on and until the end of this part that $D$ is reduced to one component}}, 
 and thus drop the $j$ indexes; 
 we indeed work around one fixed component in the subsequent sections, 
 so that all what is done below readily generalizes to the smooth $N\geq2$ case. 
 \textbf{\textit{We also keep the normalisation $\sbarD=\sbar+2$, 
 and fix the Poincaré type Kähler metric 
 $\omega_{\varphi}=\omega+dd^c\varphi$ of constant scalar curvature, all along the rest of this part}}.  
 
 \subsection{A fifth order equation on the potential $\varphi$ near the divisor}
  \textit{Localisation. ---} Choose a polydisc $(z^1,\dots,z^m)$ of holomorphic coordinates near any point in $D$, 
  such that $z^1$ is a local equation of $D$. 
  Set $Z:=\Re\big[z^1(\log z^1)\frac{\partial}{\partial z^1}\big]$ locally; 
  notice that $Z$ is bounded with bounded derivatives at any order, with respect to any Poincaré type metric. 
  Then we claim that for all $f$ in $C^{\infty}(X\backslash D)$ near $D$, 
   \begin{equation*}
    Z\cdot f =  \frac{1}{2}\partial_t\Pi_0f + \mathcal{O}(e^{-t}), 
   \end{equation*}
  where the $\mathcal{O}(e^{-t})$ is understood at any order in Poincaré type metric. 
  This follows easily from decomposition \eqref{eqn_decomp}, and 
  the estimates 
   \begin{equation*}
    p^*d_D(\Pi_0f)_t(Z)=\mathcal{O}(e^{-t}) \qquad\text{and}\qquad dt(Z)=\frac{1}{2}+\mathcal{O}(e^{-t}) 
   \end{equation*}
  at any order. 
  For the first one, use that $(\Pi_0f)_t$ has bounds at any order on $D$, uniform in $t$, 
  and that in our open subset of work, $z^j=p^*(z^j|_D)+\mathcal{O}(e^{-t})$ for $j=2,\dots,N$, at any order. 
  Now $t=\log\big(-\log(|\sigma|^2)\big)+\mathcal{O}(e^{-t})$ where $\sigma$ is some (global) section associated to $D$; 
  again, the error term is understood at any order, 
  and since $\log\big(-\log(|\sigma|^2)\big)=\log\big(-\log(|z^1|^2)\big)+\mathcal{O}(e^{-t})$ with this $\mathcal{O}(e^{-t})$ understood likewise (see e.g. \cite{auv1}, proof of Prop. 1.2), 
  and $Z$ is bounded at any order, 
  we are done with the following elementary computation: 
   \begin{align*}
    Z\cdot \log\big(-\log(|z^1|^2)\big)&= \Re\Big[z^1(\log z^1)\frac{\partial \log\big(-\log(|z^1|^2)\big)}{\partial z^1} \Big]   \\
                                       &= \Re\Big[z^1(\log z^1) \frac{1}{z^1\log(|z^1|^2)}\Big]
                                        = \Re\Big[\frac{\log z^1}{\log(|z^1|^2)}\Big]                                             
                                        = \frac{1}{2}.
   \end{align*}

   ~
  
 \noindent
 \textit{The equation. ---} Mimicking what is done is Section \ref{sctn_mdl}, 
 we differentiate the equation $\scal(\omega_{\varphi})=\sbar$ with respect to $Z$, 
 and use $\omega=\omega_X-dd^ct+\mathcal{O}(e^{-t})$ at any order near $D$: 
  \begin{align*}
   0 &= Z\cdot \scal(\omega_{\varphi})  = \Delta_{\varphi}\big\langle\omega_{\varphi}, \mathcal{L}_{Z}(\omega+dd^c\varphi)\big\rangle_{\varphi}
                                       -2\big\langle\varrho_{\varphi}, \mathcal{L}_{Z}(\omega+dd^c\varphi)\big\rangle_{\varphi}                    \\
                                     &= -2\lich_{\varphi}\big(Z\cdot(\varphi-t)\big)
                                       +\Delta_{\varphi} \langle\omega_{\varphi}, \mathcal{L}_{Z}\omega_X \rangle_{\varphi}
                                       -2 \langle\varrho_{\varphi}, \mathcal{L}_{Z}\omega_X \rangle_{\varphi}+\mathcal{O}(e^{-t}) ,
  \end{align*}
 with $\lich_{\varphi}$ the Lichnerowicz associated to $\omega_{\varphi}$ 
 (see e.g. \cite[p.192]{ap} -- recall that $\omega_{\varphi}$ has constant scalar curvature). 
 One checks moreover that in Poincaré type metrics, 
 $\mathcal{L}_{Z}\omega_X=\mathcal{O}(e^{-t})$ at any order, 
 and thus $\Delta_{\varphi}\langle\omega_{\varphi}, \mathcal{L}_{Z}\omega_X\rangle_{\varphi}$ and 
 $\langle\varrho_{\varphi}, \mathcal{L}_{Z}\omega_X \rangle_{\varphi}$ are $\mathcal{O}(e^{-t})$ at any order near $D$. 
 Furthermore as $Z\cdot(\varphi-t)= \frac{1}{2}\partial_t(\Pi_0\varphi)-\frac{1}{2}+\mathcal{O}(e^{-t})$ 
 with the $\mathcal{O}$ understood at any order, 
 we get: 
  \begin{equation}  \label{eqn_lichnearD}
   \lich_{\varphi}\big(\partial_t(\Pi_0\varphi)\big) = \mathcal{O}(e^{-t})
  \end{equation}
 near $D$ at any order. 
 Observe that \eqref{eqn_lichnearD} makes sense \textit{globally} near $D$, 
 and since $D$ is compact we can indeed patch together the local equations, and sum them up into this single equation. 
 
 ~
 
 We now analyse the operator $\lich_{\varphi}$ in more detail, 
 in order to deduce asymptotics on $\partial_t(\Pi_0\varphi)$ from \eqref{eqn_lichnearD}.

 \subsection{Asymptotics of the Lichnerowicz operator of $\omega_{\varphi}$}
  Recall that near the divisor, 
   \begin{align*}
     \omega_{\varphi}  = \big(1+(\partial_t^2-&\partial_t)\Pi_0\varphi\big)  dt\wedge 2 e^{-t}\eta  \\
                        &+dt\wedge d_D^c\partial_t(\Pi_0\varphi) 
                        +d_D\partial_t(\Pi_0\varphi) \wedge 2e^{-t}\eta 
                        +p^*\omega_t^{\varphi} +\mathcal{O}(e^{-t}), 
   \end{align*}
  where $\omega_t^{\varphi}=\omega_D+dd^c_D(\Pi_0\varphi)_t$, and with the $\mathcal{O}(e^{-t})$ at any order with respect to $\omega$. 
  Now according to Proposition \ref{prop_almcnstt}, 
  $\partial_t(\Pi_0\varphi)=o(1)$ at any order with respect to $\omega$, 
  thus if we use less precise asymptotics, we can simplify the previous formula into
   \begin{equation}   \label{eqn_asympomega}
    \omega_{\varphi} = dt\wedge 2 e^{-t}\eta +p^*\omega_t^{\varphi}+o(1), 
   \end{equation}
  at any order with respect to $\omega$, and with $|\partial^{\ell}_t\omega_t^{\varphi}|_{\omega_D}=o(1)$ for all $\ell\geq0$. 
 In other words, though the component of $\omega_{\varphi}$ which is parallel to $D$ is not constant, 
 its derivatives in $t$ of positive order tend to 0 near $D$, 
 and likewise,  the mixed terms and their time derivatives at any order tend to 0 near $D$. 
 Consequently, at the level of Ricci forms, one has:
  \begin{equation}   \label{eqn_asymprhophi}
   \varrho_{\varphi} = -dt\wedge 2 e^{-t}\eta +p^*\varrho(\omega_t^{\varphi})+o(1),  
  \end{equation}
 with the $o(1)$ understood as above. 
 From these asymptotics we deduce on $\lich_{\varphi}$: 
  \begin{prop}  \label{prop_asymplich}
   Denote by 
   $\lich_{\omega_t^{\varphi}}$ the Lichnerowicz operator associated to $\omega_t^{\varphi}$ on $D$ for all $t\geq A$. 
   For any $\kappa$, one has on 
   $\lich^0_{\varphi}:=\Pi_0\circ\lich_{\varphi}\circ q^*:C^{\kappa+4,\alpha}\big([A,\infty)\times\R\big)\to C^{\kappa,\alpha}\big([A,\infty)\times\R\big)$ 
   the asymptotics
    \begin{equation}   \label{eqn_asymplich2}
     \lich^0_{\varphi}       = \frac{1}{2}\Big(\frac{\partial}{\partial t}-\frac{\partial^2}{\partial t^2}\Big)^2
                              +\Big(\frac{\partial}{\partial t}-\frac{\partial^2}{\partial t^2}\Big)
                              +\lich_{\omega_t^{\varphi}} 
                              + \Delta_{\omega_t^{\varphi}}\circ\Big(\frac{\partial}{\partial t}-\frac{\partial^2}{\partial t^2}\Big) 
                              +o(1). 
    \end{equation}
  \end{prop}
  
 
 Here the Hölder spaces are defined with respect to $dt^2+h$, with $h$ fixed on $D$. 
 
 ~
 
  \prf. ---   
   We do it for the map: $C^{4,\alpha} \to C^{0,\alpha}$. 
   Let $f\in C^{4,\alpha}\big([A,\infty)\times D\big)$, seen as $q^*f$ near $D$. 
   As $\omega_{\varphi}$ has constant scalar curvature, 
   $\lich_{\varphi}f=\frac{1}{2}\Delta_{\varphi}^2f+\langle\varrho_{\varphi},dd^cf\rangle_{\varphi}$. 
   By formula \eqref{eqn_asympomega} and the identity $\Delta_{\varphi}f=-\Lambda_{\varphi}dd^cf=-m\frac{dd^cf\wedge \omega_{\varphi}^{m-1}}{\omega_{\varphi}^m}$, 
   we can write $\Delta_{\varphi}f= (\partial_t-\partial_t^2 )f+ p^*\Delta_{\omega_t^{\varphi}}f+\vareps(f)$, 
   with $|\vareps(f)|_{C^{2,\alpha}(p^{-1}([t,\infty)\times D))}\leq\|f\|_{C^{4,\alpha}}\overline{\vareps}(t)$, 
   where $\overline{\vareps}(t)=o(1)$ at any order is independent of $f$. 
   One more application of $\Delta_{\varphi}$ yields 
   $\Delta_{\varphi}^2f= (\partial_t-\partial_t^2 )^2f+ p^*\Delta_{\omega_t^{\varphi}}^2f
                         +(\partial_t-\partial_t^2 )p^*\Delta_{\omega_t^{\varphi}}f+ p^*\Delta_{\omega_t^{\varphi}}(\partial_t-\partial_t^2 )f+\vareps(f)$ 
   where this time, $|\vareps(f)|_{C^{0,\alpha}(p^{-1}([t,\infty)\times D))}\leq\|f\|_{C^{4,\alpha}}\overline{\vareps}(t)$.  
   Moreover as the time derivatives of $\omega_t^{\varphi}$ are $o(1)$ at any order, 
   we replace $(\partial_t-\partial_t^2 )p^*\Delta_{\omega_t^{\varphi}}f$ by another $p^*\Delta_{\omega_t^{\varphi}}(\partial_t-\partial_t^2 )f$, 
   for the price of adding to $\vareps(f)$ a comparable error term. 
   
   Likewise, on the curvature term, 
   formulae \eqref{eqn_asympomega} and \eqref{eqn_asymprhophi}, 
   the pointwise inner product identity $\langle\varrho_{\varphi},dd^cf\rangle_{\varphi} = 
                 (\Lambda_{\varphi}\varrho_{\varphi})(\Lambda_{\varphi}dd^cf)
                 -m(m-1)\frac{\omega_{\varphi}^{m-2}\wedge\varrho_{\varphi}\wedge dd^cf}{\omega_{\varphi}^m}$,
   its analogues for the $\omega_{t}^{\varphi}$ , and the differentiation formula 
   $dd^cf=(\partial_t^2-\partial_t)f dt\wedge 2e^{-t}\eta+dd^c_Df
          +dt\wedge d^c_D\partial_tf+d_D\partial_tf\wedge 2e^{-t}\eta+\mathcal{O}(\|f\|_{C^{4,\alpha}}e^{-t})$ 
   give: 
   $\langle\varrho_{\varphi},dd^cf\rangle_{\varphi}=(\partial_t-\partial_t^2 )f
    + \big\langle\varrho(\omega_{t}^{\varphi}),dd^c_Df_t\big\rangle_{\omega_{t}^{\varphi}}+\vareps(f)$, 
   with $|\vareps(f)|_{C^{0,\alpha}([t,\infty)\times D)}\leq\|f\|_{C^{4,\alpha}}\overline{\vareps}(t)$. 
   
   Sum these expansions;  
   as $\lich_{\omega_{t}^{\varphi}}=\frac{1}{2}\Delta_{\omega_{t}^{\varphi}}^2
       +\big\langle\varrho(\omega_{t}^{\varphi}),dd^c_D\cdot\big\rangle_{\omega_{t}^{\varphi}}
       +\frac{1}{2}\big\langle d_D\scal(\omega_{t}^{\varphi}),d_D\cdot\big\rangle_{\omega_{t}^{\varphi}}$,  
   as $d_D\scal(\omega_{t}^{\varphi})=o(1)$ on $D$ at any order as well as all its time derivatives,   
   and denoting by $\lich_1$ the operator in the right-hand side of \eqref{eqn_asymplich2}, one has: 
   $\lich_{\varphi}f= \lich_1 f+ \vareps(f) $ with $\vareps(f) \in C^{0,\alpha}([A,\infty)\times D)$,  
   and $\|\vareps(f)\|_{C^{0,\alpha}([t,\infty)\times D)}\leq \|f\|_{C^{4,\alpha}}\overline{\vareps}(t)$. 
   Now conclusion follows by applying $\Pi_0$ to the latter equality, 
   and from the fact that $\Pi_0\lich_1 f=\lich_1 f$, as $f$ is $\mathbb{S}^1$-invariant.    
  \cqfd. 
  
  ~
  
  These asymptotics are not sufficient to conclude however, as the $\omega_t^{\varphi}$ may \textit{vary};  
  we hence slightly change our point of view in next section to address this difficulty.

 \subsection{Twisting the equation near the divisor}
  \noindent
  \textit{Conjugated Lichnerowicz operator. }--- 
  By definition, $\omega_D=\omega_X|_D$, and $\omega$ induces $\omega_D$ on $D$. 
  It is more useful for what follows to choose first $\omega_D$ as a limit point of $(\omega_{t}^{\varphi})$; 
  there might be no smooth $\tilde{\omega}_X$ in $[\omega_X]$ with restriction $\omega_D$ on $D$, 
  but we can still construct a Poincaré-type $\omega$ of class $[\omega_X]$, 
  and such that $\omega=dt\wedge 2e^{-t}\eta+p^*\omega_D+\mathcal{O}(e^{-t})$ at any order near $D$, see \cite[\S3.2.1]{auv1}; 
  of course, this does not affect what precedes. 
  \textit{\textbf{From now on, we keep these $\omega_D$ and $\omega$}}. 
  
  According to Propositions \ref{prop_asympparam} and \ref{prop_almcnstt}, 
  we can write $\omega_t^{\varphi}=\sigma_t^*\omega_D+o(1)$ for $t\geq A$, 
  with this $o(1)$ understood at any order in $t$ and $D$, 
  and $\partial_t^{\ell} \sigma_t=o(1)$ in $C^{\infty}(D)$ for all $\ell>0$. 
  Consider the following diffeomorphism of $[A,\infty)\times D$:
   \begin{equation*}   \label{eqn_sigmatilde}
    \tilde{\sigma} (t,w) := \big(t,\sigma_t(w)\big),
   \end{equation*}
  where $w$ is the variable on $D$. 
  We let $\tilde{\sigma}$ act on functions on $[A,\infty)\times D$ by pull-back, 
  with inverse action by push-forward. 
  Recall that $\lich^0_{\varphi}=\Pi_0\circ\lich_{\varphi}\circ q^*$, and define: 
   \begin{equation*}
    \lich_{\varphi}^{\tilde{\sigma}}:=(\tilde{\sigma}^{-1})^*\circ \lich^{0}_{\varphi}\circ\tilde{\sigma}^*, 
     \qquad \text{i.e.} 
     \qquad \lich_{\varphi}^{\tilde{\sigma}}u:=(\tilde{\sigma})_*\big[\Pi_0\lich_{\varphi}\big(q^*(\tilde{\sigma})^*u\big)\big], 
   \end{equation*}
  so that $\lich_{\varphi}^{\tilde{\sigma}}\big((\tilde{\sigma})_*v\big)=(\tilde{\sigma})_*\big(\Pi_0\lich_{\varphi}v\big)$ -- here we omit $q^*$. 
  According to \eqref{eqn_lichnearD}, 
  we thus have at any order, as $t$ goes to $\infty$: 
   \begin{equation}      \label{eqn_lichbarv}
    \lich_{\varphi}^{\tilde{\sigma}}\big((\tilde{\sigma})_*(\partial_t\Pi_0\varphi)\big) = \mathcal{O}(e^{-t}). 
   \end{equation}
  
 ~
 
  \noindent
  \textit{Asymptotic mapping properties of $\lich_{\varphi}^{\tilde{\sigma}}$. } --- 
  Observe the following, which follows from the properties of $\tilde{\sigma}$, 
  and a verification similar to that of Proposition \ref{prop_asymplich}: 
  \begin{prop}    \label{prop_o1}
   Denote by $\lich_{\Delta^*\times D}^{0}$ the operator 
    \begin{equation*}
     \lich_{\Delta^*\times D}^{0}
         =  \frac{1}{2}\Big(\frac{\partial}{\partial t}-\frac{\partial^2}{\partial t^2}\Big)^2
               + \Big(\frac{\partial}{\partial t}-\frac{\partial^2}{\partial t^2}\Big)
               + \lich_{\omega_D} 
               + \Delta_{\omega_D}\circ\Big(\frac{\partial}{\partial t}-\frac{\partial^2}{\partial t^2}\Big). 
    \end{equation*}
   on $[0,\infty)\times D$. 
  Then 
   \begin{equation*}
    \lich_{\varphi}^{\tilde{\sigma}}- \lich_{\Delta^*\times D}^{0}=o(1), 
   \end{equation*}
  that is, the coefficients of this difference tend to $0$ in $C^{\infty}(D)$ at any order in $t$, as $t\to\infty$. 
 \end{prop}
  \begin{rmk}
   As notation suggests, $\lich_{\Delta^*\times D}^{0}$ is nothing by the Lichnerowicz operator of 
   $dt\wedge 2e^{-t}d\vthet+\omega_D$ on $\Delta^*\times D$, restricted to $\mathbb{S}^1$-invariant functions. 
  \end{rmk}

  The interest of Proposition \ref{prop_o1} lies in the following technical result,
  which results from the study of $\lich^0_{\Delta^*\times D}$, 
  and is the analytical key-step in our study of the asymptotics of $\partial_t(\Pi_0\varphi)$, 
   and thus of those of $\varphi$; 
  let $\chi$ be a smooth cut-off function on $\R$, 
  with $\chi=0$ on $(-\infty,\frac{1}{3}]$ and $\chi=1$ on $[\frac{2}{3},\infty)$. 
   \begin{prop}           \label{prop_isom}
    For any $\kappa\geq 0$, $\alpha\in(0,1)$ and $\delta\in\R$, define  
     \begin{equation*}
      C^{\kappa+4,\alpha}_{0,\delta}\big([A,\infty)\times D\big)
                                 =\big\{e^{\delta (t-A)}u\in C^{\kappa+4,\alpha}\big([A,\infty)\times D\big)\big|\, 
                                        u(0,\cdot)=\partial_tu(0,\cdot)\equiv0\big\}, 
     \end{equation*}
    and $C^{\kappa,\alpha}_{\delta}\big([A,\infty)\times D\big)=e^{-\delta (t-A)}C^{\kappa,\alpha}\big([A,\infty)\times D\big)$. 
    Then there exist $\delta_0\in(0,1]$,  
    and functions $\psi_1, \dots,\psi_r\in C^\infty([A_0,\infty)\times D)$, where $r=\dim(\ker \lich_{\omega_D})$, 
    such that for all $\delta\in(0,\delta_0)$ and all $A\geq A_0$ large enough, 
    the $\psi_j$ are linearly independent on $[A+1,\infty)\times D$, and
     \begin{equation}   \label{eqn_isom1}
      \lich_{\varphi}^{\tilde{\sigma}} : 
        C^{\kappa+4,\alpha}_{0,\delta}\big([A,\infty)\times D\big) \oplus \spn\big( \chi(\cdot-A)\psi_j\big)_{j=1,\dots,r} 
      \longrightarrow C^{\kappa,\alpha}_{\delta}\big([A,\infty)\times D\big)  
     \end{equation}
    is an isomorphism. 
    
    Moreover, for any non-trivial $r$-tuple $(\lambda_1,\dots,\lambda_r)$, 
    $\sum_{j=1}^r \lambda_j\psi_j$ \emph{does not} tend to 0 as $t$ goes to infinity. 
    
    One also has an isomorphism 
     \begin{equation}   \label{eqn_isom2}
      \lich_{\varphi}^{\tilde{\sigma}} : L^{2,\kappa+4}_0\big([A,\infty)\times D\big) 
                                         \longrightarrow L^{2,\kappa} \big([A,\infty)\times D\big), 
     \end{equation}
    where $L^{2} \big([A,\infty)\times D\big)=\big\{u\in L^2_{\loc}\big|\,\int_A^{\infty}e^{A-t}dt\int_D |u|^2 \vol^D<\infty\big\}$, 
    $L^{2,\lambda}$ is the subspace of functions in $u\in L^{2,\lambda}_{\loc}$ with $u,\dots,\nabla^{\lambda}u\in L^2$, 
    and $L^{2,\lambda}_0$ the subspace of functions $u$ such that $u(A,\cdot)=(\partial_t u)(A,\cdot)=0$. 
   \end{prop}

  We devote next section to the proof of this result. 
  For now, we use it to establish the desired asymptotics on $\varphi$. 
  
   ~
   
  \prf \textit{of Theorem \ref{thm_asmpttcsPKcsc} from Proposition \ref{prop_isom}. } ---
   As observed above in equation \eqref{eqn_lichbarv},  
   if one sets $v:=\tilde{\sigma}^*(\partial_t\Pi_0\varphi)$, 
   then $\lich_{\varphi}^{\tilde{\sigma}}(v)\in C^{\kappa,\alpha}_{\delta}$ for any $\delta\leq1$, $(\kappa,\alpha)\in \N\times(0,1)$. 
   Taking now $A$ and $\delta_0$ as in the statement of Proposition \ref{prop_isom}, 
   one has $\lich_{\varphi}^{\tilde{\sigma}}\big(\chi(t-A)v\big)
             \in C^{\kappa,\alpha}_{\delta}\big([A,\infty)\times D\big)\subset L^{2,\kappa}\big([A,\infty)\times D\big)$, 
   and $\chi(t-A)v \in L^{2,\kappa+4}_0\big([A,\infty)\times D\big)$. 
   
   Pick $\delta\in(0,1]$, $\delta<\delta_0$. 
   According to isomorphisms \eqref{eqn_isom1} and \eqref{eqn_isom2} in Proposition \ref{prop_isom}, 
   there exist $w\in C^{\kappa+4,\alpha}_{0,\delta}\big([A,\infty)\times D\big)$, 
   and $\lambda_1,\dots,\lambda_r$, such that 
    \begin{equation*}
      \chi(t-A)v = w + \sum_{j=1}^r \lambda_j\psi_j. 
    \end{equation*}
   
   But $\partial_t\Pi_0\varphi$ tends to $0$ as $t$ goes to $\infty$ (Proposition \ref{prop_o1}), 
   hence $v=\tilde{\sigma}^*(\partial_t\Pi_0\varphi)$ does so. 
   Since this holds as well for $w$, we get that $\sum_{j=1}^r \lambda_j\psi_j$ tends to 0 as $t$ goes to $\infty$. 
   By Proposition $\ref{prop_isom}$, this implies $\lambda_1=\cdots=\lambda_r= 0$, 
   that is: $v\in C^{\kappa+4,\alpha}_{\delta}\big([A,\infty)\times D\big)$. 
   Such a statement is stable by pushing forward with $\tilde{\sigma}$, so that 
   $\partial_t(\Pi_0\varphi)$ is $C^{\kappa+4,\alpha}_{\delta}$ near $D$. 
   This holds for all $\kappa>0$; 
   after integrating along $t$, and adding the $\Pi_{\perp}$ component, we get the final statement: 
    \begin{equation*}
     \varphi = p^*\psi_D + \mathcal{O}(e^{-\delta t}),  
    \end{equation*}
   for some $\psi_D\in C^{\infty}(D)$, 
   with the $\mathcal{O}$ understood at any order near the divisor.  
   We easily see that the metric $g_D^{\psi_D}$, of Kähler form $\omega_D^{\psi_D}:=\omega_D+dd^c_D\psi_D$, 
   has constant scalar curvature on $D$ (see \cite[\S2.2]{auv2}), 
   and that 
    \begin{equation*}
     g_{\varphi}= dt^2+4e^{-2t}\eta^2 + p^*g_{\psi_D} + \mathcal{O}(e^{-\delta t}) 
    \end{equation*}
   near $D$ at any order. 
   \cqfd

   
  \subsection{Proof of Proposition \ref{prop_isom}}
  
   We subdivide this proof into three steps. 
   We first prove an analogous statement (Lemma \ref{lem_mapprop}) for the model operator $\lich_{\Delta^*\times D}^{0}$. 
   We then come back to $\lich_{\varphi}^{\tilde{\sigma}}$, 
   and exploit its asymptotic convergence to $\lich_{\Delta^*\times D}^{0}$ to deal with its Fredholm properties 
   and compute its index between relevant spaces (Lemma \ref{lem_fredholm}).  
   By contrast, its geometric origin 
   -- recall we came to $\lich_{\varphi}^{\tilde{\sigma}}$ from the study of the constant scalar curvature Poincaré type Kähler metric $\omega_{\varphi}$ -- 
   is also used to analyse its kernel in conclusion,   
   where we exhibit the functions $\psi_1, \dots, \psi_r$ of the statement of Proposition \ref{prop_isom}, 
   and deal with their asymptotic linear independence (Lemma \ref{lem_psi_j}). 
   
   \subsubsection{Mapping properties of $\lich_{\Delta^*\times D}^{0}$}
    We first state the following, 
    on which Proposition \ref{prop_isom} is partly modelled: 
    \begin{lem}   \label{lem_mapprop}
     \begin{enumerate}
     \item The map $\lich_{\Delta^*\times D}^{0}: L^{2,\kappa+4}_0\big([0,\infty)\times D\big) \to L^{2,\kappa} \big([0,\infty)\times D\big)$ 
           is an isomorphism for any $\kappa\geq 0$. 
     \item There exists $\vareps>0$ such that 
      \begin{equation*}   \label{eqn_isom3}
       \lich_{\Delta^*\times D}^{0}: C^{\kappa+4,\alpha}_{0,\delta}\big([0,\infty)\times D\big)
                                                   \longrightarrow C^{\kappa,\alpha}_{\delta}\big([0,\infty)\times D\big)
      \end{equation*}
     is an isomorphism for all $\delta\in (-\frac{1}{2}-\vareps, 0)$ for any $(\kappa,\alpha)\in\N\times (0,1,)$. 
     \item $\delta=0$ is a critical weight for $\lich_{\Delta^*\times D}^{0}$, 
     and if $(f_1,\dots,f_r)$ denotes a basis of $\ker\lich_{\omega_D}$, 
     then there exists $\delta_0>0$ such that 
      \begin{equation}   \label{eqn_isom4}
       \lich_{\Delta^*\times D}^{0}: C^{\kappa+4,\alpha}_{0,\delta}\big([0,\infty)\times D\big) \oplus \spn\big(\chi(t) f_j\big)
                                                   \longrightarrow C^{\kappa,\alpha}_{\delta}\big([0,\infty)\times D\big)
      \end{equation} 
     is an isomorphism for all $\delta\in(0,\delta_0)$ and any $(\kappa,\alpha)\in \N\times(0,1)$.  
     \end{enumerate}
    \end{lem}

    \begin{rmk}
     1. Observe that as $\lich_{\Delta^*\times D}^{0}$ is invariant by translation in the $t$ direction, 
        we can translate these statements on $[A,\infty)\times D$ for any $A$. 
        
     2. The map \eqref{eqn_isom4} is well-defined, 
        as $\lich_{\Delta^*\times D}^{0}\big(\chi(t) f_j\big)=\lich_{\omega_D}(f_j)=0$ on $\{t\geq1\}$ for $j=1,\dots,r$.  
    \end{rmk}

   ~
   
   \prf \textit{of Lemma \ref{lem_mapprop}}. --- 
    We start by points 1. and 2. 
    We will use several times the following inequality: 
    if $f(0)=0$ and $\int_0^{\infty}|f|^2 e^{2\delta t}dt<\infty$, then 
     \begin{equation}  \label{eqn_ineq1}
      \int_0^{\infty}|\partial_t f|^2e^{2\delta t} dt \geq \delta^2 \int_0^{\infty}|f|^2 e^{2\delta t}dt
     \end{equation}
    (see \cite[Lemme 6.1]{biq}); this follows from Cauchy-Schwarz inequality applied to 
     \begin{equation*}  \label{eqn_ibp01}
      0= \frac{1}{2}\int_0^{\infty}\partial_t (f^2 e^{2\delta t}) dt 
       = \int_0^{\infty}(f\partial_tf)e^{2\delta t} dt + \delta \int_{0}^{\infty} f^2e^{2\delta t} dt
     \end{equation*}
    where $f$ is smooth with compact support and vanishes at $t=0$, plus a density argument. 
    Now, to simplify expressions, we denote $\lich^{0}_{\Delta^*\times D}$ by $\lich$;  
    we recall that 
     $ \lich =  \frac{1}{2}\big(\partial_t-\partial_t^2\big)^2
               + \big(\partial t-\partial_t^2\big)
               + \lich_{\omega_D} 
               + \Delta_{\omega_D}\circ\big(\partial t-\partial_t^2\big)$,
    so that 
     \begin{align*}
      \int_{[0,\infty)\times D} & u\lich u e^{2\delta t} \, dt\vol_D =                                          \\
               &\frac{1}{2} \int_{[0,\infty)\times D}  u(\partial_t-\partial_t^2)^2 u \,e^{2\delta t} dt\vol_D  
                + \int_{[0,\infty)\times D}  u(\partial_t-\partial_t^2) u \,e^{2\delta t} dt\vol_D              \\
               &+ \int_{[0,\infty)\times D}  u\lich_{\omega_D} u \,e^{2\delta t} dt\vol_D                        
          \qquad\,+ \int_{[0,\infty)\times D}  u(\partial_t-\partial_t^2)\Delta_D u \,e^{2\delta t} dt\vol_D.    
     \end{align*}
   We deal successively with the different summands to estimate the positivity of $\int_{[0,\infty)\times D}  u\lich u \, e^{2\delta t} dt\vol_D$. 
   We assume from now on that $u\in L^{2,4}_{0,\delta}\big([0,\infty]\times D\big)$. 
   
    ~
    
  \noindent 
  \textit{First summand: $\int_{[0,\infty)\times D}  u(\partial_t-\partial_t^2)^2 u \,e^{2\delta t} dt\vol_D$.} 
   We claim that:
    \begin{align*}
      \int_{[0,\infty)\times D}  u(\partial_t-\partial_t^2)^2 u \,e^{2\delta t} dt\vol_D
       =&\int_{[0,\infty)\times D} e^{2\delta t} (\partial_t^2 u)^2 \, dt\vol_D                                     \\
        &-(1+2\delta)(1+4\delta)\int_{[0,\infty)\times D} e^{2\delta t} (\partial_t u)^2 \, dt\vol_D                \\
        & +2\delta^2(1+2\delta)^2\int_{[0,\infty)\times D}  u^2 e^{2\delta t} \, dt\vol_D
    \end{align*}

  Indeed (we assume $u$ smooth, and vanishing near infinity for convenience), if we notice that 
  $\partial_t-1=-e^{t}\circ\partial_t\circ e^{-t}$, we get: 
   \begin{align*}
    \int_0^{\infty}  u&(\partial_t-\partial_t^2)^2 u  \,e^{2\delta t} dt                                            \\
     &= \int_0^{\infty}  u\big[\partial_t(1-\partial_t)\big](\partial_t-\partial_t^2) u \,e^{2\delta t} dt
      = -\int_0^{\infty}  ue^{(1+2\delta) t}\partial_t\big[e^{-t}(\partial_t-\partial_t^2)\partial_t u\big] \, dt \\
     &= \int_0^{\infty}  e^{2\delta t}(\partial_t u) \big[(\partial_t-\partial_t^2)\partial_t u\big] \, dt 
        +(1+2\delta) \int_0^{\infty}e^{2\delta t} u \big[(\partial_t-\partial_t^2)\partial_t u\big] \, dt
   \end{align*}
  by integration by parts, after using that $u(0)=0$ to get rid of the boundary term at $t=0$.                                              
  Now:
   \begin{align*}
    \int_0^{\infty}  e^{2\delta t}(\partial_t u) \big[(\partial_t-\partial_t^2)\partial_t u\big] \, dt
     &= -\int_0^{\infty}  e^{(1+2\delta) t}(\partial_t u) \partial_t\big(e^{-t}\partial_t^2 u\big) \, dt                                  \\
     &= \int_0^{\infty} e^{2\delta t}(\partial_t^2 u)^2 \,dt + (1+2\delta)\int_0^{\infty} e^{2\delta t}(\partial_t u)(\partial_t^2 u) \,dt
   \end{align*}
  after another integration by parts and using that $\partial_tu|_{t=0}=0$ to get rid of the boundary term at $t=0$,  
  hence 
   \begin{align*}
    \int_0^{\infty}  u(&\partial_t-\partial_t^2)^2 u \,e^{2\delta t} dt                                             \\
     = \int_0^{\infty} & e^{2\delta t}(\partial_t u)^2 \,dt                                                          
       +  (1+2\delta) \bigg[ \underbrace{\int_0^{\infty}e^{2\delta t} u \big[(\partial_t-\partial_t^2)\partial_t u\big] \, dt}_{:=A}
                             + \underbrace{\int_0^{\infty} e^{2\delta t}(\partial_t u)(\partial_t^2 u) \,dt}_{:=B} \bigg].           
   \end{align*}
  Then 
   \begin{equation*}
    A= -\int_0^{\infty}  e^{1+2\delta t}  u \partial_t\big(e^{-t}\partial_t^2 u\big) \, dt 
     = B + (1+2\delta)\underbrace{\int_0^{\infty}  e^{2\delta t} u \partial_t^2 u \, dt }_{:=C}
   \end{equation*}
  (again, no boundary terms in the integration by parts), 
   \begin{equation*}
    B = \frac{1}{2}\int_0^{\infty} e^{2\delta t} \partial_t\big((\partial_t u)^2\big)\,dt 
      =-\delta \int_0^{\infty} e^{2\delta t} (\partial_t u)^2 \,dt
   \end{equation*}
  (no boundary term, $\partial_tu|_{t=0}=0$), and 
   \begin{align*}
    C=& -\int_0^{\infty} \partial_t\big( e^{2\delta t} u\big) \partial_t u \, dt 
     =  -\int_0^{\infty} e^{2\delta t}(\partial_t u)^2 \,dt  - 2\delta\int_0^{\infty}  e^{2\delta t} u \partial_t u \, dt  \\
     =&  -\int_0^{\infty} e^{2\delta t}(\partial_t u)^2 \,dt  + 2\delta^2 \int_0^{\infty}  e^{2\delta t} u^2 \, dt,    
   \end{align*}
  as $2\delta\int_0^{\infty}  e^{2\delta t} u \partial_t u \,dt
      =\delta\int_0^{\infty}  e^{2\delta t} \partial_t (u^2)
      =-2\delta^2 \int_0^{\infty}  e^{2\delta t} u^2 \, dt$. 
  The claim now readily follows from gathering these expressions for $A$, $B$ and $C$, 
  and integrating along $D$. 
   
  ~
  
  \noindent
  \textit{Second summand: }$\int_{[0,\infty)\times D}  u(\partial_t-\partial_t^2) u \,e^{2\delta t} dt\vol_D$. 
   We now see that 
    \begin{align*}
     \int_{[0,\infty)\times D}  u(\partial_t-\partial_t^2) u \,&e^{2\delta t} dt\vol_D                         \\
         = & \int_{[0,\infty)\times D}  e^{2\delta t}(\partial_t u)^2  \, dt\vol_D
            -\delta(1+2\delta)\int_{[0,\infty)\times D}  e^{2\delta t} u^2  \, dt\vol_D. 
    \end{align*}
   We proceed exactly as above; 
   assuming again that $u$ is smooth and vanishes for large $t$, we have
    \begin{align*}
     \int_0^{\infty}e^{2\delta t}u(\partial_t-\partial_t^2) u \, dt  
      =& - \int_0^{\infty} e^{(1+2\delta) t}u\partial_t\big(e^{-t}\partial_t u\big) \,dt                                           \\
      =& \int_0^{\infty} e^{2\delta t}(\partial_tu)^2\,dt + (1+2\delta)\int_0^{\infty} e^{2\delta t}u\partial_tu\,dt               \\
      =& \int_0^{\infty} e^{2\delta t}(\partial_tu)^2\,dt + \frac{1}{2}(1+2\delta)\int_0^{\infty} e^{2\delta t}\partial_t(u^2)\,dt \\
      =& \int_0^{\infty} e^{2\delta t}(\partial_tu)^2\,dt - \delta(1+2\delta)\int_0^{\infty} e^{2\delta t} u^2 \,dt,
    \end{align*}
   where we only used that $u(0,\cdot)=0$ to get rid of boundary terms. 
   We conclude as above. 
   
    ~
   
   \noindent
   \textit{Third summand:} $\int_{[0,\infty)\times D}  u\lich_{\omega_D} u \,e^{2\delta t} dt\vol_D$. 
    We only prove that this summand is nonnegative, independently of $\delta$. 
    Integrating first along $D$, as $\lich_{\omega_D}=(\mathcal{D}_D)^{*_D}\mathcal{D}_D$, this is straightforward: 
     \begin{align*}
      \int_{[0,\infty)\times D}  u\lich_{\omega_D} u \,e^{2\delta t} dt\vol_D
       =&\int_0^{\infty}e^{2\delta t} dt\vol_D\bigg(\int_D u(\mathcal{D}_D)^{*_D}\mathcal{D}_D u \vol_D\bigg)    \\
       =&\int_0^{\infty}e^{2\delta t} dt\vol_D\bigg(\int_D \big|\mathcal{D}_D u\big|_D^2 \vol_D \bigg). 
     \end{align*}

    ~
   
   \noindent
   \textit{Last summand:} $\int_{[0,\infty)\times D}  u(\partial_t-\partial_t^2)\Delta_D u \,e^{2\delta t} dt\vol_D$. 
    We use now a different approach to see that, for all $\delta$,  
     \begin{equation*}
      \int_{[0,\infty)\times D}  u(\partial_t-\partial_t^2)\Delta_D u \,e^{2\delta t} dt\vol_D
       \geq -\delta(1+\delta)\int_{[0,\infty)\times D} |d_Du|_D^2  \,e^{2\delta t} dt\vol_D. 
     \end{equation*}
    Consider for this an $L^2$ orthonormal basis $(\varphi_j)_{j\geq 0}$ of eigenfunctions of $\Delta_D$, 
    and call $\mu_j$ the nonnegative eigenvalue attached to $\varphi_j$, i.e. $\Delta_D\varphi_j=\mu_j\varphi_j$. 
    Set moreover $u=\sum_{j=0}^{\infty}u_j\varphi_j$; 
    the $u_j$ are thus functions of $t$, and as $u|_{t=0}=(\partial_t u)|_{t=0}=0$, 
    we have $u_j(0)=\partial_tu_j(0)=0$ for all $j\geq0$. 
    This decomposition yields
     \begin{align*}
      \int_{[0,\infty)\times D}  u(\partial_t-\partial_t^2)&\Delta_D u \,e^{2\delta t} dt\vol_D
        = \sum_{j=0}^{\infty} \mu_j\vl(D)\int_0^{\infty} e^{2\delta t} u_j(\partial_t-\partial_t^2)u_j\,dt                   \\
       &= \sum_{j=0}^{\infty} \mu_j\vl(D)\bigg(\int_0^{\infty} e^{2\delta t} (\partial_tu_j)^2\,dt
                                               -\delta(1+2\delta)\int_0^{\infty} e^{2\delta t} u_j^2\,dt\bigg),  
     \end{align*}
    since for all $j$, 
    $\int_0^{\infty} e^{2\delta t} u_j(\partial_t-\partial_t^2)u_j\,dt
     =\int_0^{\infty} e^{2\delta t} (\partial_tu_j)^2\,dt
       -\delta(1+2\delta)\int_0^{\infty} e^{2\delta t} u_j^2\,dt $, 
    see the paragraph ``\textit{second summand}'' above. 
    Now the introductory inequality \eqref{eqn_ineq1} gives us 
    $\int_0^{\infty} e^{2\delta t} (\partial_tu_j)^2\,dt\geq\delta^2\int_0^{\infty} e^{2\delta t} u_j^2\,dt$ for all $j\geq0$, 
    so that: 
     \begin{align*}
      \int_{[0,\infty)\times D}  u(\partial_t-\partial_t^2) \Delta_D u \,e^{2\delta t} dt\vol_D
        &\geq \big(\delta^2-\delta(1+2\delta)\big)\sum_{j=0}^{\infty}\mu_j\vl(D)\int_0^{\infty} e^{2\delta t} u_j^2\,dt       \\
        &=    -\delta(1+\delta)\int_{[0,\infty)\times D} |d_Du|_D^2  \,e^{2\delta t} dt\vol_D.
     \end{align*}
    
   ~
   
   \noindent
   \textit{Recapitulation. }
   According to the previous four paragraphs, we have for all $\delta\in\R$
    \begin{align*}
     \int_{[0,\infty)\times D}   e^{2\delta t} u\lich u  \, dt\vol_D                                                                
       \geq  & \frac{1}{2}\int_{[0,\infty)\times D} e^{2\delta t} (\partial_t^2 u)^2 \, dt\vol_D                                    \\  
             & +\frac{1}{2}(1-6\delta-8\delta^2)\int_{[0,\infty)\times D} e^{2\delta t} (\partial_t u)^2 \, dt\vol_D                \\
             & +\delta(1+2\delta)(2\delta-1)(\delta+1)\int_{[0,\infty)\times D}  u^2 e^{2\delta t} \, dt\vol_D                      \\
             & -\delta(1+\delta)\int_{[0,\infty)\times D} |d_Du|_D^2  \,e^{2\delta t} dt\vol_D                                     ,
    \end{align*}
   hence
    \begin{align*}
     \int_{[0,\infty)\times D}   e^{2\delta t} u\lich u  \, dt\vol_D
       \geq  &  \frac{1}{2}(1-7\delta)(1+\delta)\int_{[0,\infty)\times D} e^{2\delta t} (\partial_t u)^2 \, dt\vol_D                \\  
             & +\delta(1+2\delta)(2\delta-1)(\delta+1)\int_{[0,\infty)\times D}  u^2 e^{2\delta t} \, dt\vol_D                      \\
             & -\delta(1+\delta)\int_{[0,\infty)\times D} |d_Du|_D^2  \,e^{2\delta t} dt\vol_D,                                      
    \end{align*}
   where we get to the second inequality after applying the version of \eqref{eqn_ineq1} integrated along $D$ to the 
   summand $\int_{[0,\infty)\times D} e^{2\delta t} (\partial_t^2 u)^2 \, dt\vol_D $; 
   this is justified since $(\partial_tu)|_{t=0}=0$. 
   Since $(1-7\delta)(1+\delta)>0$ for $\delta\in\big(-1,\frac{1}{7}\big)$, 
   and in particular for $\delta\in(-1,0)$, we can apply \eqref{eqn_ineq1} to the summand 
   $\frac{1}{2}(1-7\delta)(1+\delta)\int_{[0,\infty)\times D} e^{2\delta t} (\partial_t u)^2 \, dt\vol_D$, 
   and get for those $\delta$: 
    \begin{equation*}  \label{eqn_ineq2}
     \begin{aligned}
      \int_{[0,\infty)\times D} e^{2\delta t}  u\lich u \, dt\vol_D 
          \geq & \frac{1}{2}\delta(1+\delta)(\delta-1)(\delta+2)\int_{[0,\infty)\times D}  u^2 e^{2\delta t} \, dt\vol_D        \\
               & -\delta(1+\delta)\int_{[0,\infty)\times D} |d_Du|_D^2  \,e^{2\delta t} dt\vol_D,  
     \end{aligned}
    \end{equation*}
   It readily follows that $\lich$ is an isomorphism $L^{2,4}_{\delta,0} \to L^2_{\delta}$ for all $\delta \in (-1,0)$. 
   Now as  $\ker \lich_{\omega_D}$ is never trivial (it contains al least constants on $D$), $\delta=0$ is clearly a critical value.    
   To see what its index becomes for small $\delta>0$, consider functions $u_j$, $j=0,\dots N$ on $D$, 
   such that $\lich\big(\sum_{j=0}^Nt^ju_j\big)=0$, with $u_N\not\equiv0$ . 
   The $N$th degree term in  $\lich\big(\sum_{j=0}^Nt^ju_j\big)$ is $\lich_{\omega_D}u_N$, hence $\lich_{\omega_D}u_N=0$.  
   Now if $N>0$, the $(N-1)$th degree term is $N(u_N+\Delta_{\omega_D}u_N)+\lich_{\omega_D}u_{N-1}$, which thus vanishes; 
   integrating it against $u_{N}\vol^D$, 
   we get: $0=N\int_D(u_N^2+|du_N|_D^2)\vol_D +\int_Du_N\lich_{\omega_D}u_{N-1}\vol_D
             =N\int_D(u_N^2+|du_N|_D^2)\vol_D +\int_Du_{N-1}\lich_{\omega_D}u_N\vol_D
             =N\int_D(u_N^2+|du_N|_D^2)\vol_D$, 
   so $u_N\equiv0$. 
   Hence $N=0$, and $\lich u_0=\lich_{\omega_D}u_0$: $u_0\in \ker\lich_{\omega_D}$. 
   Since $\lich$ is elliptic and invariant by translation, its index thus becomes $-\dim \ker \lich_{\omega_D}$ 
   for $\delta>0$ small, $\delta\in(0,\delta_0)$ say (see \cite[Theorem 1.2]{lm}). 
   We hence get an operator  
    \begin{equation*}
     \lich: L^{2,4}_{\delta,0}\big([0,\infty)\times D\big) \oplus \chi\ker\lich_{\omega_D} 
            \longrightarrow L^{2}_{\delta}\big([0,\infty)\times D\big)
    \end{equation*}
   of 0 index, which is an isomorphism as its domain lies inside $L^{2,4}_{-1/2,0}$, for $\delta\in(0,\delta_0)$. 
   
   The analogous statements with Hölder spaces instead of Sobolev spaces are deduced from these statements. 
   \hfill$\blacksquare$
   
  \subsubsection{Asymptotic kernel and Fredholm properties of $\lich_{\varphi}^{\tilde{\sigma}}$} 
   We keep the notation $\lich$ for $\lich_{\Delta^*\times D}^0$, 
   and take $\delta_0$ as in Lemma \ref{lem_mapprop}, which we assume $\leq1$. 
   \begin{lem}  \label{lem_fredholm}
    For all $\delta\in (-1,\delta_0)$, $\delta\neq0$, and large $A$,
    the operator 
     \begin{equation*}   \label{eqn_lmfrdhlm}
      \lich_{\varphi}^{\tilde{\sigma}} : C^{\kappa+4,\alpha}_{0,\delta}\big([A,\infty)\times D\big) 
                                   \longrightarrow C^{\kappa,\alpha}_{\delta}\big([A,\infty)\times D\big) 
     \end{equation*}
    is Fredholm, 
    with the same index as $\lich$; 
    in particular, it has index $-\dim\ker\lich_{\omega_D}$ for $\delta\in(0,\delta_0)$. 
    Moreover, $\lich_{\varphi}^{\tilde{\sigma}}$ has zero kernel for $\delta>-1$. 
   \end{lem}
   
    \prf. --- Since the coefficients of $\lich_{\varphi}^{\tilde{\sigma}}$ and $\lich$ differ by some $o(1)$ at any order, 
     we know that $\lich_{\varphi}^{\tilde{\sigma}}$ is Fredholm for the same $\delta$ as $\lich$, 
     and that the difference between their indices does not depend on $\delta$, see \cite[Thm. 6.1]{lm}. 
     Now, for instance, 
     $\lich : C^{\kappa+4,\alpha}_{0,\delta}\big([A,\infty)\times D\big) \to C^{\kappa,\alpha}_{\delta}\big([A,\infty)\times D\big)$
     is an isomorphism for $\delta\in(-1,0)$, independently of $A$. 
     Therefore, up to increasing $A$, the same assertion holds for $\lich_{\varphi}^{\tilde{\sigma}}$; 
     fixing such an $A$, $\lich_{\varphi}^{\tilde{\sigma}}$ and $\lich$ have same index, 0, for $\delta=\frac{1}{2}$, 
     and hence for all $\delta$. 
     From this and the Fredholm assertion we deduce that: 
      \begin{itemize}
       \item $\lich_{\varphi}^{\tilde{\sigma}}$ has no critical weight in $(-1,0)$, 
             and as a result has constant kernel and cokernel for $\delta$ in this range; 
             it is thus an isomorphism for $\delta\in(0,1)$; 
       \item $\lich_{\varphi}^{\tilde{\sigma}}$ has index $-\dim\ker\lich_{\omega_D}$ for $\delta\in(0,\delta_0)$, 
             and this corresponds to (minus) a cokernel dimension, 
             since $C^{\kappa+4,\alpha}_{0,\delta}\subset C^{\kappa+4,\alpha}_{0,-1/2}$ for such $\delta$. 
             \hfill $\blacksquare$
      \end{itemize}

  \subsubsection{Making explicitly $\lich_{\varphi}^{\tilde{\sigma}}$ into an isomorphism for small positive weights}
   According to the latter lemma, if $r$ denotes $\dim\ker\lich_{\omega_D}$, 
   there exist $r$ linearly independent functions $\psi_1^0, \dots,\psi_r^0$ 
   (\textit{resp.} $\psi_1, \dots,\psi_r$), which are in $\bigcap_{\vareps>0}C^{k+4,\alpha}_{-\vareps,0}$ such that 
   $\lich$ (\textit{resp.} $\lich_{\varphi}^{\tilde{\sigma}}$) is an isomorphism 
   from $C^{\kappa+4,\alpha}_{\delta,0}\oplus\spn(\psi_j^0)_{j=1,\dots,r}$ 
   (\textit{resp.} $C^{\kappa+4,\alpha}_{\delta,0}\oplus\spn(\psi_j^0)_{j=1,\dots,r}$) to $C^{\kappa,\alpha}_{\delta}$, 
   for all $\delta\in(0,\delta_0)$. 
    
   Now as seen at the end of the proof of Lemma \ref{lem_mapprop}, the $\psi_j^0$ are easy to determine: 
   one can 
   take $\psi_j^0= \chi(t)f_j$, $j=1,\dots,r$, where $(f_1,\dots,f_r)$ denotes a basis of $\ker\lich_{\omega_D}$.

   On the other hand, we have to look for other candidates for the $\psi_j$, since in general, 
   $\lich_{\varphi}^{\tilde{\sigma}}(f_j)$ is only $o(1)$, and not $\mathcal{O}(e^{-\delta t})$ with $\delta>0$ -- unless one of those $f_j$ is a constant, 
   which only settles the case when $\ker\lich_{\omega_D}$ is reduced to $\R$. 
   The good candidates involve the $f_j$ however: 
    \begin{lem}  \label{lem_psi_j}
     For $j=1,\dots,r$, set $Z_j:=\nabla^{\omega}(p^*f_j)$. Then
      \begin{equation}   \label{eqn_lmpsij}
       \lich_{\varphi}^{\tilde{\sigma}}\big[(\tilde{\sigma})_*\big(p^*f_j+\Pi_{0}(Z_j\cdot\varphi)\big)\big]=\mathcal{O}(e^{-t})
      \end{equation}
     at any order for $j=1,\dots,r$. 
     Moreover, the $(\tilde{\sigma})_*\big(f_j+\Pi_{0}(Z_j\cdot\varphi)\big)$ remain asymptotically linearly independent, 
     in the sense that $\sum_{j=1}^r \lambda_j\big[(\tilde{\sigma})_*\big(f_j+\Pi_{0}(Z_j\cdot\varphi)\big)\big]\to 0$ 
     as $t$ goes to $\infty$ implies $\lambda_1=\cdots=\lambda_r=0$. 
    \end{lem}

    \prf. --- 
    Let us settle the linear independence assertion. 
    We chose $\omega_D$ as a limit point of $\omega|_D+dd^c_D(\Pi_0\varphi)$, 
    and built $\omega$ so that $\omega|_D=\omega_D$, 
    which means that there exists $(t_{\ell})$ going to $\infty$ so that $(d_D\Pi_0\varphi)|_{\{t=t_{\ell}\}}$ tends to $0$ in $C^{0}(D)$. 
    Moreover $\Pi_{0}(Z_j\cdot\varphi)=Z_j\cdot(\Pi_0\varphi)+Z_j\cdot(\Pi_{\perp}\varphi)-\Pi_{\perp}(Z_j\cdot\varphi)=Z_j\cdot(\Pi_0\varphi)+\mathcal{O}(e^{-t})$ 
    at any order, for $j=1,\dots,r$. 
    Consequently, for any $(\lambda_{j})_{1\leq j\leq r}$, as the $Z_j$ are tangent to $D$, 
    $(\tilde{\sigma})^*\sum_{j=1}^r \lambda_j\big[(\tilde{\sigma})_*\big(f_j+\Pi_{0}(Z_j\cdot\varphi)\big)\big]\big|_{t=t_{\ell}}\to\sum_{j=1}^r \lambda_jf_j$ 
    as $\ell$ goes to $\infty$.  
    Now if $\sum_{j=1}^r \lambda_j\big[(\tilde{\sigma})_*\big(f_j+\Pi_{0}(Z_j\cdot\varphi)\big)\big]\to0$, 
    then $\sum_{j=1}^r \lambda_jf_j=0$, hence $\lambda_1=\cdots=\lambda_r=0$. 
    
    Notice that \eqref{eqn_lmpsij} is equivalent to $\lich_{\varphi}\big(f_j+\Pi_{0}(Z_j\cdot\varphi)\big)=\mathcal{O}(e^{-t})$. 
    We get these latter equations in a way similar to the starting point of this part, $\lich_{\varphi}\big(\partial_t(\Pi_0\varphi)\big)=\mathcal{O}(e^{-t})$. 
    Namely, $j$ being fixed, we choose in $X$ a neighbourhood $U$ of coordinates $(z^1,\dots,z^m)$ around some point of $D$ where
    $D$ is given by $z^1=0$, 
    we extend $Z_j|_D=\nabla^{\omega_D}f_j$ in $U$ independently of $z^1$ and denote this holomorphic extension by $\tilde{Z}_j$; 
    we have: $\tilde{Z}_j=Z_j+\mathcal{O}(e^{-t})$ at any order. 
    Since $\scal(\omega_{\varphi})$ is constant, 
     \begin{equation*}
      0=\tilde{Z}_j\cdot\scal(\omega_{\varphi})
       = \Delta_{\varphi}\big(\Lambda_\varphi\mathcal{L}_{\tilde{Z}_j}\omega_{\varphi}\big)
        - 2\big(\mathcal{L}_{\tilde{Z}_j}\omega_{\varphi},\varrho(\omega_{\varphi})\big)_{\varphi}. 
     \end{equation*}
    We will thus be done if we prove that $\mathcal{L}_{\tilde{Z}_j}\omega_{\varphi}=dd^c\big(p^*f_j+\Pi_0(\tilde{Z}_j\cdot\varphi)\big)$ on $U$ 
    up to some $\mathcal{O}(e^{-t})$ at any order, 
    as $\lich_{\varphi}=\frac{1}{2}\Delta_{\varphi}^2+(\varrho_{\varphi},dd^c\cdot)_{\varphi}$, 
    and as replacing $\tilde{Z}_j$ by $Z_j$ (which is globally defined around $D$) in the expression $\lich_{\varphi}\big(\Pi_0(\tilde{Z}_j\cdot\varphi)\big)$ 
    only gives rise to an error term which is $\mathcal{O}(e^{-t})$ at any order. 
    
    Now $\omega_{\varphi}=\tilde{\omega}_D+dd^c(\varphi-t)+\mathcal{O}(e^{-t})$ near $D$, where $\tilde{\omega}_D$ extends $\omega_D$ in $U$ independently of $z^1$, 
    and with the $O$ at any order. 
    Thus by Cartan's formula, 
    $\mathcal{L}_{\tilde{Z}_j}\omega_{\varphi}=d\big(\tilde{\omega}_D(\tilde{Z}_j,\cdot)\big)+\mathcal{L}_{\tilde{Z}_j}dd^c(\varphi-t)+\mathcal{O}(e^{-t})$.  
    Observe that $\omega_D(Z_j,\cdot)=g_D(J_D\nabla^{\omega_D}f_j,\cdot)=d^c_Df_j$ on $D$ by definition of $Z_j$, 
    and that $d^c(p^*f_j)=d^c_{X}f_j+\mathcal{O}(e^{-t})$ near $D$.  
    Moreover $\mathcal{L}_{\tilde{Z}_j}dd^c(\varphi-t)=dd^c\big(\tilde{Z}_j\cdot(\varphi-t)\big)=dd^c\big(\Pi_0(\tilde{Z}_j\cdot\varphi)\big)$, 
    as $dt(\tilde{Z}_j)$ and $\Pi_{\perp}(\tilde{Z}_j\cdot\varphi)$ are $\mathcal{O}(e^{-t})$ at any order. 
    Summing these estimates thus gives $\mathcal{L}_{\tilde{Z}_j}\omega_{\varphi}=dd^c\big(p^*f_j+\Pi_0(\tilde{Z}_j\cdot\varphi)\big)+\mathcal{O}(e^{-t})$ near $D$ on $U$, 
    as wanted. 
    \hfill $\blacksquare$
    
    ~
    
    To complete the proof of Proposition \ref{prop_isom}, just set $\psi_j:=\chi(A-t)(\tilde{\sigma})_*\big(f_j+\Pi_{0}(Z_j\cdot\varphi)\big)$, 
    for $j=1,\dots,r$. 
    \cqfd

 ~
 
\section{The extremal case}    \label{sctn_ext}
  \subsection{Extremal Kähler metrics on the model $\Delta^*\times(Y\backslash E)$}   \label{subsec_extmodel}
  \subsubsection{Potentials of extremal Kähler metrics of Poincaré type on $\Delta^*\times(Y\backslash E)$}
   We come back in this section to the point of view of Part \ref{sctn_mdl}, 
   and recall that $\omega_{\Delta^*}=-dd^ct=dt\wedge 2e^{-t}d\vthet$, 
   that $\omega_{\YE}$ is a fixed Kähler metric of Poincaré type of class $[\omega_Y]$ on $Y\backslash E$, 
   and that $\omega_0=\omega_{\Delta^*}+\omega_{\YE}$ on $\Delta^{*}\times (Y\backslash E)$. 
   
   As we will see below, the following class of potentials on $\Delta^*\times(Y\backslash E)$ is useful when working on extremal metrics: 
   we say that $\varphi\in \mathscr{K}'(\omega_0)$ 
   if $\partial_{\vthet}\varphi\equiv0$,
   $|\varphi|\leq C(\mathfrak{u}_Y+|t|)$, 
   $d\varphi$ is bounded at any order with respect to $\omega_0$, 
   and $\omega_{\varphi}=\omega_0+dd^c\varphi\geq c\omega_0$ for some $c>0$. 
   We say that $\omega_{\varphi}$ is \textit{extremal} if 
   $\mathsf{K}_{\varphi}:=\nabla^{\varphi}\scal(\omega_{\varphi})$ is holomorphic on $\Delta^*\times(Y\backslash E)$. 
   
   Observe now that if one takes $\varphi\in \mathscr{K}'(\omega_0)$ instead of $\varphi\in \mathscr{K}(\omega_0)$, 
   and if $\lich_{\varphi}(\dot{\varphi}-1)=0$, 
   \textit{without assuming that $\omega_{\varphi}$ has constant scalar curvature}, 
   as $\dot{\varphi}$ is bounded at any order for $\omega_0$ 
   and $\mathcal{D}_{\varphi}[e^{-t}(\dot{\varphi}-1)]=0$ is automatic as underlined in the proof of Lemma \ref{lem_D_u}, 
   then the proof of Lemma \ref{lem_Dv_isL2}, 
   which is the major step in the proof of Theorem \ref{thm_splttng}, remains valid. 
   In this regard, the aim of this section is: 
    \begin{thm}   \label{thm_splttngext}
     Let $\varphi\in \mathscr{K}'(\omega_0)$ such that $\omega_{\varphi}$ is extremal. 
     Then $\varphi=at+\psi$, with $a<1$ and $\psi\in \mathscr{E}(Y\backslash E)$. 
     Therefore, $\omega_{\varphi}=(1-a)\omega_{\Delta^{*}}+\omega_Y^{\psi}$, 
     where $\omega_Y^{\psi}=\omega_{\YE}+dd^c_Y\psi$ is thus extremal on $Y\backslash E$, of class $[\omega_Y]$, 
     and of Poincaré type if $E\neq \varnothing$. 
    \end{thm}
   
   The assertion ``$\psi\in \mathscr{E}(Y\backslash E)$'' means: 
   $\psi\in C^{\infty}_{\rm loc}(\YE)$, 
   $|\psi|\leq C(1+\mathfrak{u}_Y)$, and $d\psi$ is bounded at any order with respect to $\omega_{\YE}$. 
   
   Next paragraph is devoted to the proof; 
   for now, as evoked above, another elementary but crucial step in proving Theorem \ref{thm_splttngext} similarly to Theorem \ref{thm_splttng} is: 
    \begin{lem}
     Let $\varphi\in\mathscr{K}'(\omega_0)$ such that $\omega_{\varphi}$ is extremal. 
     Then $\lich_{\varphi}(\dot{\varphi}-1)=0$. 
    \end{lem}
   \prf. --- Let $\varphi$ be a potential as in the statement. 
    Then $\mathsf{K}_{\varphi}$ is holomorphic and bounded at any order, 
    so by Lemma \ref{lem_holovf} is tangent to $D$ and constant along $\Delta^*$: 
    $\kphi=\mathsf{K}_D$ for a fixed (holomorphic) $\mathsf{K}_D\in C^{\infty}(TD)$. 
    As in Part \ref{sctn_mdl}, 
    let $Z$ be the locally defined $\Re\big(z(\log z)\frac{\partial}{\partial z}\big)$ with $z$ the coordinate on $\Delta^*$. 
    Then as $\partial_{\vthet}\scal(\omega_{\varphi})=0$, 
    $\langle\kphi,Z\rangle_{\varphi}=Z\cdot\scal(\omega_{\varphi})=\frac{1}{2}\partial_t\scal(\omega_{\varphi})$. 
    Now $\mathcal{L}_{Z}\omega_{\varphi}=\frac{1}{2}dd^c(\dot{\varphi}-1)$, 
    and the infinitesimal variation of the scalar curvature along a deformation $\phi$ of the potential $\varphi$ 
    is given by $-2\lich_{\varphi}\phi+ \langle d\scal(\omega_{\varphi}),d\phi\rangle_{\varphi}$, 
    hence $Z\cdot\scal(\omega_{\varphi})=-\lich_{\varphi}(\dot{\varphi}-1)+ \frac{1}{2}\big\langle d\scal(\omega_{\varphi}),d(\dot{\varphi}-1)\big\rangle_{\varphi}$.  
    Therefore $\lich_{\varphi}(\dot{\varphi}-1)=0$ is equivalent to 
    $\big\langle d\scal(\omega_{\varphi}),d(\dot{\varphi}-1)\big\rangle_{\varphi}=\partial_t\scal(\omega_{\varphi})$, 
    and this latter follows from 
    $\big\langle d\scal(\omega_{\varphi}),d(\dot{\varphi}-1)\big\rangle_{\varphi}
     =e^{t}\big\langle d\scal(\omega_{\varphi}),d[e^{-t}(\dot{\varphi}-1)]\big\rangle_{\varphi}
       +(\dot{\varphi}-1)\langle d\scal(\omega_{\varphi}),dt\rangle_{\varphi}
     =e^{t}\cdot e^{-t}\partial_t\scal(\omega_{\varphi})+(\dot{\varphi}-1)dt(\mathsf{K}_{\varphi})$, 
    and $dt(\mathsf{K}_{\varphi})=0$, 
    as $\mathsf{K}_{\varphi}$ is tangent to $D$.  
    Here we used identity \eqref{eqn_nablau=holo}: $\nabla^{\varphi}\big(e^{-t}(\dot{\varphi}-1)\big)=e^{-t}\partial_t$,   
    holding in general, see the proof of Lemma \ref{lem_D_u}. 
    \cqfd

   \subsubsection{Proof of the splitting theorem for extremal metrics on $\Delta^*\times(Y\backslash E)$}
   Fix $\varphi\in \mathscr{K}'(\omega_0)$ such that $\omega_{\varphi}$ is extremal. 
   So far one has all the required ingredients to reach the statement : 
   $\int_{\Delta^*\times(Y\backslash E)}e^t|\mathcal{D}_{\varphi}(\dot{\varphi}-1)|_{\varphi}^2\vol^{\varphi}<\infty$. 
   To pass from this to the statement $\int_{\Delta^*\times(Y\backslash E)}e^t|\mathcal{D}_{\varphi}(\dot{\varphi}-1)|_{\varphi}^2\vol^{\varphi}=0$, 
   we use a subsequence similar to that of the proof of Lemma \ref{lem_Dv_is0}, adapted as follows: 
    \begin{enumerate}
     \item to have a subsequence bounded at any order on compact subsets, 
           replace $(\varphi_j)$ by $\big(\varphi_j-\varphi_j(0,\cdot)\big)$ 
           (this does not affect the attached metrics); 
           without loss of generality, we assume it converges at any order on compact subsets; 
     \item call $\varphi_{\infty}$ the $C^{\infty}_{\loc}$-limit;  
           it verifies $|\varphi_{\infty}|\leq C(\mathfrak{u}_Y+|t|)$, $d\varphi_{\infty}$ is bounded at any order for $\omega_0$, 
           and $\omega_0+dd^c\varphi_{\infty}\geq c \omega_0$, as uniform local bounds pass to the limit. 
           One gets by dominated convergence $\mathcal{D}_{\varphi_{\infty}}(\dot{\varphi_{\infty}}-1)=0$.  
           Again, $\ddot{\varphi_{\infty}}\equiv0$ follows from this, 
           which in the current situation provides 
           $\varphi_{\infty}=\psi_2t+\psi_1$ with $\psi_1,\psi_2$ functions on $Y\backslash E$. 
           As $d\varphi_{\infty}$ is bounded, this implies $d\psi_2\equiv 0$, and finally $\varphi_{\infty}=at+\psi_{\infty}$, 
           with $\psi_{\infty}$ a function on $Y\backslash E$ and $a$ a constant, necessarily $<1$, since then, 
           $\omega_{\varphi_{\infty}}=(1-a)\omega_{\Delta^*}+(\omega_Y+dd^c_Y\psi_{\infty})$;  
     \item compute $\lim_{+\infty}\mathcal{K}_{\varphi}
                    =\mathcal{K}_{\varphi_{\infty}}(0)=2\pi\int_{Y\backslash E}(\dot{\varphi_{\infty}}-1)^2|dt|_{\varphi_{\infty}}^2\frac{(\omega_0^{\varphi_{\infty}})^{m-1}}{(m-1)!} 
                    =2\pi(1-a)\vl(Y\backslash E)$.          
           Symmetrically, $\lim_{-\infty}\mathscr{K}_{\varphi}=2\pi(1-b)\vl(Y\backslash E)$ for some $b < 1$. 
           Hence we are done if we prove that $a=b$, 
           since as for $\varphi_{\infty}$, 
           we then deduce that $\varphi=at+\psi$, $\psi\in \mathscr{E}(Y\backslash E)$, from $\mathcal{D}_{\varphi}(\dot{\varphi}-1)=0$. 
    \end{enumerate}
   One proves that $a=b$ as follows. 
   Considering the sequence of times $t_j$ along which $\varphi_{\infty}$ arises, 
   one has $\omega_{\varphi}=(1-a)\omega_{\Delta^*}+\omega_Y^{\psi_{\infty}}+\vareps_j$ at $t=t_j$, 
   with $\vareps_j$ uniformly bounded at any order for $\omega_0$, and $\vareps_j\to0$ in $C^{\infty}_{\rm loc}(\YE)$ 
   as $j\to\infty$. 
   As a result, 
   $\scal(\omega_{\varphi})(t_j,\cdot)$, which is uniformly bounded along $Y\backslash E$, 
   tends to $\frac{2}{1-a}+\scal(\omega_Y^{\psi_{\infty}})$ in $C^{0}_{\rm loc}(Y\backslash E)$; 
   likewise, $\omega_{t_j}^{\varphi}$, uniformly bounded along $Y\backslash E$, tends to $\omega_Y^{\psi_{\infty}}$ in $C^{0}_{\rm loc}(Y\backslash E)$. 
   Thus, if $\Delta_{0,s}=\{0\leq t\leq s\}\subset\Delta^*$ for $s\geq 0$, as $\mathsf{K}_{\varphi}$ is tangent to $D$: 
    \begin{align*}
     0 &=\int_{\Delta_{0,t_j}\times (Y\backslash E)} e^t dt(\mathsf{K}_{\varphi})\vol^{\varphi}
        = \int_{\Delta_{0,t_j}\times (Y\backslash E)} d\scal(\omega_{\varphi})\wedge d^ce^t\wedge\frac{\omega_{\varphi}^{m-1}}{(m-1)!}  \\
       &=2\pi\Big(\int_{\{t_j\}\times (Y\backslash E)} \scal(\omega_{\varphi}) \frac{(\omega_{t_j}^{\varphi})^{m-1}}{(m-1)!}
                   -\int_{\{0\}\times(Y\backslash E)} \scal(\omega_{\varphi}) \frac{(\omega_{0}^{\varphi})^{m-1}}{(m-1)!}   \Big)    
        \qquad\text{as} \quad                                                                            dd^ce^t=0             \\
       &\xrightarrow{j\to\infty} 2\pi\Big(\vl(Y\backslash E)\big(\frac{2}{1-a}+\sbar_{Y\backslash E}\big)
                                          -\int_{\{0\}\times(Y\backslash E)} \scal(\omega_{\varphi}) \frac{(\omega_{0}^{\varphi})^{m-1}}{(m-1)!}\Big),
    \end{align*}
   by dominated convergence, 
   i.e. $\frac{2}{1-a}=\frac{1}{\vl(Y\backslash E)}\int_{\{0\}\times(Y\backslash E)} \scal(\omega_{\varphi}) \frac{(\omega_{0}^{\varphi})^{m-1}}{(m-1)!}-\sbar_{Y\backslash E}$. 
   The same holds symmetrically for $\frac{2}{1-b}$, 
   and thus $a=b$, which ends this proof. 
  \cqfd

 \subsection{Applications to extremal Poincaré type metrics on $X\backslash D$}   \label{subsec_PKext}
  On a compact Kähler manifold $(X,\omega_X)$ with a simple normal crossing divisor $D=\sum_{j=1}^N D_j$, 
  we say that a Poincaré type Kähler metric is \textit{extremal} if the gradient of its scalar curvature is real holomorphic. 
  As in the compact case, one can check that this corresponds to being a critical point of the Calabi functional, i.e. the squared $L^2$ norm of the scalar 
  curvature. 
  We set, for $j=1,\dots, N$, $E_j=\sum_{\ell\neq j} D_j\cap D_{\ell}$, 
  and use below the fibrations introduced in Section \ref{subsec_fib}. 
 
 \subsubsection{Existence of extremal metrics on the divisor}
   We start exploiting Theorem \ref{thm_splttngext} with a statement analogous to Theorem \ref{thm_cscKdivisor}: 
    \begin{prop}   \label{prop_extKdivisor}
     Assume that there exists an extremal Poincaré type Kähler metric of class $[\omega_X]$ on $\XD$. 
     Then for any $j$, 
     there exists an extremal (Poincaré type) Kähler metric of class $[\omega_X|_{D_j}]$ on $D_j\backslash E_j$. 
    \end{prop}
  \prf. --- 
   Proceeding as in the proof of Proposition \ref{prop_almcnstt},  
   after fixing $j$, 
   we construct a sequence $(t_k)$ going to $\infty$ such that $\varphi_k:=\big((\Pi_{0,j}\varphi)(\cdot+t_k,\cdot)-(\Pi_{0,j}\varphi)(t_k,x_0)\big)$, 
   with $x_0\in D_j\backslash E_j$, 
   converges at any order on compact subsets of $\R\times(D_j\backslash E_j)$ to some $\mathbb{S}^1$-invariant $\varphi_{\infty}\in \mathscr{K}'(\omega_{\Delta^*}+\omega_{D_j\backslash E_j})$
   -- the normalisation is required for the $C^0$-bound. 
   As in the constant scalar curvature case, 
   equations pass to the limit and 
   the resulting $\omega_{\varphi_{\infty}}$ is extremal on $\Delta^{*}\times (D_j\backslash E_j)$. 
   Then Theorem \ref{thm_splttngext} tells us that $\varphi_{\infty}$ splits as $at+\psi$, $\psi\in \mathscr{E}(D_j\backslash E_j)$, 
   and $\omega_{D_j}^{\psi}=\omega_{D_j\backslash E_j}+dd^c_{D_j}\psi$ is an extremal metric, of Poincaré type if $E_j\neq\varnothing$, 
   and of class $[\omega_X|_{D_j}]$. 
  \cqfd

  \begin{rmk}    \label{rmk_ExtCvgce}
   A by-product of this proof is: 
   under the same assumptions as in Proposition \ref{prop_extKdivisor}, fix $j\in \{1,\dots,N\}$, and $x_0\in D_j\backslash E_j$. 
   Then for any $(t_k)$ going to $\infty$, 
   any subsequence of $\big((\Pi_{0,j}\varphi)(\cdot+t_k,\cdot)-(\Pi_{0,j}\varphi)(\cdot+t_k,x_0)\big)_k$ converging
   in $C^{\infty}_{\rm loc}\big(\R\times (D_j\backslash E_j)\big)$ has limit of shape $at+\psi$, 
   with $\psi\in \mathscr{E}(D_j\backslash E_j)$, 
   $|a|, |d\psi|_{\omega_{D_j\backslash E_j}}\leq \sup_{\XD} |d\varphi|_{\omega}$, 
   $\omega_{D_j}^{\psi}=\omega_{D_j\backslash E_j}+dd^c_{D_j}\psi$ an extremal Poincaré type metric on $D_j\backslash E_j$, 
   and $1-a\leq c^{-1}$, $\omega_{D_j}^{\psi}\geq c$, where $\omega_{\varphi}\geq c\omega$ on $\XD$. 
   Therefore any subsequence of $t\mapsto \partial_{t_j}\Pi_{0,j}\varphi(t,\cdot)$ converging in $C^{\infty}_{\rm loc}(D_j\backslash E_j)$
   tends to a constant $a$ as above, 
   whereas for any all $\ell\geq 2$, $t\mapsto \partial_{t_j}^{\ell}\Pi_{0,j}\varphi(t,\cdot)$ tends to 0 in $C^{\infty}_{\rm loc}(D_j\backslash E_j)$ 
   as $t\to \infty$. 
  \end{rmk}

 \subsubsection{A numerical constraint}
  Here we use notations of Part \ref{sctn_asympttcs}: for fixed $j\in \{1,\dots,N\}$, 
  the volume $\vl(D_j)$ is computed with respect to the (Poincaré) class induced on $D_j$, 
  as well as the mean scalar curvature $\sbarDj$. 
  We use the content of previous paragraph to prove a constraint on extremal Poincaré type metrics on $\XD$, 
  generalizing the constraint ``$\sbarDj>\sbar$'' of the constant scalar curvature case \cite{auv2}.  
  Notice that in this special case though, the analytical background was less involved; 
  in particular, the analysis of the model case was not required. 
  The extremal constraint states as: 
  \begin{prop}  \label{prop_extcnstrt}
   Assume that $\omega_{\varphi}$ is an extremal Kähler metric of Poincaré type of class $[\omega_X]$, 
   Then for all $j=1,\dots,N$ indexing a component of $D$, one has: 
     \begin{equation}     \label{eqn_extcnstrt}
      \sbarDj+\frac{1}{4\pi\vl(D_j)}\bigg(-\int_{X\backslash D} \mathsf{K}_{\varphi}\cdot e^{t_j} \frac{\omega_{\varphi}^m}{m!}
                                       + \int_{X\backslash D} \scal_{\varphi}\Delta_{\varphi} e^{t_j} \frac{\omega_{\varphi}^m}{m!}
                                       \bigg)>0,
     \end{equation}
    where we set $\scal_{\varphi}=\scal(\omega_{\varphi})$, and recall that $\mathsf{K}_{\varphi}=\nabla^{\varphi}\scal_{\varphi}$.
  \end{prop}
  
  \begin{rmk}
   When $\scal_{\varphi}=\sbar$, one recovers the obstruction $\sbar_{D_j}>\sbar$ in this way: 
   if $\scal_{\varphi}=\sbar$,  $\mathsf{K}_{\varphi}\equiv0$, and 
   $\int_{X\backslash D} \scal_{\varphi}\Delta_{\varphi} e^{t_j} \frac{\omega_{\varphi}^m}{m!}
    = \sbar\int_{X\backslash D}\Delta_{\varphi} e^{t_j} \frac{\omega_{\varphi}^m}{m!}
    = -4\pi\sbar\vl(D_j)$. 
   One can moreover rewrite constraint \eqref{eqn_extcnstrt} as: $\int_{\XD}\lich_{\varphi}(e^{t_j})\frac{\omega_{\varphi}^m}{m!}<0$. 
  \end{rmk}
  
  
  \prf \textit{of Proposition \ref{prop_extcnstrt}. }--- 
   We first make an easy but crucial observation 
   about the statement. 
   If indeed $dd^ce^{t_j}$ is bounded for Poincaré type metrics (see \cite[Part 5]{auv2}), 
   making $\Delta_{\varphi}e^{t_j}$ integrable for the volume form $\omega_{\varphi}^m$, 
   the analogue about $\mathsf{K}_{\varphi}\cdot e^{t_j}$ deserves a slight explanation. 
   Now, as an $L^2$ holomorphic vector field on $X\backslash D$, 
   $\kphi$ extends smoothly through $D$, and its normal component vanishes along the divisor (\cite{auv1}, proof of Lemma 5.2). 
   In other words, in an open set $U$ of local coordinates $(z^1,\dots,z^m)$ centred at some point of $D_j$ such that $D_j$ is given by $z^1=0$, 
   $\kphi$ can be written as $\Re\big(z^1f_1\frac{\partial}{\partial z^1}+f_2\frac{\partial}{\partial z_2}+\cdots+f_m\frac{\partial}{\partial z_m}\big)$, 
   with $f_1, \dots,f_m$ holomorphic on $U$. 
   Moreover in $U$, $e^{t_j}=-\log(|z^1|^2)+f_{j,U}$ with $f_{j,U}$ a smooth function by construction, and therefore, 
   $\kphi\cdot e^{t_j}=-\Re(f_1)+\kphi\cdot f_{j,U}$. 
   From this we get that $\kphi\cdot e^{t_j}$ is bounded near $D_j$, hence on $\XD$, and is thus $L^1$ for $\omega_{\varphi}^m$.  
   
   This being said, let us come to the proof of the proposition itself. 
   By Remark \ref{rmk_ExtCvgce}, 
   on $p_j^{-1}\big(\{t_k\}\times (D_j\backslash E_j)\big)$, one has:
    \begin{equation}    \label{eqn_aspttcsomegaphi}
     \omega_{\varphi}= a_kdt_j\wedge 2e^{-t_j}\eta_j+p_j^*\omega_{t_k}^{\varphi}+\vareps(t_k,z), 
    \end{equation}
   where: $a_k\geq c$ ($c$ such that $\omega_{\varphi}\geq 2c\omega$),  
   the $\omega_{t_k}^{\varphi}:=\omega_{D_j}+dd^c_{D_j}[\Pi_{0,j}\varphi(t_k,\cdot)]$ 
   are uniformly bounded on $(D_j\backslash E_j)$ at any order in the Poincaré sense and 
   positively bounded below on compact subsets, 
   and all the $\partial_{t_j}^{\ell}\vareps$, $\ell\geq0$, 
   tend to 0 in $C^{\infty}_{\rm loc}(D_j\backslash E_j)$ 
   as $k$ goes to $\infty$. 
   Consequently, a direct computation yields, on $p_j^{-1}\big(\{t_k\}\times (D_j\backslash E_j)\big)$:  
    \begin{equation}   \label{eqn_aspttcsScalomegaphi}
     \scal_{\varphi} = -2a_kdt_j\wedge 2e^{-t_j}\eta_j+p_j^*\varrho(\omega_{t_k}^{\varphi})+\vareps_1(t_k,z), 
    \end{equation}
   with $\vareps_1$ tending to 0 in the same sense as $\vareps$ above.  
   We can now write, for $s$ large and seeing $\{t_j\leq s\}$ as a subset of $\XD$: 
    \begin{align*}
     \int_{\{t_j\leq s\}} \scal_{\varphi}&dd^c e^{t_j}\wedge\frac{\omega_{\varphi}^{m-1}}{(m-1)!}               \\
      = &\int_{\{t_j\leq s\}} d\Big(\scal_{\varphi}d^c e^{t_j}\wedge\frac{\omega_{\varphi}^{m-1}}{(m-1)!} \Big)
        - \int_{\{t_j\leq s\}} d\scal_{\varphi}\wedge d^c e^{t_j}\wedge\frac{\omega_{\varphi}^{m-1}}{(m-1)!}    \\
      = &\int_{\{t_j = s\}} \scal_{\varphi}d^c e^{t_j}\wedge\frac{\omega_{\varphi}^{m-1}}{(m-1)!} 
        - \int_{\{t_j\leq s\}} \kphi\cdot e^{t_j}\frac{\omega_{\varphi}^{m}}{m!},   
    \end{align*}
   where we used Stokes' theorem to pass from the integral of shape $\int_{\{t_j\leq s\}} d\alpha$ in the second line 
   to the integral of shape $\int_{\{t_j = s\}}\alpha$ in the third line 
   (there is no interference with the other $D_{\ell}$ here, 
   as the $\alpha$ in play, and its differential, are $L^1$), 
   and applied a hermitian identity to $\kphi=\nabla^{\varphi}\scal_{\varphi}$ to replace $d\scal_{\varphi}\wedge d^c e^{t_j}\wedge\frac{\omega_{\varphi}^{m-1}}{(m-1)!}$ 
   by $\kphi\cdot e^{t_j}\frac{\omega_{\varphi}^{m}}{m!}$ in the last summand. 
   Now 
   $\int_{\{t_j\leq s\}} \scal_{\varphi}dd^c e^{t_j}\wedge\frac{\omega_{\varphi}^{m-1}}{(m-1)!}
    =-\int_{\{t_j\leq s\}} \scal_{\varphi}\Delta_{\varphi} e^{t_j}\vol^{\varphi}$ 
   tends to $-\int_{\XD} \scal_{\varphi}\Delta_{\varphi} e^{t_j}\vol^{\varphi}$ as $s$ goes to $\infty$; 
   likewise, $\int_{\{t_j\leq s\}} \kphi\cdot e^{t_j}\frac{\omega_{\varphi}^{m}}{m!}$ tends to $\int_{\XD} \kphi\cdot e^{t_j}\vol^{\varphi}$ as $s$ goes to $\infty$. 
   Finally, 
   taking $s=t_k$ and using the asymptotics \eqref{eqn_aspttcsomegaphi} and \eqref{eqn_aspttcsScalomegaphi} with the evoked uniform bounds, 
   we get that $\int_{\{t_j = t_k\}} \scal_{\varphi}d^c e^{t_j}\wedge\frac{\omega_{\varphi}^{m-1}}{(m-1)!}
                 = 4\pi\sbar_{D_j}\vl(D_j)-4\pi \int_{D_j\backslash E_j}a_k^{-1}\vol^{\omega_{t_k}^{\varphi}}+o(1)$ as $k$ goes to $\infty$, 
   and the last integral is $\geq c^{-1}\vl(D_j)$, hence the result. 
   \cqfd
   
  \begin{rmk}  \label{rmk_finala}
   We also get from this proof that $(a_k)_k$ tends to the inverse of the left-hand-side of \eqref{eqn_extcnstrt}, $a_j$ say, 
   \emph{that depends only on $\omega_{\varphi}$}, this for all sequence $(t_k)$. 
   We hence sharpen Remark \ref{rmk_ExtCvgce} by saying that: \emph{$t\mapsto\partial_{t_j}\Pi_{0,j}\varphi(t,\cdot)$ 
   tends to $(1-a_j)$ in $C^{\infty}_{\loc}(D_j\backslash E_j)$}.  
  \end{rmk}

 \subsubsection{Asymptotics of extremal Kähler metrics of Poincaré type (smooth divisor)}
  When $D$ is smooth, we have a perfect analogue of Theorem \ref{thm_asmpttcsPKcsc}: 
   \begin{thm} \label{thm_asmpttcsPKext}
    Assuming $D$ smooth, let $\omega_{\varphi}$ be an extremal metric of Poincaré type on $\XD$, of class $[\omega_X]$. 
    One has near $D_j$ the asymptotics of Riemannian metrics: 
     \begin{equation*}   \label{eqn_asympExtOmega}
      g_{\varphi} = a_j\big(dt_j^2 + (2e^{-t_j}\eta)^2\big) + p^*g_{D_j}   +\mathcal{O}(e^{-\delta t_j})
     \end{equation*}
    with $g_{D_j}$ an extremal Kähler metric on $(D_j,J_{D_j})$, such that $[g_{D_j}(J_{D_j}\cdot,\cdot)]=[\omega|_{D_j}]$, 
    and where $a_j>0$ is the inverse of the left-hand side of \eqref{eqn_extcnstrt}; 
    this holds for all $j$. 
   \end{thm}

  \prf. --- 
   The ingredients are the same as in the constant scalar curvature case. 
   Namely, \textit{1.} starting from an equation analogous to \eqref{eqn_lichnearD}, 
   \textit{2.} we twist it after parametrising judiciously the ``family of almost extremal metrics'' $(\omega_t^{\varphi})_t$, 
   and then \textit{3.} lead the appropriate asymptotic analysis of the arising Lichnerowicz-modelled operator. 
   Before examining those three points, 
   and more precisely why they adapt to the extremal situation, 
   we assume that the divisor is reduced to one component, 
   and we proceed to the following 
   normalisation: 
   call $a_D=a_D(\omega_{\varphi})$ the inverse of the left hand-side of \eqref{eqn_extcnstrt}; 
   up to working with $a_D^{-1}\omega_{\varphi}$,
   which near $D$ can be written as $a_D^{-1}\omega-dd^ct+dd^c\big(a_D^{-1}\varphi-(a_D^{-1}-1)t\big)$ -- 
   recall that $\omega_{\varphi}=\omega+dd^c(\varphi-t)$ --,  
   we can assume that $a_D=1$, 
   as $a_D(\lambda\omega_{\varphi})=\lambda a_D(\omega_{\varphi})$ for all $\lambda>0$, by \eqref{eqn_extcnstrt}.  
   We then recover, by Remarks \ref{rmk_ExtCvgce} and \ref{rmk_finala}, 
   the extremal analogue of Proposition \ref{prop_almcnstt}.

   The first point goes as follows.
   Again we work on an open subset $U$ of coordinates $(z^1,\dots,z^m)$ centred at some point of $D$, 
   such that $z^1$ is the local equation of $D$, 
   and use the holomorphic local vector field 
   $Z=\Re\big(z^1(\log z^1)\frac{\partial}{\partial z^1}\big)$ to differentiate $\scal_{\varphi}$. 
   On the one hand, 
   $Z\cdot \scal_{\varphi}=d\scal_{\varphi}(Z)$;  
   on the other hand, 
    \begin{align*}
     Z\cdot \scal_{\varphi}    =&-2\Big(\frac{1}{2}\Delta_{\varphi}^2\big(Z\cdot(\varphi-t)\big)
                                       +\big\langle\varrho_{\varphi},dd^c[Z\cdot(\varphi-t)]\big\rangle_{\varphi}\Big)                         \\
                                &+\Delta_{\varphi} (\Lambda_{\varphi}\mathcal{L}_Z\omega_X )
                                      - 2  \langle\varrho_{\varphi},\mathcal{L}_Z\omega_X \rangle_{\varphi}         +\mathcal{O}(e^{-t})             \\
                               =&-2\lich_{\varphi}[Z\cdot(\varphi-t)]+\big\langle d\scal_{\varphi}, d(Z\cdot(\varphi-t))\big\rangle_{\varphi} +\mathcal{O}(e^{-t}),
    \end{align*}
   as $\lich_{\varphi}=\frac{1}{2}\Delta_{\varphi}^2+\langle\varrho_{\varphi},dd^c\cdot\rangle_{\varphi}
                       +\frac{1}{2}\big\langle d\scal_{\varphi}, d\cdot\big\rangle_{\varphi}$, 
   and $\mathcal{L}_Z\omega_X=\mathcal{O}(e^{-t})$ at any order. 
   Now $d\scal_{\varphi}(Z)
        =d\scal_{\varphi}\big(\frac{1}{2}\rho_1(\log\rho_1)\frac{\partial}{\partial\rho_1}\big)+\mathcal{O}(e^{-t})$ 
   with $z^1=\rho_1e^{i\theta_1}$, 
   and thus $Z=\frac{1}{2}\big(\rho_1(\log\rho_1)\frac{\partial}{\partial\rho_1}+\theta_1\frac{\partial}{\partial\theta_1}\big)$. 
   Set $\mathsf{K}_{D}=\kphi|_D$. 
   Then $\mathsf{K}_{\varphi}=\mathsf{K}_D+\mathcal{O}(e^{-t})$ in $U$, 
   and therefore 
   $d\scal_{\varphi}(Z)
    =\frac{1}{2}\big\langle\mathsf{K}_{\varphi},\rho_1(\log\rho_1)\frac{\partial}{\partial\rho_1}\big\rangle_{\varphi}+\mathcal{O}(e^{-t})
    =\frac{1}{2}d_D\partial_t(\Pi_0\varphi)(\mathsf{K}_D)+\mathcal{O}(e^{-t})$, as  
   $dt\big(\rho_1(\log\rho_1)\frac{\partial}{\partial\rho_1}\big)=1+\mathcal{O}(e^{-t})$ and $\eta\big(\rho_1(\log\rho_1)\frac{\partial}{\partial\rho_1}\big)=\mathcal{O}(e^{-t})$,
   whereas $d\big(Z\cdot(\varphi-t)\big)=\frac{1}{2}d\partial_t(\Pi_0\varphi)+\mathcal{O}(e^{-t})$, 
   hence:  
    \begin{equation*}
     \big\langle d\scal_{\varphi}, d\big(Z\cdot(\varphi-t)\big)\big\rangle_{\varphi}
                           = \frac{1}{2}\mathsf{K}_{\varphi}\cdot \partial_t(\Pi_0\varphi) +\mathcal{O}(e^{-t})
                           = \frac{1}{2}d_D\partial_t(\Pi_0\varphi)(\mathsf{K}_D)+\mathcal{O}(e^{-t})                                           
    \end{equation*}
   Conclude that $\lich_{\varphi}[Z\cdot(\varphi-t)]=\mathcal{O}(e^{-t})$;  
   as above, rewrite the latter as $\lich_{\varphi}[\partial_t(\Pi_0\varphi)]=\mathcal{O}(e^{-t})$, 
   which makes sense globally near $D$. 
   
   For second point, we mostly have to see that the extremal condition can replace 
   the constant scalar curvature condition in the construction of Part \ref{sctn_almcscKfamily}. 
   Observe that:
    \begin{itemize}
     \item[$\diamond$] 
           on a compact manifold, 
           the extremal condition is a non-linear 6th order elliptic equation on the Kähler potential, with $C^{\infty}$-bounded coefficients 
           as soon as the metrics are bounded in $C^{\infty}$, and bounded below; 
     \item[$\diamond$] Calabi's theorem is stated on extremal Kähler metrics. 
    \end{itemize}
   Consequently, 
   one recovers by the same methods the exact analogues of the results of Part \ref{sctn_almcscKfamily}, with the vectors $Z$ 
   of this part in the $J_Y\mathfrak{k}_0(\omega)\oplus\big(\bigoplus_{\lambda>0}\mathfrak{h}^{(\lambda)}\big)$ for $\omega$ extremal, 
   where $Z_1$ holomorphic is in $\mathfrak{h}^{(\lambda)}$ iff $\mathcal{L}_{J_Y\nabla^{\omega}\scal(\omega)}Z_1=\lambda J_YZ_1$. 
   
   In passing, we fix $\omega_D=g_D(J_D\cdot,\cdot)$ a limit point of $(\omega_t^{\varphi})$, 
   which is extremal by Remark \ref{rmk_ExtCvgce}; 
   this way, $\nabla^{g_D}\scal(\omega_D)=\mathsf{K}_D:=\kphi|_D$. 
   We assume also that $\omega|_D=\omega_D$. 
   
   Now for the third and last point, 
   the splitting formula 
   \eqref{eqn_asymplich2} is still valid 
   (since $\omega_{\varphi}=dt\wedge2e^{-t}\eta+p^*\omega_t^{\varphi}+o(1)$, 
    $d\scal_{\varphi}=p^*d\scal(\omega_t^{\varphi})+o(1)$, which is parallel to $D$ up to an $o(1)$). 
   Lemmas \ref{lem_mapprop} (dealing with the model operator) and 
   \ref{lem_fredholm} (on the conjugated geometric operator) thus remain valid. 
   To recover the full analogue of Proposition \ref{prop_isom}, 
   we need to identify the non-asymptotically small functions whose image by $\lich_{\varphi}^{\tilde{\sigma}}$ are exponentially small; 
   in other words, we need the analogue of equations \eqref{eqn_lmpsij}. 
   We can actually 
   use the same candidates as those of Lemma \ref{lem_psi_j},  
   with the same notations as in this lemma;   
   similar computations yield: 
    \begin{equation*}
     Z_j\cdot\scal_{\varphi}=-2\lich_{\varphi}[f_j+\Pi_0(Z_j\cdot\varphi)]+\big\langle d(f_j+\Pi_0(Z_j\cdot\varphi)),d\scal_{\varphi}\big\rangle_{\varphi}
                                     +\mathcal{O}(e^{-t}); 
    \end{equation*}
   if we prove that $d\scal_{\varphi}(Z_j)=\big\langle d(f_j+\Pi_0(Z_j\cdot\varphi)),d\scal_{\varphi}\big\rangle_{\varphi}+\mathcal{O}(e^{-t})$, 
   which is somehow the most delicate point of our argument, 
   we will thus be done in the same way as in the constant scalar curvature case. 
   
   Fix in $X$ an open neighbourhood $U$ of holomorphic coordinates of any point in $D$,  
   and extend vector fields on $U\cap D$, such as $Z_j|_D$, to $U$, in the natural way. 
   Set also $\alpha\sim\beta$ if $\alpha=\beta+\mathcal{O}(e^{-t})$ at any order for $\omega$ in $U$. 
   
   Observe now that on the one hand, $d\scal_{\varphi}(Z_j) = \langle\kphi,Z_j\rangle_{\varphi} \sim \langle\mathsf{K}_D,Z_j|_D\rangle_{g_t^{\varphi}}$, 
   as $\kphi\sim\mathsf{K}_D$, $Z_j\sim Z_j|_D $ and 
   $g_{\varphi}=\big(1+(\partial_t^2-\partial_t)\Pi_0\varphi)(dt^2+4e^{-2t}\eta^2)+dt\cdot d_D\Pi_0\varphi
                +2e^{-t}\eta\cdot d_D^c\Pi_0\varphi+p^*g_t^{\varphi}+\mathcal{O}(e^{-t})$, 
   with $g_t^{\varphi}=\omega_t^{\varphi}(\cdot,J_D\cdot)$.  
   On the other hand, 
   $\big\langle d(f_j+\Pi_0(Z_j\cdot\varphi)),d\scal(\omega_{\varphi})\big\rangle_{\varphi}
     = \kphi\cdot(f_j+\Pi_0(Z_j\cdot\varphi))
     \sim \mathsf{K}_D\cdot f_j + \kphi\cdot\Pi_0(Z_j\cdot\varphi) 
     \sim \langle\mathsf{K}_D,Z_j|_D\rangle_{\omega_D}+ \kphi\cdot\Pi_0(Z_j\cdot\varphi)$. 
   Comparing those two expressions, 
   our next task is to show that $\kphi\cdot\Pi_0(Z_j\cdot\varphi)\sim dd^c_D(\Pi_0\varphi)(\mathsf{K}_D,J_DZ_j|_D)$. 
   
   By definition, 
   $dd^c_D(\Pi_0\varphi)(\mathsf{K}_D,J_DZ_j|_D)
    =\mathsf{K}_D\cdot[d^c_D\Pi_0\varphi(J_DZ_j|_D)]
     - (J_DZ_j|_D)\cdot[d^c_D\Pi_0\varphi(\mathsf{K}_D)]
     - d^c_D\Pi_0\varphi([\mathsf{K}_D,J_DZ_j|_D])$ on $D$. 
   But as $\kphi\sim \mathsf{K}_{\varphi}$ and $Z_j\sim Z_j|_D$, 
   $\mathsf{K}_D\cdot[d^c_D\Pi_0\varphi(J_DZ_j|_D)]=\mathsf{K}_D\cdot[d_D\Pi_0\varphi(Z_j|_D)]\sim\kphi\cdot(d\Pi_0\varphi(Z_j))$;  
   moreover, as $\mathsf{K}_D=\nabla^{g_D}\scal(\omega_D)$ and $Z_j|_D\in J_D\mathfrak{k}_0(\omega_D)$, $[\mathsf{K}_D, J_DZ_j|_D]=0$. 
   Finally, $d^c_D\Pi_0\varphi(\mathsf{K}_D)\sim -(J\kphi)\cdot\varphi$, 
   and on $\XD$, $\mathcal{L}_{J\kphi}\omega_{\varphi}=0$, 
   whereas near $D$, $\omega_{\varphi}\sim p^*\omega_D+dd^c(\varphi-t)$, 
   $J\kphi \sim J_D\mathsf{K}_D$ and $\mathcal{L}_{J_D\mathsf{K}_D}\omega_D=0$, 
   so that $\mathcal{L}_{J\kphi}\omega_{\varphi}\sim dd^c\big((J\kphi)\cdot(\varphi-t)\big)$, 
   and $dd^c\big((J\kphi)\cdot\varphi\big)\sim dd^c\big((J\kphi)\cdot(\varphi-t)\big)\sim 0$. 
   According to the weighted $\ddbar$-lemma of \cite{auv1}
   (or more precisely to the proof of \cite[Lemma 3.10]{auv1}, which is a local version near the divisor), 
   this implies that $(J\kphi)\cdot\varphi \sim c $ for some constant $c$, 
   so that $(JZ_j)\cdot[(J\kphi)\cdot\varphi]\sim 0$, hence $(J_DZ_j|_D)\cdot[d^c_D\Pi_0\varphi(\mathsf{K}_D)]\sim0$. 
  \cqfd
   
\vspace{1cm}

 \small

 \bibliographystyle{amsalpha}
 \bibliography{biblioExtPK}

\providecommand{\bysame}{\leavevmode\hbox to3em{\hrulefill}\thinspace}
\providecommand{\MR}{\relax\ifhmode\unskip\space\fi MR }
\providecommand{\MRhref}[2]{%
  \href{http://www.ams.org/mathscinet-getitem?mr=#1}{#2}
}
\providecommand{\href}[2]{#2}
\begin{thebibliography}{CDS12b}

\bibitem[AH12]{ah}
Vestislav Apostolov and Hongnian Huang, \emph{A splitting theorem for extremal
  kaehler metrics}, Preprint arXiv:1212.3665[math.DG], 2012.

\bibitem[AP06]{ap}
Claudio Arezzo and Frank Pacard, \emph{Blowing up and desingularizing constant
  scalar curvature {K}\"ahler manifolds}, Acta Math. \textbf{196} (2006),
  no.~2, 179--228.

\bibitem[APS11]{aps}
Claudio Arezzo, Frank Pacard, and Michael Singer, \emph{Extremal metrics on
  blowups}, Duke Math. J. \textbf{157} (2011), no.~1, 1--51. \MR{2783927
  (2012k:32024)}

\bibitem[Auv11]{auv1}
Hugues Auvray, \emph{The space of {P}oincar\'e type {K}\"ahler metrics on the
  complement of a divisor}, Preprint arXiv:1109.3159 [math.DG], 2011.

\bibitem[Auv13]{auv2}
\bysame, \emph{Metrics of {P}oincar\'e type with constant scalar curvature: a
  topological constraint}, J. Lond. Math. Soc. (2) \textbf{87} (2013), no.~2,
  607--621.

\bibitem[BDB88]{bdb}
D.~Burns and P.~De~Bartolomeis, \emph{Stability of vector bundles and extremal
  metrics}, Invent. Math. \textbf{92} (1988), no.~2, 403--407.

\bibitem[Biq97]{biq}
Olivier Biquard, \emph{Fibr\'es de {H}iggs et connexions int\'egrables: le cas
  logarithmique (diviseur lisse)}, Ann. Sci. \'Ecole Norm. Sup. (4) \textbf{30}
  (1997), no.~1, 41--96.

\bibitem[Cal82]{cal}
Eugenio Calabi, \emph{Extremal {K}\"ahler metrics}, Seminar on {D}ifferential
  {G}eometry, Ann. of Math. Stud., vol. 102, Princeton Univ. Press, Princeton,
  N.J., 1982, pp.~259--290.

\bibitem[Cal85]{cal2}
\bysame, \emph{Extremal {K}\"ahler metrics. {II}}, Differential geometry and
  complex analysis, Springer, Berlin, 1985, pp.~95--114.

\bibitem[CDS12a]{cds1}
X.~X. Chen, S.~K. Donaldson, and Song Sun, \emph{Kahler-{E}instein metrics on
  {F}ano manifolds, {I}: approximation of metrics with cone singularities},
  Preprint arXiv:1211.4566[math.DG], 2012.

\bibitem[CDS12b]{cds2}
\bysame, \emph{Kahler-{E}instein metrics on {F}ano manifolds, {II}: limits with
  cone angle less than $2\pi$}, Preprint arXiv:1212.4714[math.DG], 2012.

\bibitem[CDS13]{cds3}
\bysame, \emph{Kahler-{E}instein metrics on {F}ano manifolds, {III}: limits as
  cone angle approaches $2\pi$ and completion of the main proof}, Preprint
  arXiv:1302.0282[math.DG], 2013.

\bibitem[CT08]{ChTi}
X.~X. Chen and G.~Tian, \emph{Geometry of {K}\"ahler metrics and foliations by
  holomorphic discs}, Publ. Math. Inst. Hautes \'Etudes Sci. (2008), no.~107,
  1--107. \MR{2434691 (2009g:32048)}

\bibitem[Don01]{don}
S.~K. Donaldson, \emph{Scalar curvature and projective embeddings. {I}}, J.
  Differential Geom. \textbf{59} (2001), no.~3, 479--522.

\bibitem[Gau]{gau}
Paul Gauduchon, \emph{Calabi's extremal metrics: an elementary introduction},
  Lecture notes.

\bibitem[Hua12]{hua}
Hongnian Huang, \emph{A splitting theorem on toric varieties}, Preprint
  arXiv:1212.3729 [math.DG], 2012.

\bibitem[Lev85]{lev}
Marc Levine, \emph{A remark on extremal {K}\"ahler metrics}, J. Differential
  Geom. \textbf{21} (1985), no.~1, 73--77.

\bibitem[LM85]{lm}
Robert~B. Lockhart and Robert~C. McOwen, \emph{Elliptic differential operators
  on noncompact manifolds}, Ann. Scuola Norm. Sup. Pisa Cl. Sci. (4)
  \textbf{12} (1985), no.~3, 409--447.

\bibitem[Mab04]{mab}
Toshiki Mabuchi, \emph{Stability of extremal {K}\"ahler manifolds}, Osaka J.
  Math. \textbf{41} (2004), no.~3, 563--582.

\bibitem[Sch02]{sch}
Georg Schumacher, \emph{Asymptotics of complete {K}\"ahler-{E}instein
  metrics---negativity of the holomorphic sectional curvature}, Doc. Math.
  \textbf{7} (2002), 653--658 (electronic).

\bibitem[Sz{\'e}06]{sze1}
G{\'a}bor Sz{\'e}kelyhidi, \emph{Extremal metrics and {K}-stability}, Ph.D.
  thesis, 2006.

\bibitem[Sz{\'e}07]{sze2}
\bysame, \emph{Extremal metrics and {$K$}-stability}, Bull. Lond. Math. Soc.
  \textbf{39} (2007), no.~1, 76--84.

\bibitem[Tia97]{tia}
Gang Tian, \emph{K\"ahler-{E}instein metrics with positive scalar curvature},
  Invent. Math. \textbf{130} (1997), no.~1, 1--37.

\bibitem[Tia12]{tia2}
\bysame, \emph{K-stability and {K}\"ahler-{E}instein metrics}, Preprint
  arXiv:1211.4669 [math.DG], 2012.

\bibitem[TY90]{ty}
Gang Tian and Shing-Tung Yau, \emph{Complete {K}\"ahler manifolds with zero
  {R}icci curvature. {I}}, J. Amer. Math. Soc. \textbf{3} (1990), no.~3,
  579--609.

\bibitem[Wu06]{wu}
Damin Wu, \emph{Higher canonical asymptotics of {K}\"ahler-{E}instein metrics
  on quasi-projective manifolds}, Comm. Anal. Geom. \textbf{14} (2006), no.~4,
  795--845.

\bibitem[Yau93]{yau}
Shing-Tung Yau, \emph{Open problems in geometry}, Differential geometry:
  partial differential equations on manifolds ({L}os {A}ngeles, {CA}, 1990),
  Proc. Sympos. Pure Math., vol.~54, Amer. Math. Soc., Providence, RI, 1993,
  pp.~1--28.

\end{thebibliography}

~

\small \textsc{CMLA, École Normale Supérieure de Cachan, UMR 8536} \\
\indent 61 avenue du Président Wilson, 94230 Cachan, France        \\
\indent \url{hugues.auvray@ens-cachan.fr}

 \end{document}